\documentclass[10pt]{article}
\usepackage{amsmath,amssymb,amsthm,comment}
\usepackage{mathtools}
\usepackage{color}
\usepackage[nohead,margin=1.1in]{geometry}

\usepackage{booktabs}
\usepackage{threeparttable}
\usepackage{graphicx,psfrag,epstopdf}
\usepackage{algorithmic,algorithm}
\usepackage{url}
\graphicspath{{./pics/}}
\usepackage{tcolorbox}
\allowdisplaybreaks
\theoremstyle{plain}
\newtheorem{theorem}{Theorem}[section]
\newtheorem{remark}{Remark}[section]

\newtheorem{lemma}{Lemma}[section]
\newtheorem{assumption}{Assumption}[section]
\newtheorem{prop}{Proposition}[section]

\newtheorem{corollary}{Corollary}[section]
\numberwithin{equation}{section}
\def\E{\mathbb{E}}

\begin{document}
\title{On the Convergence of A Data-Driven Regularized  
\\ Stochastic Gradient Descent for Nonlinear Ill-Posed Problems}

\author{Zehui Zhou\thanks{Department of Mathematics,
Rutgers University, 110 Frelinghuysen Road, Piscataway, NJ 08854-8019 (\texttt{zz569@math.rutgers.edu})}}

\date{}
\maketitle

\begin{abstract}

Stochastic gradient descent (SGD) is a promising method for solving large-scale inverse problems, due to its excellent scalability with respect to data size. 
In this work, we analyze a new data-driven regularized stochastic gradient descent for the efficient numerical solution of a class of nonlinear ill-posed inverse problems in infinite dimensional Hilbert spaces. 
At each step of the iteration, the method randomly selects one equation from the nonlinear system combined with a corresponding equation from the learned system based on training data to obtain a stochastic estimate of the gradient and then performs a descent step with the estimated gradient. 
We prove the regularizing property of this method under the tangential cone condition and \textit{a priori} parameter choice and then derive the convergence rates under the additional source condition and range invariance conditions.
Several numerical experiments are provided to complement the analysis.


\end{abstract}

\section{Introduction}

This work is concerned with the numerical solution of the system of nonlinear ill-posed operator equations
\begin{equation}\label{eqn:nonlin}
  F_i(x) = y_i^\dag,\quad i =1,\ldots, n,
\end{equation}
where each $F_i: \mathcal{D}(F_i)\to Y$ is a nonlinear mapping
with its domain $\mathcal{D}(F_i)\subset X$, and $X$ and $Y$ are Hilbert spaces with inner
products $\langle\cdot,\cdot\rangle$ and norms $\|\cdot\|$ respectively.
The number 
$n$ of nonlinear equations in \eqref{eqn:nonlin} can potentially be very large.
The notation $y_i^\dag \in Y$ denotes the exact data (corresponding to the reference solution
$x^\dag\in X$ to be defined below, i.e., $y^\dag=F(x^\dag)$). Equivalently, problem \eqref{eqn:nonlin} can be rewritten as
\begin{equation}\label{eqn:nonlin-gen}
  F(x) = y^\dag,
\end{equation}
with the operator $F: \bigcap_{i=1}^n \mathcal{D}(F_i)\subset X \to Y^n$ ($Y^n$ denotes the product space $Y\times \cdots \times Y$) and $y^\dag\in Y^n$ defined by
\begin{equation*}
  F(x)=\frac{1}{\sqrt{n}}\left(\begin{array}{c}F_1(x)\\ \ldots\\ F_n(x)\end{array}\right)\quad \mbox{and}\quad y^\dag=\frac{1}{\sqrt{n}}\left(\begin{array}{c}y_1^\dag\\ \ldots\\ y_n^\dag\end{array}\right),
\end{equation*}
respectively. The scaling $n^{-\frac12}$ above is introduced for the convenience of later discussions. In practice,
instead of the exact data $y^\dag$, we have access only to {the noisy data $y^\delta=y^\dag+\xi$ with the noise $\xi=\frac{1}{\sqrt{n}}\left(\begin{array}{c}\xi_1\\ \ldots\\ \xi_n\end{array}\right)$ of a noise level
$\delta\ge 0$, namely
\begin{equation*}
  \|\xi\|=\|y^\delta - y^\dag\| \leq \delta\,.
\end{equation*}}
Nonlinear inverse problems of the form \eqref{eqn:nonlin} arise in a broad range of practical applications, e.g., inverse scattering and electrical impedance tomography. 
Stochastic Gradient Descent (SGD), first proposed by Robbins and Monro \cite{RobbinsMonro:1951}, which is a randomized version of the classical Landweber method \cite{Landweber:1951}, is a very popular stochastic iterative method for solving nonlinear ill-posed inverse problems \cite{HankeNeubauerScherzer:1995,KaltenbacherNeubauerScherzer:2008,ItoJin:2015,JinZhouZou:2020} and has also attracted strong interest in machine learning \cite{SutskeverMartens:2013,BottouCurtisNocedal:2018}, due to its excellent scalability with respect to the truly massive data set (i.e., large $n$).
However, analyzing SGD from the perspective of regularization theory to solve ill-posed inverse problems remains largely under-explored, despite their computational appeals. 
The theoretical analysis of SGD-type algorithms for ill-posed inverse problems has started only recently.
Existing works on linear and nonlinear inverse problems \cite{JinLu:2019,JinZhouZou:2020,JinZhouZou:2021} focus on the standard SGD combined with \textit{a priori} stopping rules, which has been proved to be a regularized numerical method, meanwhile several works discuss different variants of SGD with various acceleration strategies \cite{LeRouxSchmidt:2012,DefazioBachLacoste:2014,KingmaBa:2015,NguyenLiu:2017,JinZhouZou:2022}. 
Few of these works use \textit{a priori} training data for the inverse problem in the lens of regularization theory. 
However, the lack of \textit{a priori} knowledge of the true solution may pose some challenges to SGD, e.g., without suitable assumptions on the true solution, the iterations may converge to a solution far away from the exact solution due to its high sensitivity to initial guess and may lead to overfitting due to its ability to quickly adapt to the noisy data. 

In this work, we are interested in the convergence analysis of a variant of SGD for problem \eqref{eqn:nonlin} given in Algorithm \ref{alg:datasgd} which incorporates \textit{a priori} knowledge for the problem. 
In the algorithm, the index $i_k$ of the equation  at the $k$th iteration is drawn uniformly from the index set $\{1,\ldots,n\}$, $\eta_k>0$ is the step size, and $\lambda_k^\delta>0$ is the regularization parameter. The data-driven model $G: X \to Y^n $ in the regularization term, given by  
\begin{equation*}
  G(x)=\frac{1}{\sqrt{n}}\left(\begin{array}{c}G_1(x)\\ \ldots\\ G_n(x)\end{array}\right),
\end{equation*}
is learned by the prior information of the problem, i.e., a set of data pairs $\{x^{(i)},y^{(i)}\}_{i=1}^N$, using neural networks.

\begin{algorithm}
  \caption{Data-driven regularized stochastic gradient descent method for problem \eqref{eqn:nonlin}.\label{alg:datasgd}}
  \begin{algorithmic}[1]
    \STATE Given initial guess $x_1$.
    \FOR{$k=1,2,\ldots$}
      \STATE Randomly draw an index $i_k$;
      \STATE Update the iterate $x_k^\delta$ by
      \begin{equation}\label{eqn:datasgd}
        x_{k+1}^\delta = x_k^\delta - \eta_k \big(F_{i_k}'(x_k^\delta)^*(F_{i_k}(x_k^\delta)-y_{i_k}^\delta)+\lambda_k^\delta G_{i_k}'(x_k^\delta)^*(G_{i_k}(x_k^\delta)-y_{i_k}^\delta)\big);
      \end{equation}
      \STATE Check the stopping criterion.
    \ENDFOR
  \end{algorithmic}
\end{algorithm}

Algorithmically, this data-driven regularized stochastic gradient descent method can be viewed as a randomized version of the data-driven iteratively regularized Landweber method proposed in \cite{AspriBanert:2020}, which is given by
\begin{equation}\label{eqn:dataLandweber}
  x_{k+1}^\delta = x_k^\delta - \eta_k \big(F'(x_k^\delta)^*(F(x_k^\delta)-y^\delta)+\lambda_k^\delta G'(x_k^\delta)^*(G(x_k^\delta)-y^\delta) \big)
\end{equation}
with the step size $\eta_k\equiv 1$.
The $k$-th step of \eqref{eqn:dataLandweber} can be viewed as the gradient descent applied to the following functional
\begin{equation*}
  J(x) = \frac{1}{2}\big(\|F(x)-y^\delta\|^2 +\lambda_k^\delta\|G(x)-y^\delta\|^2\big) = \frac{1}{n}\sum_{i=1}^n \frac{1}{2} \big(\|F_i(x)-y_i^\delta\|^2+\lambda_k^\delta\|G_i(x)-y_i^\delta\|^2\big).
\end{equation*}
Compared with the corresponding Landweber method \eqref{eqn:dataLandweber}, the data-driven regularized SGD \eqref{eqn:datasgd} employs only one randomly selected equation from the true model and data-driven model at each iteration to obtain the gradient estimate. 
Thus, the computational complexity is independent of the data size (which can be very large), which indicates excellent scalability with respect to the problem scale.

For linear and nonlinear ill-posed inverse problems, there exists a relatively thorough understanding of the Landweber method \cite{Louis:1989,VainikkoVeretennikov:1986,HankeNeubauerScherzer:1995,KaltenbacherNeubauerScherzer:2008}, including various data-driven regularized Landweber methods \cite{Scherzer:1998,AspriBanert:2020}. 
It has been shown that a regularization term (based on an initial guess from the data) in \cite{Scherzer:1998} stabilizes the algorithm by enabling the iterations to converge to a solution closest to the initial guess without making additional assumptions about the true solution (which is necessary for the standard Landweber method \cite{HankeNeubauerScherzer:1995}), however, it provides even slower practical convergence rates than that for the standard Landweber method. 
Motivated by this observation, a damping factor {based on a data-driven model} is introduced into the standard Landweber iteration in \cite{AspriBanert:2020}, where a strong convergence and stability for the algorithm are presented.
Intuitively, as a randomized version of the data-driven iteratively regularized Landweber method, data-driven regularized SGD defined in Algorithm \ref{alg:datasgd} is expected to enjoy similar desirable properties.

In this work, we contribute to the convergence analysis of the data-driven regularized SGD defined in Algorithm \ref{alg:datasgd} for a class of nonlinear inverse problems of the form \eqref{eqn:nonlin} from the perspective of regularization theory. 
Under the classical tangential cone condition, we prove the regularizing property of this algorithm when combined with \textit{a priori} rules on the parameter choice; see {Theorems \ref{thm:conv-exact} (for exact data) and \ref{thm:conv-noisy} (for noisy data)}. 
Further, under suitable source condition, range invariance condition and its stochastic variant, we achieve the convergence rates of this algorithm with polynomially decaying step size and regularization parameter schedules, which are comparable with that for the Landweber method in \cite{HankeNeubauerScherzer:1995} and the standard SGD for both linear and nonlinear cases in \cite{JinLu:2019} and \cite{JinZhouZou:2020}; see {Theorems \ref{thm:err-total-ex} (for exact data) and \ref{thm:err-total} (for noisy data)}.
The analysis draws on strategies for handling the data-driven damping term of the data-driven regularized Landweber method in \cite{AspriBanert:2020} and estimating the general error of the standard SGD in \cite{JinZhouZou:2020}.

Throughout, we denote the $(k-1)$-th iterates for the exact data $y^\dag$ and the noisy data $y^\delta$ by $x_k$ and $x_k^\delta$ respectively. Let $x^*$ be any solution to problem \eqref{eqn:nonlin}, we define the errors $e_k:=x_k-x^*$ and $e_k^\delta:=x_k^\delta-x^*$. 
The notation $c$, with or without a subscript, denotes
a generic non-negative constant, which may differ at each occurrence, but it is always independent of the noise level $\delta$ and the iteration number $k$.
The rest of the paper is organized as follows. First, the main results {(Theorems \ref{thm:conv-noisy} and \ref{thm:err-total})} along with relevant discussions are presented in Section \ref{sec:main}. Then, the detailed proofs and remarks on the regularizing property {(Theorem \ref{thm:conv-noisy})} and convergence rate analysis {(Theorem \ref{thm:err-total})} are given in Sections \ref{sec:conv} and \ref{sec:rate} respectively. 
{For both main results concerning noisy data, the corresponding theorems derived from exact data, which are based on simpler settings and therefore easier to analyze, are discussed first and then extended to the noisy case.}
Several numerical experiments showing the advantages of the data-driven SGD over the standard SGD and Landweber method are provided in Section \ref{sec:numer} to complement the analysis. Finally, this work is concluded with further discussions in Section \ref{sec:conc}. For better readability, a set of supplementary estimates as well as lengthy technical proofs of several results are deferred to the appendix \ref{app:estimate}. 

\section{Main results and discussions}\label{sec:main}

Suitable conditions are crucial for analyzing the convergence of the data-driven SGD in Algorithm \ref{alg:datasgd} for nonlinear inverse problems.
Both the nonlinearity of the forward operator and the source condition of the solution are often employed to establish the regularizing property and convergence rate analysis \cite{HankeNeubauerScherzer:1995,JinZhouZou:2020,AspriBanert:2020}. 
Since the solution to problem \eqref{eqn:nonlin} may be nonunique, the reference solution $x^\dag$ is taken to be the minimum norm solution (with respect to the initial guess), which is known to be unique under Assumption \ref{ass:sol}(ii) below \cite{HankeNeubauerScherzer:1995,JinZhouZou:2020,AspriBanert:2020}.
Below we shall make several assumptions on the nonlinear forward operator $F_i$, the data-driven operator $G_i$, and the reference solution $x^\dag$.

\begin{assumption}\label{ass:sol}
Let $\mathcal{B}_\rho(x^\dag)\subset \bigcap_{i=1}^n\mathcal{D}(F_i)$ be a closed ball of sufficiently large radius $\rho\geq\|x_1-x^\dag\|$ and center $x^\dag$, where $x_1$ denotes the initial guess and $x^\dag$ denotes the reference solution with respect to $x_1$. 
The following conditions hold:
\begin{itemize}
  \item[$\rm(i)$] {The operators $F_i$ and $G_i$}, $i=1,\ldots,n$, have continuous and locally uniformly bounded
  Frech\'{e}t derivatives $F'_i:x\in\mathcal{D}(F_i)\subset X \rightarrow F'_i(x)\in \mathcal{L}(X,Y)$ and $G'_i:x\in X \rightarrow G'_i(x)\in \mathcal{L}(X,Y)$ on $\mathcal{B}_\rho(x^\dag)$ respectively. 
  We denote $$\max_i\sup_{x\in  \mathcal{B}_\rho(x^\dag)}\|F_i'(x)\|\leq L_F \mbox{ and } \max_i\sup_{x\in  \mathcal{B}_\rho(x^\dag)}\|G_i'(x)\|\leq L_G$$
  with Lipschitz constants $L_F$ and $L_G$.
  \item[$\rm(ii)$] {$\rm{(Tangential\; cone\; condition).}$ There exists an $\eta_F \in[0,1)$ such that, for any $i=1,\dots,n$, and any $x,\tilde x\in  \mathcal{B}_\rho(x^\dag)$,
  \begin{equation*}
    \|F_i(x)-F_i(\tilde x)-F_i'(\tilde x)(x-\tilde x)\|\leq \eta_F\|F_i(x)-F_i(\tilde x)\|.
  \end{equation*}}
  \item[$\rm(iii)$] The data-driven operator $G$ can only partially 
  explain the model for the true data, hence
  \begin{equation*}
     C_{min}\leq\|G(x^*)-y^\dag\|\leq C_{max}
  \end{equation*}
  with some constants $C_{max}\geq C_{min}>0$ for any solution $x^*$ to problem \eqref{eqn:nonlin} in $\mathcal{B}_\rho(x^\dag)$.
  \item[$\rm(iv)$]$\rm{(Range\; invariance\; condition).}$ For the operator $H=F$ or $G$, we define
\begin{equation*}
  K_{H,i}=H_i'(x^\dag),  \quad K_{H} = \frac{1}{\sqrt{n}}\left(K_{H,1}, \cdots , K_{H,n}\right)^T \quad \mbox{and}\quad B_{H}=K_H^*K_H = \frac{1}{n}\sum_{i=1}^nK_{H,i}^*K_{H,i}.
\end{equation*} 
There exists a family of locally uniformly bounded operators $R_{H,x}^i$ such that for any $x\in  \mathcal{B}_\rho(x^\dag)$,
  \begin{equation*}
    H_i'(x) = R_{H,x}^i H_i'(x^\dag)=R_{H,x}^i K_{H,i} .
  \end{equation*}
  Let $R_{H,x}=\mathrm{diag}(R_{H,x}^i):Y^n\to Y^n$, then {\rm(}with $\|\cdot\|$ denoting the operator norm on $Y^n${\rm)}
  \begin{equation*}
    \|R_{H,x}-I\|\leq c_{H}\|x-x^\dag\|.
  \end{equation*}
\item[$\rm(v)$] 
{The operator $K_F(\cdot)$ is compact, with $\{\sigma_j, \varphi_j, \psi_j\}_{j=1}^\infty$ being the singular values and vectors such that $K_F(\cdot)=\sum_{j=1}^\infty \sigma_j\langle\varphi_j,\cdot\rangle\psi_j$.}
There exists a compact operator $R$ given by $R(\cdot)=\sum_{j=1}^\infty \hat\sigma_j\langle\psi_j,\cdot\rangle\tilde{\psi_j}$ with $\{\tilde{\psi_j}\}_{j=1}^\infty$ being an orthonormal sequence in $Y^n$ such that $\|R\|\leq c_R$ and $K_G=RK_F$. That is the compact operator $K_G(\cdot)=\sum_{j=1}^\infty \tilde{\sigma_j}\langle\varphi_j,\cdot\rangle\tilde{\psi_j}$, where $\tilde\sigma_j=\sigma_j\hat\sigma_j$.

\end{itemize}
\end{assumption}

\begin{assumption}[Source condition]\label{ass:source}
There exist some $\nu\in (0,\frac12)$ and $w\in X$ such that
  \begin{equation*}
    x^\dag - x^\delta_1 = x^\dag - x_1=(F'(x^\dag)^*F'(x^\dag))^\nu w {\quad \mbox{and}\quad \|w\|< \infty}.
  \end{equation*}
\end{assumption}

The conditions in Assumptions \ref{ass:sol}(i)(ii)(iv) and \ref{ass:source} are standard for analyzing iterative regularization methods for nonlinear inverse problems \cite{HankeNeubauerScherzer:1995,EnglHankeNeubauer:1996,KaltenbacherNeubauerScherzer:2008,JinZhouZou:2020}. 
Assumptions \ref{ass:sol}(i)(ii)(iv) have been verified for a class of nonlinear inverse problems \cite{HankeNeubauerScherzer:1995}, e.g., nonlinear integral (Hammerstein) equations and parameter identification for PDEs. 
{ Especially, the tangential cone condition in (ii), which ensures the convergence of many iterative methods, is satisfied locally for the inverse problem of determining the diffusion coefficient in a parabolic partial differential equation \cite{DeCezaroZubelli2013}, time-domain full waveform inversion (FWI) in both the acoustic regime \cite{EllerRieder2021} and the elastic regime \cite{EllerGriesmaierRieder2022}, and the electrical impedance tomography (EIT) problem under suitable criteria \cite{Kindermann2022}.}
{Both the tangential cone condition in (ii) and the range invariance condition in (iv) describe some restrictions on the nonlinearity of the operators. 
Roughly speaking, it requires $F$ to be not far from linear operators on a close ball $\mathcal{B}_\rho(x^\dag)$; see Lemmas \ref{lem:linear} and \ref{lem:R} for the consequences. In particular, the radius $\rho$ can be specified as 
$\rho=\big(e^{n\sum_{j=1}^{k(\delta)} c_j^\delta}\big(\|e_1\|^2+(C_{max}+\delta)^2+n\delta^2\sum_{j=1}^{k(\delta)} d_j\big)-(C_{max}+\delta)^2\big)^\frac12<\infty$ (with the constants $c_j^\delta=2\eta_j\lambda_j^\delta \max(1,L_G^2)(\tfrac32+2\eta_j\lambda_j^\delta L_G^2)$  and $d_j=\frac{(1+\eta_F)^2}{2(1-L_F^2\eta_j-\eta_F)}\eta_j$)
under the assumptions in Theorem \ref{thm:conv-noisy}. These assumptions guarantee that all iterates $x_k^\delta$ (before stopping) are contained in $\mathcal{B}_\rho(x^\dag)$; see Corollaries \ref{cor:mono}, \ref{cor:mono_delta} for details.
Smaller $\eta_F$ corresponds to a lower degree of nonlinearity. 
In particular, when the inverse problem is linear, the constant $\eta_F=0$.} Assumptions \ref{ass:sol}(iii) and (v) assume that the data-driven operator $G$ can catch some important features of $F$, but is not able to fully approximate the forward operator for the true data. {Specifically, }(v) suggests the singular value decomposition of $K_F$ and $K_G$ , i.e., the Frech\'{e}t derivatives of $F$ and $G$ at the reference solution $x^\dag$ respectively, which share the same orthonormal basis of $X$ with the singular value $\tilde\sigma_j\leq c_R\sigma_j$ for any $j\in\mathbb{N}$.
This assumption is used to derive a simplified recursion of the error for the data-driven SGD iterate; see Section \ref{sec:rate}. 
In fact, as $G$ is an approximation of $F$, we can always design a model $G$ such that  $K_G(\cdot)=\sum_{j=1}^{\infty} \tilde{\sigma_j}\langle\tilde\varphi_j,\cdot\rangle\tilde{\psi_j}$ with the singular value $0\leq\tilde\sigma_j\leq \mathcal{O}(\sigma_j)$ and the orthonormal sequence $\{\tilde\varphi_j\}_{j=1}^{\infty}$ in $X$ satisfying $\sup_{j}\|\tilde\varphi_j-\varphi_j\|\leq\epsilon_G$ for sufficient small $\epsilon_G>0$; see consequences with this assumption in Remarks \ref{rem:basispertu1}--\ref{rem:basispertu5}. It is worth noting that this approximate basis $\{\tilde\varphi_j\}_{j=1}^{\infty}$ can be independent of noisy observations $y^\delta$.
{In practice, (iii) is satisfies by any bounded data-driven operator, while (v) can be nearly fulfilled by many types of data-driven models, including data-driven reduced order models \cite{XieMohebujjaman2018,MouKocSan2021}, neural networks combined with model reduction \cite{BhattacharyaHosseiniStuart:2021,JinZhouZou:2024} and autoencoder neural networks \cite{Kashima2016}.}
Assumption \ref{ass:source} is the classical source condition, which represents a type of smoothness of the initial error $x^\dag-x_1$ (with respect to the operator $F'(x^\dag)^*F'(x^\dag)$). 
Without this condition, the convergence of the iterative methods can be arbitrarily slow; see \cite{EnglHankeNeubauer:1996}. 
In this work, we focus on the case when $\nu\in(0,\frac12)$, where the problems have non-smooth initial errors, in the sense that the initial errors contain several relatively high-frequency components. 
For the problems with smooth initial errors (when $\nu\geq \frac12$), both the Landweber method and SGD suffer from an undesirable saturation phenomenon, i.e., the achievable accuracy does not further improve as $\nu$ grows, since the pleasant smoothness of the initial error will not be carried over to the second (and subsequent) iterations; see \cite{EnglHankeNeubauer:1996}. 

{When analyzing the convergence behavior of an iterative method, the choice of algorithmic parameters also plays an essential role. 
We shall give two classes of choices for the algorithmic parameters, including the step size schedule $\{\eta_k\}_{k=1}^\infty$ and the regularization parameter schedule $\{\lambda_k^\delta\}_{k=1}^\infty$, in the following assumption.}

\begin{assumption}\label{ass:stepsize}
The step sizes $\{\eta_k\}_{k\geq 1}$ and the regularization parameters $\{\lambda_k^\delta\}_{k=1}^\infty$ satisfy one of the following conditions.
\begin{itemize}
\item[$\rm(i)$] {$L_F^2\eta_k <1$,} $\sum_{k=1}^\infty \eta_k =\infty$ and $\sum_{k=1}^\infty \eta_k\lambda_k^{\delta}<\infty$.
\item[$\rm(ii)$] $\eta_k=\eta_0k^{-\alpha}$ and $\lambda_k^\delta\leq \lambda_0^\delta k^{-(1-\alpha)}$, with $\alpha\in(0,1)$ { and $\eta_0(L_F^2+L_G^2 \lambda_0^\delta) <1$}. 
When Assumptions \ref{ass:source} and \ref{ass:stoch} hold, the restriction on $\lambda_k^\delta$ can be relaxed to $\lambda_k^\delta\leq \lambda_0^\delta k^{-\frac12(1-\alpha+(1+\theta)\min(2\nu(1-\alpha),\alpha))}$ with some small $\theta\in(0,\max(1,(2\nu)^{-1}-1,\alpha^{-1}-2))$ defined in Assumption \ref{ass:stoch}.
\end{itemize}
\end{assumption}

{The choice in Assumption \ref{ass:stepsize}(i) used for establishing the regularizing property in Theorem \ref{thm:conv-noisy} is more general than that in (ii) (without Assumptions \ref{ass:source} and \ref{ass:stoch}) used for deriving the convergence rates in Theorem \ref{thm:err-total}. }
The latter choice is often known as the polynomially decaying step size and regularization parameter schedules in the literature. 
When Assumptions \ref{ass:source} and \ref{ass:stoch} hold, the relaxed assumption on the regularization parameter schedule $\{\lambda_k^\delta\}_{k=1}^\infty$ in (ii) makes $\sum_{k=1}^\infty \eta_k\lambda_k^{\delta}=\infty$, which is contrary to the assumption in (i), but still enables the algorithm to achieve the same convergence rate as that obtained under the stronger assumption in (ii); see Theorems \ref{thm:err-total-ex} and \ref{thm:err-total} for both exact and noisy data.

Due to the random choice of the index $i_k$ at each iteration, the data-driven SGD iterate $x_k^\delta$ is random. 
We denote the filtration generated by the random indices $\{i_1,\ldots,i_{k-1}\}$ up to the $(k-1)$-th iteration by $\mathcal{F}_k$. 
Among different ways to measure the convergence, we consider the mean squared norm defined by $\E[\|\cdot\|^2]$, where the expectation $\E[\cdot]$ is with respect to the filtration $\mathcal{F}_k$. 
Note that the iterate $x_k^\delta$ is measurable with respect to $\mathcal{F}_k$. 
The first result presents the regularizing property of the data-driven SGD for problem \eqref{eqn:nonlin} under the tangential cone condition and \textit{a priori} parameter choice. 
The additional assumptions in Theorem \ref{thm:conv-noisy} on the regularization parameter $\lambda_k^\delta$, comparing with that for the standard SGD \cite{JinZhouZou:2020}, is due to the presence of data-driven operators in the regularization term which may lead to learning errors (as the data-driven operator can only partially explain the exact model) at each iteration. 
It is worth noting that these assumptions in Theorem \ref{thm:conv-noisy} are more relaxed than that for the data-driven iteratively regularized Landweber method \cite{AspriBanert:2020}. 
In particular, we adopt \textit{a priori} selection scheme for the regularization parameter $\lambda_k^\delta$ which is independent of the residuals of the algorithm and subsumed by the assumptions in \cite{AspriBanert:2020}. {In addition,  the conditions on the forward operator $F$ in Theorem \ref{thm:conv-noisy} are assumed to hold within the closed ball $x^*\in\mathcal{B}_\rho(x^\dag)$, rather than the entire space as assumed in \cite{JinZhouZou:2020}.}

\begin{theorem}[Convergence for noisy data]\label{thm:conv-noisy}
Let Assumptions \ref{ass:sol}(i)-(iii) and \ref{ass:stepsize}(i) be fulfilled  {with $L_F^2\eta_k<1-\eta_F$ for any $k\geq1$}. If the condition $\lim_{\delta\to 0^+}\lambda_k^\delta=\lambda_k^0$ holds for any $k\in\mathbb{N}$ and the stopping index $k(\delta)\in\mathbb{N}$ is chosen such that
\begin{equation*}
 \lim_{\delta\to0^+}k(\delta)=\infty \quad\mbox{and}\quad \lim_{\delta\to0^+}\delta^2\sum_{i=1}^{k(\delta)}\eta_i = 0,
\end{equation*} 
then for the data-driven SGD iterate $x_k^\delta$ in \eqref{eqn:datasgd}, there exists a solution {$x^*\in\mathcal{B}_\rho(x^\dag)$} to problem \eqref{eqn:nonlin} such that
\begin{equation*}
  \lim_{\delta\to 0^+}\E[\|x_{k(\delta)}^\delta - x^*\|^2]=0.
\end{equation*}
Further, if $\mathcal{N}(F'(x^\dag))\subset\mathcal{N}(F'(x))$ and $\mathcal{N}(F'(x^\dag))\subset\mathcal{N}(G'(x))$ (with $\mathcal{N}(\cdot)$ denoting the kernel of the linear operator) {for any $x\in\mathcal{B}_\rho(x^\dag)$}, then
\begin{equation*}
  \lim_{\delta\to0^+}\E[\|x_{k(\delta)}^\delta -x^\dag\|^2]=0.
\end{equation*}
\end{theorem}

Next, we make an assumption, which is a stochastic variant of the range invariance condition stated in Assumption \ref{ass:sol}(iv), on the degree of nonlinearity of the operators $F$ and $G$ in the sense of expectation. This assumption is crucial for deriving the convergence rates of the data-driven SGD in Section \ref{sec:rate} due to the presence of conditionally dependent factors in the computational variance; see the proof of Lemma \ref{lem:bound-N} (and Lemma \ref{lem:R}).
\begin{assumption}[Stochastic range invariance condition]\label{ass:stoch}
With the notations defined in Assumption \ref{ass:sol}(iv), for the operator $H=F$ or $G$, there exist some small $\theta\in(0,1)$ and $c_H>0$ such that, for function 
$Q(x^\delta_k)=K_H(x^\delta_k-x^\dag)$ or $Q(x^\delta_k)=H(x^\delta_k)-y^\delta$, $x_k^\delta\in\mathcal{B}_\rho(x^\dag)$ and
$z_t= tx_k^\delta + (1-t)x^\dag$, $t\in[0,1]$, there hold
\begin{align*}
  \E[\|(R_{H,z_t}-I)Q(x_k^\delta)\|^2]^\frac12 &\leq c_H\E[\|x_k^\delta-x^\dag\|^2]^\frac{\theta}{2}\E[\|Q(x_k^\delta)\|^2]^\frac12,\\
  \E[\|(R_{H,z_t}^*-I)Q(x_k^\delta)\|^2]^\frac12 &\leq c_H\E[\|x_k^\delta-x^\dag\|^2]^\frac{\theta}{2}\E[\|Q(x_k^\delta)\|^2]^\frac12.
\end{align*}
\end{assumption}

The second result presents the convergence rates for the data-driven SGD under the additional source condition, range invariance conditions and \textit{a priori} parameter choice. 
It shares a similar general strategy to the error analysis of the standard SGD in \cite{JinLu:2019} and \cite{JinZhouZou:2020} for linear and nonlinear inverse problems respectively. 
We decompose the total error $\E[\|x_k^\delta-x^\dag\|^2]$ into two components, i.e., the mean iterate error $\|\E[x_k^\delta]-x^\dag\|^2$ dominated by the approximation error and data error and its computational variance $\E[\|x_k^\delta-\E[x_k^\delta]\|^2]$ caused by the randomness of gradient estimates, by the standard bias-variance decomposition:
\begin{equation}\label{eqn:bias-var}
  \E[\|x_k^\delta - x^\dag\|^2] = \|\E[x_k^\delta]-x^\dag\|^2 + \E[\|x_k^\delta-\E[x_k^\delta]\|^2].
\end{equation}
Since the data-driven operator in the algorithm  originally introduces learning errors into both components at each iteration, the analysis differs significantly from that of the standard SGD \cite{JinLu:2019,JinZhouZou:2020}; see
Theorems \ref{thm:err-mean} and \ref{thm:err-var} for the bias and variance respectively. 
Similar to the observation for the standard SGD in \cite{JinZhouZou:2020}, these two components closely interact with each other due to the nonlinearity of the operators, resulting in a coupled system of recursions on the estimates of $\E[\|e_k^\delta\|^2]$ and $\E[\|F'(x^\dag)e_k^\delta\|^2]$.
Finally, we obtain an error analysis in the following result by mathematical induction. 
\begin{theorem}{[Convergence rates for noisy data]}\label{thm:err-total}
Let Assumptions \ref{ass:sol}, \ref{ass:source}, \ref{ass:stepsize}(ii) and \ref{ass:stoch} be fulfilled with $\|w\|$, $\eta_0$ and $\lambda_0^\delta$
being sufficiently small, $\|F'(x^\dag)^*F'(x^\dag)\|\leq1$, and $x_k^\delta$ be the data-driven SGD iterate defined by \eqref{eqn:datasgd}. Then the error
$e_k^\delta = x_k^\delta - x^\dag$ satisfies
\begin{align*}
  \E[\|e_k^\delta\|^2] & \leq c^* k^{-\min(2\nu(1-\alpha),\alpha)}\|w\|^2\quad \mbox{and}\quad  \E[\|F'(x^\dag)e_k^\delta\|^2] \leq c^*k^{-\min((1+2\nu)(1-\alpha),1)}\|w\|^2
\end{align*}
for all $k\leq k^*=[(\frac{\delta}{\|w\|})^{-\frac{2}{\min((1+2\nu)(1-\alpha),1)+\epsilon}}]$ (with $[\cdot]$ denoting taking the integral part of a real number) and small 
$\epsilon\in(0,2\max((1-2\nu)(1-\alpha),1-2\alpha))$,
where the constant $c^*$ depends on $\nu$, $\alpha$, $\eta_0$, $n$, $\theta$ 
and $\epsilon$, but is independent of $k$ and $\delta$.
\end{theorem}
The presence of the parameter $\epsilon$ in the stopping index $k^*$ is caused by the data and stochastic gradient noises which introduce data and stochastic errors at each iteration into both bias and variance (as they closely interact with each other). 
When the noise level $\delta=0$, i.e., using exact data, we achieve the same upper bounds of both $\E[\|e_k\|^2]$ and $\E[\|F'(x^\dag)e_k\|^2]$ where the constant $c^*$ is independent of $\epsilon$; see Theorem \ref{thm:err-total-ex}.
When $\epsilon$ in Theorem \ref{thm:err-total} is sufficiently small, setting {$\alpha\in[\frac{2\nu}{1+2\nu},1)$ and }$k=k^*$ provides the following convergence rate in terms of the noise level:
\begin{equation*}
\E[\|e_{k^*}^\delta\|^2] \leq c^* \|w\|^\frac{2}{1+2\nu}\delta^\frac{4\nu}{1+2\nu},
\end{equation*}
which is comparable with that for the Landweber method in \cite{HankeNeubauerScherzer:1995} and the standard SGD for both linear and nonlinear cases in \cite{JinLu:2019} and \cite{JinZhouZou:2020} respectively. 

\section{Convergence of the data-driven SGD}\label{sec:conv}    

In this section, we analyze the convergence of the data-driven SGD iterate defined in Algorithm \ref{alg:datasgd} for exact and noisy data in Theorems \ref{thm:conv-exact} and \ref{thm:conv-noisy} respectively. 
{First, we present a result that suggests an almost non-expansiveness property of the iterate errors under suitable assumptions on the regularization parameters.}
This result is crucial for proving the regularizing property of the data-driven SGD iterates under
\textit{a priori} parameter choice; see the appendix \ref{app:prop:mono-exact} for the proof. Below $x^\dag$ denotes the unique reference solution to problem \eqref{eqn:nonlin} of minimal
distance to $x_1$.
{\begin{prop}\label{prop:mono-exact}
Let Assumptions \ref{ass:sol}(i)-(iii) and \ref{ass:stepsize}(i) be fulfilled with $L_F^2\eta_k<1-\eta_F$ for any $k\geq1$. Then for any data-driven SGD iterate $x_k^\delta\in\mathcal{B}_\rho(x^\dag)$ defined by \eqref{eqn:datasgd}
, the error $e_k^\delta=x_k^\delta-x^\dag$ satisfies
\begin{align}
\|e_{k+1}^\delta\|^2 
\leq &(1+nc_k^\delta)\|e_k^\delta\|^2+nc_k^\delta(C_{max}+\delta)^2+nd_k\delta^2,\label{eqn:prop_e}\\
\E[\|e_{k+1}^\delta\|^2] 
\leq &(1+c_k^\delta)\E[\|e_k^\delta\|^2]+c_k^\delta(C_{max}+\delta)^2+2(1+\eta_F)\eta_k\delta\E[\|F(x_k^\delta)-y^\delta\|^2]^\frac12\nonumber\\
&-2(1-L_F^2\eta_k-\eta_F)\eta_k\E[\|F(x_k^\delta)-y^\delta\|^2]\label{eqn:prop_Ee1}\\
\leq &(1+c_k^\delta)\E[\|e_k^\delta\|^2]+c_k^\delta(C_{max}+\delta)^2+d_k\delta^2, \label{eqn:prop_Ee2} 
\end{align}
where the constants $c_k^\delta=2\eta_k\lambda_k^\delta \max(1,L_G^2)(\tfrac32+2\eta_k\lambda_k^\delta L_G^2)$  and $d_k=\dfrac{(1+\eta_F)^2}{2(1-L_F^2\eta_k-\eta_F)}\eta_k$.
\end{prop}}

Below we analyze the convergence of the data-driven SGD iterate for exact and noisy data separately.

\subsection{Convergence for exact data}

The next result, showing that the sequence of mean squared errors $\{\E[\|e_{k}\|^2]\}_{k\geq 1}$ is a Cauchy sequence and all iterates $\{x_k\}_{k\geq 1}$ are contained in the closed ball $\mathcal{B}_\rho(x^\dag)$, follows directly from Proposition \ref{prop:mono-exact}.

\begin{corollary}\label{cor:mono}
Let Assumptions \ref{ass:sol}(i)-(iii) and \ref{ass:stepsize}(i) be fulfilled with {$L_F^2\eta_k<1-\eta_F$ for any $k\geq1$, and $\rho^2=e^{n\sum_{k=1}^\infty c_k^0}(\|e_1\|^2+C_{max}^2)-C_{max}^2$, where the constant $c_k^0=2\eta_k\lambda_k^0 \max(1,L_G^2)(\tfrac32+2\eta_k\lambda_k^0 L_G^2)$}. Then for the data-driven SGD iterate $x_k$ in \eqref{eqn:datasgd} with the exact data $y^\dag$,
the error $\{\E[\|e_{k}\|^2]=\E[\|x_k-x^\dag\|^2]\}_{k\geq 1}$ is a Cauchy sequence that converges to some constant $C_e\geq 0$, $x_k\in \mathcal{B}_\rho(x^\dag)$, and 
{\begin{align*}
\sum_{k=1}^\infty \eta_k\E[\|F(x_k)-y^\dag\|^2] &\leq \tfrac12(1-L_F^2\eta_k-\eta_F)^{-1}\big(\|e_1\|^2+(\rho^2+C_{max}^2)\sum_{k=1}^\infty c_k^0\big)<\infty.
\end{align*}}
\end{corollary}
\begin{proof}
{Under the condition $\sum_{k=1}^\infty \eta_k\lambda_k^{0}<\infty$ in Assumption \ref{ass:stepsize}(i), which gives the estimate $\sum_{k=1}^\infty c_k^0< \infty$, we specify the radius $\rho$ in Assumption \ref{ass:sol} as $\rho^2=e^{n\sum_{k=1}^\infty c_k^0}(\|e_1\|^2+C_{max}^2)-C_{max}^2<\infty$}.
First, we shall show all iterates $x_k$ in \eqref{eqn:datasgd} remain in the closed ball $\mathcal{B}_\rho(x^\dag)$ by mathematical induction. {For the case $k=1$, $x_1\in \mathcal{B}_\rho(x^\dag)$ by the relation $\|x_1-x^\dag\|^2=\|e_1\|^2\leq \rho^2$.} Now, we assume that $x_k\in \mathcal{B}_\rho(x^\dag)$ hold up to the case $k$, and prove the assertion for the case $k+1$.
{By the recursion \eqref{eqn:prop_e} in Proposition \ref{prop:mono-exact} with $\delta=0$, there holds
\begin{align*}
\|e_{k+1}\|^2 \leq (1+nc_k^0)\|e_k\|^2+nc_k^0 C_{max}^2,
\end{align*}
which directly indicates that
\begin{align*}
\|e_{k+1}\|^2+ C_{max}^2\leq (1+nc_k^0)(\|e_k\|^2+C_{max}^2)\leq \Pi_{j=1}^k (1+nc_j^0)(\|e_1\|^2+C_{max}^2).
\end{align*}
Further, by using the fact $1+x\leq e^{x}$ for any $x\geq 0$, we bound the iterate error $\|x_{k+1}-x^\dag\|^2=\|e_{k+1}\|^2$ by
\begin{align*}
\|e_{k+1}\|^2\leq e^{n\sum_{j=1}^k c_j^0}(\|e_1\|^2+C_{max}^2)-C_{max}^2\leq \rho^2,
\end{align*}
i.e., $x_{k+1}\in \mathcal{B}_\rho(x^\dag)$. Next, by the recursion \eqref{eqn:prop_Ee2} in Proposition \ref{prop:mono-exact} with $\delta=0$ and the previous result $\E[\|e_k\|^2]\leq \rho^2$ for any $k\geq1$, the difference between two successive iterate errors can be bounded by
\begin{align*}
\E[\|e_{k+1}\|^2]-\E[\|e_{k}\|^2]  
\leq &c_k^0\E[\|e_k\|^2]+c_k^0 C_{max}^2\leq c_k^0(\rho^2+ C_{max}^2),
\end{align*}
which implies that, for any $\ell>i$,
\begin{align*}
\E[\|e_{\ell}\|^2] -\E[\|e_{i}\|^2]=\sum_{k=i}^{\ell-1}(\E[\|e_{k+1}\|^2] -\E[\|e_{k}\|^2])
\leq (\rho^2+ C_{max}^2)\sum_{k=i}^{\ell-1}c_k^0.
\end{align*}
With the estimate $\sum_{k=1}^\infty c_k^0< \infty$ derived from Assumption \ref{ass:stepsize}(i), we obtain that 
\begin{align*}
\lim_{i<\ell,i\to\infty}\E[\|e_{\ell}\|^2] -\E[\|e_{i}\|^2]
\leq (\rho^2+ C_{max}^2) \lim_{i<\ell,i\to\infty}\sum_{k=i}^{\ell-1}c_k^0=0,
\end{align*}}
which implies that $\{\E[\|e_{k}\|^2]\}_{k\geq 1}$ is a Cauchy sequence. Furthermore, the fact that $\E[\|e_{k}\|^2]\geq 0$ guarantees that $\lim_{k\to \infty}\E[\|e_{k}\|^2]:=C_e\geq 0$.

{Similarly, the recursion \eqref{eqn:prop_Ee1} in Proposition \ref{prop:mono-exact} with $\delta=0$ gives
\begin{align*}
2(1-L_F^2\eta_k-\eta_F)\eta_k\E[\|F(x_k)-y^\dag\|^2] \leq& (1+c_k^0)\E[\|e_{k}\|^2] - \E[\|e_{k+1}\|^2]+c_k^0 C_{max}^2\\ 
\leq &\E[\|e_{k}\|^2] - \E[\|e_{k+1}\|^2]+c_k^0 (\rho^2+C_{max}^2),
\end{align*}
and thus
\begin{align*}
2(1-L_F^2\eta_k-\eta_F)\sum_{k=1}^\infty\eta_k\E[\|F(x_k)-y^\dag\|^2]&\leq  \|e_1\|^2+(\rho^2+C_{max}^2)\sum_{k=1}^\infty c_k^0<\infty.
\end{align*}}
This completes the proof of the corollary.
\end{proof}

The next result shows that the sequence $\{x_k\}_{k\geq 1}$ is a Cauchy sequence in $\mathcal{B}_\rho(x^\dag)$; see the appendix \ref{app:prop:Cauchy} for the proof.
\begin{prop}\label{prop:Cauchy}
Let Assumptions \ref{ass:sol}(i)-(iii) and \ref{ass:stepsize}(i) be fulfilled {with $L_F^2\eta_k<1-\eta_F$ for any $k\geq1$}.
Then for the exact data $y^\dag$,
the sequence $\{x_k\}_{k\geq 1}$ defined by \eqref{eqn:datasgd} is a Cauchy sequence in $\mathcal{B}_\rho(x^\dag)$.
\end{prop}

Now, we can state the convergence of the data-driven SGD iterate in Algorithm \ref{alg:datasgd} for the exact data $y^\dag$.
\begin{theorem}[Convergence for exact data]\label{thm:conv-exact}
Let Assumptions \ref{ass:sol}(i)-(iii) and \ref{ass:stepsize}(i) be fulfilled {with $L_F^2\eta_k<1-\eta_F$ for any $k\geq1$}.
Then for the exact data $y^\dag$, the data-driven SGD sequence
$\{x_k\}_{k\geq1}$ defined in \eqref{eqn:datasgd} converges to a solution $x^*\in\mathcal{B}_\rho(x^\dag)$ to problem \eqref{eqn:nonlin}:
\begin{equation*}
  \lim_{k\to\infty}\E[\|x_k-x^*\|^2] = 0.
\end{equation*}
Further, if
$\mathcal{N}(F'(x^\dag))\subset\mathcal{N}(F'(x))$ and $\mathcal{N}(F'(x^\dag))\subset \mathcal{N}(G'(x))$ for any $x\in\mathcal{B}_\rho(x^\dag)$, then
\begin{equation*}
  \lim_{k\to\infty}\E[\|x_k -x^\dag\|^2]=0.
\end{equation*}
\end{theorem}
\begin{proof}
The argument below follows closely \cite[Lemma 3.4 and Theorem 3.5]{JinZhouZou:2020}. For the convenience of readers, we state similar results here. 
By Lemma \ref{lem:linear} and Assumption \ref{ass:sol}(i), for any $x,\tilde x\in \mathcal{B}_\rho(x^\dag)$, there holds
\begin{equation*}
  \|(F(x)-y^\dag)-(F(\tilde x)-y^\dag)\|=\|F(x)-F(\tilde x)\| \leq (1-\eta_F)^{-1}\|F'(x)(x-\tilde x)\|\leq L_F(1-\eta_F)^{-1}\|x-\tilde x\|.
\end{equation*}
{ Thus, by Proposition \ref{prop:Cauchy} (i.e., the fact that $\{x_k\}_{k\geq 1}$ is a Cauchy sequence in $\mathcal{B}_\rho(x^\dag)$), we obtain that
$\{F(x_k)-y^\dag\}_{k\geq 1}$ is a Cauchy sequence that converges to $F(x^*)-y^\dag$ with $x^*:=\lim_{k\to \infty}x_k\in\mathcal{B}_\rho(x^\dag)$.}
Furthermore, the fact that $\E[\|F(x_k)-y^\dag\|^2]\geq 0$ guarantees that $\lim_{k\to \infty}\E[\|F(x_k)-y^\dag\|^2]:=\epsilon_r\geq 0$. There exists some $k_0\in \mathbb{N}$, such that, for any $k\geq k_0$, $\E[\|F(x_k)-y^\dag\|^2]\geq \frac12\epsilon_r$. If $\epsilon_r>0$, Assumption \ref{ass:stepsize}(i) leads to the inequality 
\begin{equation*}
  \sum_{k=1}^\infty \eta_k\E[\|F(x_k)-y^\dag\|^2] \geq \sum_{k=k_0}^\infty \eta_k\E[\|F(x_k)-y^\dag\|^2] \geq \frac12\epsilon_r\sum_{k=k_0}^\infty\eta_k=\infty,
\end{equation*}
which contradicts the result in Corollary \ref{cor:mono} that
\begin{equation*}
\sum_{k=1}^\infty \eta_k\E[\|F(x_k)-y^\dag\|^2] <\infty.
\end{equation*}
Thus, we have $\E[\|F(x^*)-y^\dag\|^2]=\lim_{k\to\infty} \E[\|F(x_k)-y^\dag\|^2]=\epsilon_r=0$
which implies that $x^*$ is a solution to problem \eqref{eqn:nonlin}.

Further, Lemma \ref{lem:linear}(ii) indicates that there exists a unique reference solution to problem \eqref{eqn:nonlin} such that 
\begin{equation*}
  x^\dag - x_1\in\mathcal{N}(F'(x^\dag))^\perp.
\end{equation*} 
If $\mathcal{N}(F'(x^\dag))\subset \mathcal{N}(F'(x_k))$ and $\mathcal{N}(F'(x^\dag))\subset \mathcal{N}(G'(x_k))$ for any $k\geq 1$, then with the definition of $x_{k}$ in \eqref{eqn:datasgd} for exact data, we have
\begin{equation*}
  x_{k+1}-x_1 = \sum_{i=1}^k(x_{i+1}-x_i)=- \sum_{i=1}^k \eta_i \big(F_{i_i}'(x_i)^*(F_{i_i}(x_i)-y_{i_i}^\dag)+\lambda_i^0 G_{i_i}'(x_i)^*(G_{i_i}(x_i)-y_{i_i}^\dag)\big)\in\mathcal{N}(F'(x^\dag))^\perp.
\end{equation*}
Combining the above two observations, we derive that 
\begin{equation*}
  x^\dag - x^* = x^\dag - x_1 + x_1-x^*\in\mathcal{N}(F'(x^\dag))^\perp.
\end{equation*}
Further, by Lemma \ref{lem:linear}(ii), there holds $x^\dag-x^* \in \mathcal{N}(F'(x^\dag))$ which implies $x^\dag-x^*=0$. This completes the proof of the theorem.
\end{proof}

\subsection{Convergence for noisy data}
The next result, showing that the iterates $\{x_k^\delta\}_{k=1}^{ k(\delta)}$ are contained in the closed ball $\mathcal{B}_\rho(x^\dag)$ (where $k(\delta)$ denotes the stopping index defined in Theorem \ref{thm:conv-noisy}), follows directly from Proposition \ref{prop:mono-exact}.

\begin{corollary}\label{cor:mono_delta}
Let assumptions in Theorem \ref{thm:conv-noisy} be fulfilled with {$$\rho^2=e^{n\sum_{j=1}^{k(\delta)} c_j^\delta}\big(\|e_1\|^2+(C_{max}+\delta)^2+n\delta^2\sum_{j=1}^{k(\delta)} d_j\big)-(C_{max}+\delta)^2,$$
where the constants $c_j^\delta=2\eta_j\lambda_j^\delta \max(1,L_G^2)(\tfrac32+2\eta_j\lambda_j^\delta L_G^2)$ and $d_j=\dfrac{(1+\eta_F)^2}{2(1-L_F^2\eta_j-\eta_F)}\eta_j$.}
Then for any $k\leq k(\delta)$, the data-driven SGD iterate $x_k^\delta$ in \eqref{eqn:datasgd} is contained in $\mathcal{B}_\rho(x^\dag)$.
\end{corollary}
\begin{proof}
Under the assumptions in Theorem \ref{thm:conv-noisy}, we specify the radius $\rho$ in Assumption \ref{ass:sol} as {$$\rho^2=e^{n\sum_{j=1}^{k(\delta)} c_j^\delta}\big(\|e_1\|^2+(C_{max}+\delta)^2+n\delta^2\sum_{j=1}^{k(\delta)} d_j\big)-(C_{max}+\delta)^2<\infty.$$}
{By the recursion \eqref{eqn:prop_e} in Proposition \ref{prop:mono-exact}, there holds
\begin{align*}
 \|e_{k+1}^\delta\|^2 
\leq &(1+nc_k^\delta)\|e_k^\delta\|^2+nc_k^\delta(C_{max}+\delta)^2+n d_k\delta^2,
\end{align*}
which implies that 
\begin{align*}
\|e_2^\delta\|^2\leq(1+nc_1^\delta)\|e_1^\delta\|^2+nc_1^\delta(C_{max}+\delta)^2+n d_1\delta^2\leq e^{nc_1^\delta}\big(\|e_1\|^2+(C_{max}+\delta)^2\big)+n\delta^2 d_1-(C_{max}+\delta)^2\leq \rho^2.
\end{align*}
Similar to the analysis in the proof of Corollary \ref{cor:mono}, for any $3\leq k+1\leq k(\delta)$, we bound $r_{k+1}:=\|e_{k+1}^\delta\|^2+(C_{max}+\delta)^2$ by 
\begin{align*}
r_{k+1}\leq&(1+nc_k^\delta)r_{k}+nd_k\delta^2
\leq \Pi_{j=k-1}^k(1+nc_j^\delta)r_{k-1}+\Pi_{j=k}^k(1+nc_j^\delta)nd_{k-1}\delta^2+nd_k\delta^2\\
\leq& \cdots \leq \Pi_{j=1}^k(1+nc_j^\delta)r_{1}+\sum_{i=1}^{k-1} \Pi_{j=i+1}^k(1+nc_j^\delta)nd_i\delta^2+nd_k\delta^2
\leq \Pi_{j=1}^k(1+nc_j^\delta)(r_{1}+n\delta^2\sum_{i=1}^{k} d_i)\nonumber\\
\leq& e^{n\sum_{j=1}^k c_j^\delta}\big(\|e_1\|^2+(C_{max}+\delta)^2+n\delta^2\sum_{j=1}^{k} d_j\big)\leq \rho^2+(C_{max}+\delta)^2,\nonumber
\end{align*}}
i.e., $x_{k+1}^\delta\in \mathcal{B}_\rho(x^\dag)$.
This completes the proof of the corollary.
\end{proof}

The following result gives the pathwise (i.e., along each realization in the filtration $\mathcal{F}_k$) stability of the data-driven SGD iterate $x_k^\delta$ with respect to
the noise level $\delta$ at $\delta=0$; see the appendix \ref{app:lem:conti-noise} for the proof.
\begin{lemma}\label{lem:conti-noise}
Let assumptions in Theorem \ref{thm:conv-noisy} be fulfilled.
For any fixed $k\in \mathbb{N}$ and any path
$(i_1,\ldots,i_{k-1})\in\mathcal{F}_k$, let $x_k$ and $x_k^\delta$ be the data-driven SGD iterates along the path for exact data $y^\dag$ and noisy data $y^\delta$ respectively. Then
\begin{equation*}
  \lim_{\delta\to 0^+}\|x_k^\delta-x_k\|=0 \mbox{ and } \lim_{\delta\to 0^+}\E[\|x_k^\delta-x_k\|^2]^\frac12=0.
\end{equation*}
\end{lemma}

Now, we can give the proof of Theorem \ref{thm:conv-noisy} which gives the regularizing property of the data-driven SGD under \textit{a priori} stopping rules.
\begin{proof}[Proof of Theorem \ref{thm:conv-noisy}]
Let $\{\delta_t\}_{t\geq 1}\subset\mathbb{R}$ be a sequence converging to zero, and let $y_t:=y^{\delta_t}$
be a corresponding sequence of noisy data. For each pair $(\delta_t,y_t)$, we denote by
$k_t=k(\delta_t)$ the stopping index. 
{First, the recursion \eqref{eqn:prop_Ee2} in Proposition \ref{prop:mono-exact} and Corollary \ref{cor:mono_delta} give that
\begin{align*}
\E[\|e_{k+1}^\delta\|^2]-\E[\|e_k^\delta\|^2] 
\leq &c_k^\delta\E[\|e_k^\delta\|^2]+c_k^\delta(C_{max}+\delta)^2+d_k\delta^2
\leq c_k^\delta\big(\rho^2+(C_{max}+\delta)^2\big)+d_k\delta^2.
\end{align*}
For any $k<k_t$, summing the above inequality with $\delta=\delta_t$ from $k$ to $k_t-1$ and applying the triangle inequality lead to
\begin{align*}
\E[\|e_{k_t}^{\delta_t}\|^2] \leq &\E[\|e_{k}^{\delta_t}\|^2] + \delta_t^2\sum_{j=k}^{k_t-1}d_j+\big(\rho^2+(C_{max}+\delta)^2\big)\sum_{j=k}^{k_t-1}c_j^\delta\\
\leq& 2\E[\|x_{k}^{\delta_t}-x_{k}\|^2] + 2\E[\|x_{k}-x^*\|^2] + \delta_t^2\sum_{j=1}^{k_t}d_j+\big(\rho^2+(C_{max}+\delta)^2\big)\sum_{j=k}^{\infty}c_j^\delta.
\end{align*}
By Theorem \ref{thm:conv-exact} and the condition $\sum_{k=1}^\infty \eta_k\lambda_k^{\delta}<\infty$ in Assumption \ref{ass:stepsize}(i), for any $\epsilon>0$, there exists some $K\in\mathbb{N}_+$, such that for any $k\geq K$, we have $\E[\|x_{k}-x^*\|^2]<\frac{\epsilon}{8}$ and $\sum_{j=k}^{\infty}c_j^\delta<\frac{\epsilon}{4\big(\rho^2+(C_{max}+\delta)^2\big)}$. 
Further, for the fixed index $K$, Lemma \ref{lem:conti-noise} and the condition on the index $k_t$, i.e., $\lim_{t\to\infty}\delta_t^2\sum_{i=1}^{k_t}\eta_i=0$, guarantee that there exists some $T\in\mathbb{N}_+$, such that for any $t\geq T$, we have 
$\E[\|x_{K}^{\delta_t}-x_{K}\|^2]< \frac{\epsilon}{8}$ and $\delta_t^2\sum_{j=k}^{k_t}d_j<\frac{\epsilon}{4}$.
Now, under the condition $\lim_{t\to\infty}k_t=\infty$, we can select $k_t>K$, then there holds 
\begin{align*}
\E[\|e_{k_t}^{\delta_t}\|^2] 
\leq& 2\E[\|x_{K}^{\delta_t}-x_{K}\|^2] + 2\E[\|x_{K}-x^*\|^2] + \delta_t^2\sum_{j=1}^{k_t}d_j+\big(\rho^2+(C_{max}+\delta)^2\big)\sum_{j=K}^{\infty}c_j^\delta<\epsilon.
\end{align*}}
This completes
the proof of the first assertion. The case for $\mathcal{N}(F'(x^\dag))\subset\mathcal{N}(F'(x))$  and $\mathcal{N}(F'(x^\dag))\subset\mathcal{N}(G'(x))$ follows similarly
as Theorem \ref{thm:conv-exact}.
\end{proof}

\section{Convergence rates}\label{sec:rate}
In this section, we prove convergence rates for the data-driven SGD under Assumptions \ref{ass:sol}, \ref{ass:source}, \ref{ass:stepsize}(ii) and
\ref{ass:stoch}. The main results are given in Theorems \ref{thm:err-total-ex} and \ref{thm:err-total} for
exact and noisy data respectively. These results represent the second main contributions
of the work. We shall employ some shorthand notation. 
\begin{equation}\label{eqn:Pi-B}
\Pi_{j}^k(B)=\prod_{i=j}^k\big(I-\eta_i(B_F+\lambda_i^\delta B_G)\big), \quad {\mbox{where} \quad B_H=K_H^* K_H = H'(x^\dag)^*H'(x^\dag) \quad \mbox{for } H=F \mbox{ or } G,}
\end{equation}
(and the convention $\Pi_j^k(B)=I$ for $j>k$), and to shorten lengthy expressions, we define for
$s,\tilde{s}\geq0$ and $j\in\mathbb{N}$,
\begin{equation*}
 \tilde s = s+\tfrac12 \quad\mbox{ and }\quad\phi_j^{s} = \|B_F^{s}\Pi_{j+1}^k(B)\|.
\end{equation*}
The rest of this section is structured as follows. In view of the standard bias-variance
decomposition \eqref{eqn:bias-var}, we first derive two important recursion formulas for the weighted bias $\|B_F^s(
\E[x_k^\delta]-x^\dag)\|$ and weighted variance $\E[\|B_F^s(x_k^\delta-\E[x_k^\delta])\|^2]$, for any $s\geq0$, in Sections \ref{ssec:mean} and
\ref{ssec:residual} respectively, and then use the recursions to derive the desired convergence rates under
\textit{a priori} parameter choice in Section \ref{ssec:rate}.

\subsection{Recursion on the mean}\label{ssec:mean}
In this part, we derive a recursion for the upper bound on the weighted error of the mean data-driven SGD iterate $\E[x_k^\delta]$.
The next result gives a useful representation of the mean $\E[e_k^\delta]$ of the error $e_k^\delta=x_k^\delta-x^\dag$; see the appendix \ref{app:lem:recursion-mean} for the proof.
\begin{lemma}\label{lem:recursion-mean}
Let Assumption \ref{ass:sol}(iv) be fulfilled. Then for the data-driven SGD iterate $x_k^\delta$ in \eqref{eqn:datasgd}, the error $e_k^\delta=x_k^\delta-x^\dag$ satisfies
\begin{equation*}
  \E[e_{k+1}^\delta] = \Pi_1^k(B)e_1^\delta - \sum_{j=1}^{k}\eta_{j}\Pi_{j+1}^k(B) (K_{F}^*\E[v_{F,j}]+ \lambda_j^\delta K_{G}^*\E[v_{G,j}]),
\end{equation*}
with the vector $v_{F,j},v_{G,j}\in Y^n$ given by
\begin{align}
v_{F,j}=& ( R_{F,x_j^\delta}^{ *}-I )(F(x_j^\delta)-y^\delta)+\big(F(x_j^\delta)-F(x^\dag)-K_{F}(x_j^\delta-x^\dag)\big)- \xi,\label{eqn:w_F}\\
v_{G,j}=& ( R_{G,x_j^\delta}^{*}-I )(G(x_j^\delta)-y^\delta)+ \big(G(x_j^\delta)-G(x^\dag)-K_{G}(x_j^\delta-x^\dag)\big)+(G(x^\dag)-y^\dag)- \xi.\label{eqn:w_G}
\end{align}

\end{lemma}

\begin{remark}
If the data-driven operator $G$ is linear, under Assumption \ref{ass:sol}(iv), the vector $v_{G,j}$ in \eqref{eqn:w_G} simplifies to $$v_{G,j}=G(x^\dag)-y^\dag-\xi,$$ 
which is independent of the iterate index $j$.
\end{remark}

\begin{corollary}\label{cor:recursion-mean}
Let Assumptions \ref{ass:sol}(iv)(v) be fulfilled. Then for the data-driven SGD iterate $x_k^\delta$ in \eqref{eqn:datasgd}, the error $e_k^\delta=x_k^\delta-x^\dag$ satisfies
\begin{equation*}
  \E[e_{k+1}^\delta] = \Pi_1^k(B)e_1^\delta - \sum_{j=1}^{k}\eta_{j}\Pi_{j+1}^k(B)K_{F}^* \E[v_{j}],
\end{equation*}
with the vector $v_j\in Y^n$ given by
\begin{align}
v_j=v_{F,j}+ \lambda_j^\delta R^*v_{G,j}\label{eqn:w}
\end{align}    
where $v_{F,j}$ and $v_{G,j}$ are defined in \eqref{eqn:w_F} and \eqref{eqn:w_G} respectively.
\end{corollary}

\begin{remark}\label{rem:basispertu1}
Without Assumption \ref{ass:sol}(v), for the compact operator $K_F(\cdot)=\sum_{j=1}^\infty \sigma_j\langle\varphi_j,\cdot\rangle\psi_j$, we may design a data-driven approximation of $F$, i.e., $G$, such that $K_G(\cdot)=\sum_{j=1}^{\infty} \tilde{\sigma_j}\langle\tilde\varphi_j,\cdot\rangle\tilde{\psi_j}$ with the singular value $\tilde\sigma_j\leq c_R\sigma_j$ and the orthonormal sequence $\{\tilde\varphi_j\}_{j=1}^{\infty}$ in $X$ satisfying $\sup_{j}\|\tilde\varphi_j-\varphi_j\|\leq\epsilon_G$ for sufficient small $\epsilon_G>0$. In this case, we decompose $K_G$ into two components by $$K_G(\cdot)=\sum_{j=1}^{\infty} \tilde{\sigma_j}\langle\varphi_j,\cdot\rangle\tilde{\psi_j}+\sum_{j=1}^{\infty} \tilde{\sigma_j}\langle\tilde\varphi_j-\varphi_j,\cdot\rangle\tilde{\psi_j}:= K_{G_m}(\cdot)+K_{G_\epsilon}(\cdot),$$
where $K_{G_m}(\cdot)=RK_F(\cdot)$ with $R(\cdot)=\sum_{j=1}^\infty \tilde\sigma_j\sigma_j^{-1}\langle\psi_j,\cdot\rangle\tilde{\psi_j}$ satisfying $\|R\|\leq c_R$, i.e., $K_{G_m}$ satisfies Assumption \ref{ass:sol}(v), and $K_{G_\epsilon}(\cdot)=\sum_{j=1}^{\infty} \tilde{\sigma_j}\langle\tilde\varphi_j-\varphi_j,\cdot\rangle\tilde{\psi_j}$ with $\|K_{G_\epsilon}\|\leq \epsilon_G\|K_G\|$.
Thus, the error $e_k^\delta=x_k^\delta-x^\dag$ satisfies
\begin{equation*}
  \E[e_{k+1}^\delta] = \Pi_1^k(B)e_1^\delta - \sum_{j=1}^{k}\eta_{j}\Pi_{j+1}^k(B) K_{F}^*\E[v_j]-\sum_{j=1}^{k}\eta_{j}\lambda_j^\delta\Pi_{j+1}^k(B) K_{G_\epsilon}^*\E[v_{G,j}],
\end{equation*}
where $v_j$, $v_{F,j}$ and $v_{G,j}$ are defined in \eqref{eqn:w}, \eqref{eqn:w_F} and \eqref{eqn:w_G} respectively.
\end{remark}

The next result gives a useful bound on the mean $\E[v_j]$ under Assumption \ref{ass:sol}.
\begin{lemma}\label{lem:bound-w}
Let Assumption \ref{ass:sol} be fulfilled. Then for $v_j$ defined in \eqref{eqn:w} and 
$e_j^\delta=x_j^\delta-x^\dag$, there holds
\begin{align*}
  \|\E[v_j]\| &\leq \big(\frac{3-\eta_F}{2(1-\eta_F)}c_F + {(c_G\E[\|e_j^\delta\|^2]^\frac12+\tfrac32)}c_Gc_R^2\lambda_j^\delta \big)\E[ \|B_{F}^\frac12 e_j^\delta\|^2]^\frac12\E[ \|e_j^\delta\|^2]^\frac12\\
  &+c_R (c_G\E[ \|e_j^\delta\|^2]^\frac12+1) \lambda_j^\delta C_{max}+\big((c_{F}+c_Gc_R\lambda_j^\delta)\E[\|e_j^\delta\|^2]^\frac12+c_R\lambda_j^\delta+1\big) \delta.
\end{align*}
\end{lemma}
{\begin{proof}
By the triangle inequality and Assumptions \ref{ass:sol}(iii)(v), there holds
\begin{align}\label{eqn:Ev}
  \|\E[v_j]\| &\leq \|\E[v_{F,j}]\|+\lambda_j^\delta \|R^*\E[v_{G,j}]\|\leq \|\E[v_{F,j}]\|+c_R\lambda_j^\delta \|\E[v_{G,j}]\|,
\end{align}
with
\begin{align}
\|\E[v_{F,j}]\|{\leq}& \|\E[( R_{F,x_j^\delta}^{ *}-I )(F(x_j^\delta)-y^\delta)]\|+ \|\E[F(x_j^\delta)-F(x^\dag)-K_{F}(x_j^\delta-x^\dag)]\|+ \|\xi\|:={\rm I}_1+{\rm I}_2+\delta,\label{eqn:EvF}\\
\|\E[v_{G,j}]\|{\leq}& \|\E[( R_{G,x_j^\delta}^{*}-I )(G(x_j^\delta)-y^\delta)]\|+ \|\E[G(x_j^\delta)-G(x^\dag)-K_{G}(x_j^\delta-x^\dag)]\|+\|G(x^\dag)-y^\dag\|+\|\xi\|\nonumber\\
:=&{\rm I}_3+{\rm I}_4+C_{max}+\delta.\label{eqn:EvG}
\end{align}
Now, we bound the terms ${\rm I}_1$--${\rm I}_4$ separately. For the first and third terms ${\rm I}_1$ and ${\rm I}_3$, by the triangle inequality, Assumption \ref{ass:sol}(iv) and Lemma \ref{lem:linear} (under Assumptions \ref{ass:sol}(i)(ii)), there holds
\begin{align*}
{\rm I}_1\leq &\E[\|( R_{F,x_j^\delta}^{ *}-I )(F(x_j^\delta)-y^\delta)\|]\leq c_{F}\E[\|e_j^\delta\|\big(\|F(x_j^\delta)-F(x^\dag)\|+\|y^\dag-y^\delta\|\big)]\\
\leq &c_{F}\E[\|e_j^\delta\|\big(\frac{1}{1-\eta_F}\|K_F e_j^\delta\|+\delta\big)]\leq c_{F}\big(\E[\|e_j^\delta\|]\delta+\frac{1}{1-\eta_F}\E[\|e_j^\delta\|\|K_F e_j^\delta\|]\big),\\
{\rm I}_3\leq &\E[\|( R_{G,x_j^\delta}^{ *}-I )(G(x_j^\delta)-y^\delta)\|]\leq c_{G}\E[\|e_j^\delta\|\|G(x_j^\delta)-y^\delta\|].
\end{align*}
Then, the Cauchy-Schwarz inequality and Lemma \ref{lem:resG_refined} (under Assumptions \ref{ass:sol}(i)(iii)(iv)) imply that
\begin{align*}
{\rm I}_1\leq & c_{F}\E[\|e_j^\delta\|^2]^\frac12 \delta+\frac{c_{F}}{1-\eta_F}\E[\|e_j^\delta\|^2]^\frac12\E[\|K_F e_j^\delta\|^2]^\frac12\\
{\rm I}_3\leq & c_{G}\E[\|e_j^\delta\|^2]^\frac12 \E[\|G(x_j^\delta)-y^\delta\|^2]^\frac12
\leq  c_{G}\E[\|e_j^\delta\|^2]^\frac12 (C_{max}+\delta)+c_{G}(c_G\E[\|e_j^\delta\|^2]^\frac12+1)\E[\|e_j^\delta\|^2]^\frac12\E[\|K_G e_j^\delta\|^2]^\frac12.
\end{align*}
Further, under the Assumption \ref{ass:sol}(v), there holds
\begin{align*}
{\rm I}_3\leq & c_{G}\E[\|e_j^\delta\|^2]^\frac12 (C_{max}+\delta)+c_{G}(c_G\E[\|e_j^\delta\|^2]^\frac12+1)\E[\|e_j^\delta\|^2]^\frac12\E[\|RK_F e_j^\delta\|^2]^\frac12\\
\leq & c_{G}\E[\|e_j^\delta\|^2]^\frac12 (C_{max}+\delta)+c_{G}c_R(c_G\E[\|e_j^\delta\|^2]^\frac12+1)\E[\|e_j^\delta\|^2]^\frac12\E[\|K_F e_j^\delta\|^2]^\frac12.
\end{align*}
For the second and fourth terms ${\rm I}_2$ and ${\rm I}_4$, it follows from the Cauchy-Schwarz inequality and Lemma \ref{lem:R} with $H=F$ and $H=G$ respectively, that
\begin{align*}
{\rm I}_2\leq \E[\|F(x_j^\delta)-F(x^\dag)-K_{F}(x_j^\delta-x^\dag)\|]\leq  \frac{c_F}{2}\E[\|K_F e_j^\delta\|\|e_j^\delta\|]\leq  \frac{c_F}{2}\E[\|K_F e_j^\delta\|^2]^\frac12\E[\|e_j^\delta\|^2]^\frac12,\\
{\rm I}_4\leq \E[\|G(x_j^\delta)-G(x^\dag)-K_{G}(x_j^\delta-x^\dag)\|]\leq  \frac{c_G}{2}\E[\|K_G e_j^\delta\|\|e_j^\delta\|]\leq  \frac{c_G}{2}\E[\|K_G e_j^\delta\|^2]^\frac12\E[\|e_j^\delta\|^2]^\frac12.
\end{align*}
Then, under the Assumption \ref{ass:sol}(v), there holds
\begin{align*}
{\rm I}_4\leq & \frac{c_G}{2}\E[\|R K_F e_j^\delta\|^2]^\frac12 \E[\|e_j^\delta\|^2]^\frac12\leq \frac{c_G c_R}{2}\E[\| K_F e_j^\delta\|^2]^\frac12 \E[\|e_j^\delta\|^2]^\frac12.
\end{align*}
Combining the preceding estimates with the identity $\|K_Fe_j^\delta\|=\|B_F^\frac{1}{2}e_j^\delta\|$ gives the desired bound.
\end{proof}}

\begin{remark}\label{rem:basispertu2}
Without Assumption \ref{ass:sol}(v), by the decomposition $K_G= K_{G_m}+K_{G_\epsilon}=RK_F+K_{G_\epsilon}$ in Remark \ref{rem:basispertu1}, {we can bound $\E[\|K_G e_j^\delta\|^2]^\frac12$ within the upper bounds of the estimates
\begin{align*}
{\rm I}_3:=&\|\E[( R_{G,x_j^\delta}^{*}-I )(G(x_j^\delta)-y^\delta)]\|\leq c_{G}\E[\|e_j^\delta\|^2]^\frac12 (C_{max}+\delta)+c_{G}(c_G\E[\|e_j^\delta\|^2]^\frac12+1)\E[\|e_j^\delta\|^2]^\frac12\E[\|K_G e_j^\delta\|^2]^\frac12,\\
{\rm I}_4:=&\|\E[G(x_j^\delta)-G(x^\dag)-K_{G}(x_j^\delta-x^\dag)]\|\leq\frac{c_G}{2}\E[\|K_G e_j^\delta\|^2]^\frac12\E[\|e_j^\delta\|^2]^\frac12
\end{align*}
provided in the proof of Lemma \ref{lem:bound-w} by
\begin{align*}
\E[\|K_G e_j^\delta\|^2]^\frac12
\leq &\E[\|RK_F e_j^\delta\|^2]^\frac12+\E[\|K_{G_\epsilon}e_j^\delta\|^2]^\frac12\\
\leq &c_R\E[\|K_F e_j^\delta\|^2]^\frac12+\epsilon_G \|K_G\|\E[\|e_j^\delta\|^2]^\frac12
\leq c_R\E[\|B_F^\frac12 e_j^\delta\|^2]^\frac12+L_G \epsilon_G \E[\|e_j^\delta\|^2]^\frac12.
\end{align*}
By the estimate \eqref{eqn:EvG}, we then obtain that
\begin{align*}
\|\E[v_{G,j}]\|
\leq &(c_G\E[\|e_j^\delta\|^2]^\frac12+\tfrac32)c_G\big(\E[\| B_{F}^\frac12e_j^\delta\|^2]^\frac12 \E[\|e_j^\delta\|^2]^\frac12+L_G\epsilon_G \E[\|e_j^\delta\|^2]\big)+(c_{G}\E[\|e_j^\delta\|^2]^\frac12 +1)(C_{max}+\delta).
\end{align*}
Thus, together with the estimates for $\|\E[v_{F,j}]\|$, as given in the proof of Lemma \ref{lem:bound-w}, we derive from \eqref{eqn:Ev} that}
\begin{align*}
\|\E[v_j]\| 
\leq& \big(\frac{3-\eta_F}{2(1-\eta_F)}c_F + {(c_G\E[\|e_j^\delta\|^2]^\frac12+\tfrac32)}c_Gc_R^2\lambda_j^\delta \big)\E[ \|B_{F}^\frac12 e_j^\delta\|^2]^\frac12\E[ \|e_j^\delta\|^2]^\frac12+c_R (c_G\E[ \|e_j^\delta\|^2]^\frac12+1) \lambda_j^\delta C_{max}\\
&+\big((c_{F}+c_Gc_R\lambda_j^\delta)\E[\|e_j^\delta\|^2]^\frac12+c_R\lambda_j^\delta+1\big) \delta+{(c_G\E[\|e_j^\delta\|^2]^\frac12+\tfrac32)}c_GL_G\epsilon_G c_R\lambda_j^\delta\E[\|e_j^\delta\|^2],
\end{align*}
where the last term in the right hand side of the above inequality, i.e., the additional component of the upper bound compared with the estimate in Lemma \ref{lem:bound-w}, tends to $0^+$ as $\epsilon_G\to 0^+$.
In particular, when the data-driven operator $G$ is linear, under Assumptions \ref{ass:sol}(i)--(iv) (where $c_G=0$), with or without Assumption \ref{ass:sol}(v), there hold $\|\E[v_{G,j}]\|\leq C_{max}+\delta$ and 
\begin{align*}
\|\E[v_j]\| \leq& \frac{3-\eta_F}{2(1-\eta_F)}c_F\E[ \|B_{F}^\frac12 e_j^\delta\|^2]^\frac12\E[ \|e_j^\delta\|^2]^\frac12+c_R  \lambda_j^\delta C_{max}+\big(c_{F}\E[\|e_j^\delta\|^2]^\frac12+c_R\lambda_j^\delta+1\big) \delta.
\end{align*}
\end{remark}
Last, we present a bound on the error  $\E[e_k^\delta]$ in the weighted norm. The two cases $s=0$ and
$s=\frac12$ will be employed for deriving convergence rates in Section \ref{ssec:rate}.
\begin{theorem}\label{thm:err-mean}
Let Assumptions \ref{ass:sol} and \ref{ass:source} be fulfilled. Then for the data-driven SGD iterate $x_k^\delta$ and $e_k^\delta=x_{k}^\delta-x^\dag$ and any $s\geq 0$, there holds
\begin{align}
\|B_F^s\E[e_{k+1}^\delta]\|\leq  \phi_0^{s+\nu} \|w\| +\sum_{j=1}^k\eta_j\phi_j ^{\tilde s}\big(C_j\E[ \|B_{F}^\frac12 e_j^\delta\|^2]^\frac12\E[ \|e_j^\delta\|^2]^\frac12+C^G_j\lambda_j^\delta C_{max}+C^F_j \delta\big),\label{eqn:recur-mean}
\end{align}
where $C_j=\frac{3-\eta_F}{2(1-\eta_F)}c_F + {(c_G\E[\|e_j^\delta\|^2]^\frac12+\tfrac32)}c_Gc_R^2\lambda_j^\delta$, $C^G_j=c_R (c_G\E[ \|e_j^\delta\|^2]^\frac12+1) $, and $C^F_j=(c_{F}+c_Gc_R\lambda_j^\delta)\E[\|e_j^\delta\|^2]^\frac12+c_R\lambda_j^\delta+1$.
\end{theorem}
\begin{proof}
By Corollary \ref{cor:recursion-mean} and triangle inequality,
\begin{equation*}
  \|B_F^s\E[e_{k+1}^\delta]\| \leq \|B_F^s\Pi_1^k(B)e^\delta_1\| + \sum_{j=1}^{k}\eta_j\|B_F^s\Pi_{j+1}^k(B)K_{F}^* \E[v_{j}]\|:={\rm I} + \sum_{j=1}^k\eta_j{{\rm I}'_j}.
\end{equation*}
It remains to bound the terms ${\rm I}$ and ${{\rm I}'_j}$. First, under Assumption \ref{ass:sol}(v), the operators $\Pi_j^k(B)$ and $B_F^s$ are commutative for any $j$ and $s$. Together with the source condition in Assumption \ref{ass:source} and the shorthand notation $\phi_j^s$, there holds
\begin{align*}
  {\rm I} = \|B_F^s\Pi_1^k(B)B_F^\nu w\| \leq \|\Pi_1^k(B)B_F^{s+\nu} \|\|w\|=\phi_0^{s+\nu} \|w\|.
\end{align*}
To bound the terms ${{\rm I}'_j}$, we have
\begin{align*}
  {{\rm I}'_j} \leq \|B_F^s\Pi_{j+1}^{k}(B) K_F^*\E[v_{j}]\| \leq \|B_F^{s+\frac12}\Pi_{j+1}^{k}(B)\|\|\E[v_{j}]\|=\phi_j ^{\tilde s}\|\E[v_{j}]\|.
\end{align*}
Then, Lemma \ref{lem:bound-w} and the shorthand notation $\phi_j^s$ complete the proof of the theorem.
\end{proof}
\begin{remark}\label{rem:basispertu3}
Without Assumption \ref{ass:sol}(v), the operators $\Pi_j^k(B)$ and $B_F^s$ may not be commutative. Using the decomposition $K_G= K_{G_m}+K_{G_\epsilon}$ in Remark \ref{rem:basispertu1}, we further decompose $\Pi_j^k(B)$ into $$\Pi_j^k(B)=\prod_{i=j}^k\big(I-\eta_i(B_F+\lambda_i^\delta B_{G_m})-\eta_i\lambda_i^\delta B_{G_\epsilon}\big).$$ 
Under Assumption \ref{ass:stepsize}(ii) which implies that $\|I-\eta_i(B_F+\lambda_i^\delta B_{G_m})\|\leq 1$ and $\eta_i\lambda_i^\delta\leq \eta_0\lambda_0^\delta\leq L_G^{-2}$, for any $x\in X$, we have
\begin{align*}
\|\Pi_j^k(B)x\|=& \|\prod_{i=j}^k\big(I-\eta_i(B_F+\lambda_i^\delta B_{G_m})-\eta_i\lambda_i^\delta B_{G_\epsilon}\big)x\| \\
\leq& \|\prod_{i=j}^k\big(I-\eta_i(B_F+\lambda_i^\delta B_{G_m})\big)x\|+\big((1+\eta_0\lambda_0^\delta\|B_{G_\epsilon}\|)^{k-j+1}-1\big)\|x\|,
\end{align*}
where $\prod_{i=j}^k\big(I-\eta_i(B_F+\lambda_i^\delta B_{G_m})\big):={\prod_{j}^k(B')}$ and $B_F^s$ are commutative, and $$(1+\eta_0\lambda_0^\delta\|B_{G_\epsilon}\|)^{k-j+1}-1
\leq(1+\eta_0\lambda_0^\delta\epsilon_G^2L_G^2)^{k-j+1}-1
\leq(1+\epsilon_G^2)^{k-j+1}-1
\to 0^+, \quad \epsilon_G\to 0^+.$$
We define $\phi_j^{'s}=\|B_F^{s}\Pi_{j+1}^k(B')\|$, following the analysis in the proof of Theorem \ref{thm:err-mean} with the decomposition of {$\prod_{j}^k(B)$} yields that 
\begin{align*}
\|B_F^s\E[e_{k+1}^\delta]\| \leq &\|B_F^s\Pi_1^k(B)e^\delta_1\| + \sum_{j=1}^{k}\eta_j\|B_F^s\Pi_{j+1}^k(B)K_{F}^* \E[v_{j}]\|+\sum_{j=1}^{k}\eta_{j}\lambda_j^\delta\|B_F^s\Pi_{j+1}^k(B) K_{G_\epsilon}^*\E[v_{G,j}]\|\\
\leq &  \phi_0^{'s+\nu} \|w\| +\sum_{j=1}^{k}\eta_j\phi_j ^{'\tilde s}\|\E[v_{j}]\|+\epsilon_G L_G\sum_{j=1}^{k}\eta_{j}\lambda_j^\delta\phi_j ^{'s} \|\E[v_{G,j}]\|+\big((1+\epsilon_G^2)^{k}-1\big)\|B_F\|^s \|e^\delta_1\|\\
&+\sum_{j=1}^{k}\eta_j\big((1+\epsilon_G^2)^{k-j}-1\big)\|B_F\|^s(\|B_F\|^{\frac12} \| \E[v_{j}]\|+\lambda_j^\delta L_G \epsilon_G \|\E[v_{G,j}]\|),
\end{align*}
where the last three terms in the right hand side of the above inequality, i.e., the additional component of the upper bound compared with the estimate in Theorem \ref{thm:err-mean}, tends to $0^+$ as $\epsilon_G\to 0^+$.
\end{remark}
\begin{remark}
Under Assumptions \ref{ass:sol} and \ref{ass:source},
\begin{enumerate}
\item[$\rm(i)$] for linear inverse problems with linear data-driven operator $G$, the recursion \eqref{eqn:recur-mean} can be simplified with $c_F=c_G=0$ to
\begin{align*}
\|B_F^s\E[e_{k+1}^\delta]\|\leq \phi_0^{s+\nu}\|w\| + c_R C_{max}\sum_{j=1}^k\eta_j\phi_j^{\tilde s}\lambda_j^\delta +\sum_{j=1}^k\eta_j\phi_j^{\tilde s}(c_R\lambda_j^\delta+1)\delta,
\end{align*}
where the three terms on the right hand side represent the approximation error, learning error and data error respectively.
\item[$\rm(ii)$] for nonlinear inverse problems with linear data-driven operator $G$, the recursion \eqref{eqn:recur-mean} can be simplified with $c_G=0$ to
\begin{align*}
\|B_F^s\E[e_{k+1}^\delta]\|\leq & \big(\phi_0^{s+\nu} \|w\| +\frac{3-\eta_F}{2(1-\eta_F)}c_F\sum_{j=1}^k\eta_j\phi_j ^{\tilde s}\E[ \|B_{F}^\frac12 e_j^\delta\|^2]^\frac12\E[ \|e_j^\delta\|^2]^\frac12\big)\\
&+c_R C_{max}\sum_{j=1}^k\eta_j\phi_j ^{\tilde s}\lambda_j^\delta+\sum_{j=1}^k\eta_j\phi_j ^{\tilde s}(c_{F}\E[\|e_j^\delta\|^2]^\frac12+c_R\lambda_j^\delta+1) \delta.
\end{align*} 
The estimate of the mean $\E[e_k^\delta]$, which includes an additional stochastic error when compared to that for the linear case in (i) and  \cite{JinLu:2019}, also depends on the variance of the iterate $x_k^\delta$ via the terms
$\E[\|B_F^\frac12 e_j^\delta\|^2]$ and $\E[\|e_j^\delta\|^2]$. 
\end{enumerate}
Compared with the estimate on the mean error of the standard SGD for both linear \cite{JinLu:2019} and nonlinear \cite{JinZhouZou:2020} inverse problems, the data-driven SGD introduces a new error, i.e., the learning error, that related to $C_{max}$, which represents a new phenomena for data-driven algorithms.
\end{remark}

\subsection{Stochastic error}\label{ssec:residual}
Now, we turn to the weighted computational variance $\E[\|B_F^s(x_k^\delta-\E[x_k^\delta])\|^2]=\E[\|B_F^s(e_k^\delta-\E[e_k^\delta])\|^2]$, which arises due to the random choice of
the index $i_k$ at the $k$th data-driven SGD iteration. First, we give an upper bound on the variance in terms of suitable iteration noises
$N_{j,1}$ and $N_{j,2}$ (defined in \eqref{eqn:N} below); see the appendix \ref{app:lem:recursion-var} for the proof.

\begin{lemma}\label{lem:recursion-var}
Let {Assumption \ref{ass:sol}(iv)} be fulfilled. Then for the data-driven SGD iterate $x_k^\delta$ and
$e_j^\delta=x_j^\delta-x^\dag$, there holds
\begin{align*}
\E[\|B_F^s(e_{k+1}^\delta-\E[e_{k+1}^\delta])\|^2] \leq& \big(\sum_{j=1}^k\eta_j\phi_j^{\tilde{s}}(2\E[\|N_{j,1}\|^2]^\frac12+\E[\|N_{j,2}\|^2]^\frac12)\big)\big(\sum_{j=1}^k\eta_j\phi_j^{\tilde{s}}\E[\|N_{j,2}\|^2]^\frac12\big)\\
&+ \sum_{j=1}^k\eta_j^2(\phi_j^{\tilde{s}})^2\E[\|N_{j,1}\|^2],
\end{align*}
with the random variables $N_{j,1}$ and $N_{j,2}$ given by
\begin{equation}\label{eqn:N}
\begin{aligned}
  N_{j,1} &= (K_F+\lambda_j^\delta R^*K_G)e_j^\delta-
   ( K_{F,i_j}+\lambda_j^\delta R^* K_{G,i_j})e_j^\delta\varphi_{i_j},\\
  N_{j,2} & = \E[v_{F,j}]-v_{F,j,i_j}\varphi_{i_j}+ \lambda_j^\delta R^*(\E[v_{G,j}]-v_{G,k,i_j}\varphi_{i_j}),
\end{aligned}
\end{equation}
{where the random variables $v_{F,k}$ and $v_{F,k}$ are defined in \eqref{eqn:w_F} and \eqref{eqn:w_G} respectively, and $v_{F,k,i_k}$ and $v_{G,k,i_k}$ are given by
\begin{align}
v_{F,k,i_k}=&( R_{F,x_k^\delta}^{i_k *}-I )(F_{i_k}(x_k^\delta)-y_{i_k}^\delta)+ \big(F_{i_k}(x_k^\delta)-F_{i_k}(x^\dag)-K_{F,i_k}(x_k^\delta-x^\dag)\big)- \xi_{i_k},\label{eqn:w_Fik}\\
v_{G,k,i_k}=&( R_{G,x_k^\delta}^{i_k *}-I )(G_{i_k}(x_k^\delta)-y_{i_k}^\delta)+ \big(G_{i_k}(x_k^\delta)-G_{i_k}(x^\dag)-K_{G,i_k}(x_k^\delta-x^\dag)\big)+\big(G_{i_k}(x^\dag)-y_{i_k}^\dag\big)- \xi_{i_k},\label{eqn:w_Gik},
\end{align}} 
and $\varphi_i=(0,\ldots,0,
{n^\frac12},0,\ldots,0)$ denotes the $i$th Cartesian coordinate in $\mathbb{R}^n$ scaled by $n^\frac12$.
\end{lemma}

\begin{remark}\label{rem:basispertu4}
Without Assumption \ref{ass:sol}(v), by the decomposition of $K_G$ in Remark \ref{rem:basispertu1}, the random variables $M_{j,1}$ and $M_{j,2}$ in the proof of Lemma \ref{lem:recursion-var} (see the appendix \ref{app:lem:recursion-var}) can be decompose into
\begin{align*}
M_{j,1}=&K_F^* N_{j,1}+\lambda_j^\delta K_{G_\epsilon}^*(K_G e_j^\delta-K_{G,i_j}e_j^\delta\varphi_{i_j} ):=K_F^* N_{j,1}+\lambda_j^\delta K_{G_\epsilon}^*N_{j,1'},\\
M_{j,2}=&K_F^* N_{j,2}+\lambda_j^\delta K_{G_\epsilon}^*(\E[v_{G,j}]-v_{G,k,i_j}\varphi_{i_j}):=K_F^* N_{j,2}+\lambda_j^\delta K_{G_\epsilon}^*N_{j,2'},
\end{align*}
where $$\E[\|M_{j,t}\|^2]^\frac12\leq \E[\|K_F^* N_{j,t}\|^2]^\frac12+\lambda_j^\delta \|K_{G_\epsilon}^*\|\E[\|N_{j,t'}\|^2]^\frac12\leq\E[\|K_F^* N_{j,t}\|^2]^\frac12+\lambda_j^\delta L_G\epsilon_G \E[\|N_{j,t'}\|^2]^\frac12\to\E[\|K_F^* N_{j,t}\|^2]^\frac12$$
as $\epsilon_G\to 0^+$, for any $t=1,2$.
Further, with the decomposition of $\Pi_j^k(B)$ in Remark \ref{rem:basispertu3}, the weighted computational variance is bounded by
\begin{align*}
\E[\|B_F^s(e_{k+1}^\delta-\E[e_{k+1}^\delta])\|^2] \leq& \sum_{j=1}^k \eta_j^2\E[\|B_F^s\Pi_{j+1}^k(B')M_{j,1}\|^2]+ 2\sum_{j=1}^k\sum_{i=1}^j\eta_i\eta_j\E[\|B_F^s\Pi_{i+1}^k(B')M_{i,1}\|\| B_F^s\Pi_{j+1}^k(B')M_{j,2}\|]\\
& + \sum_{j=1}^k\sum_{i=1}^k\eta_i\eta_j\E[\| B_F^s\Pi_{i+1}^k(B')M_{i,2}\|\|B_F^s\Pi_{j+1}^k(B')M_{j,2}\|]+{\rm I}_\epsilon,
\end{align*}
where 
\begin{align*}
{\rm I}_\epsilon=&\|B_F\|^{2s}\sum_{j=1}^k \eta_j^2\epsilon^{(j)}\E[\|M_{j,1}\|^2]+ 2\|B_F\|^{2s}\sum_{j=1}^k\sum_{i=1}^j\eta_i\eta_j\epsilon^{(i)}\epsilon^{(j)}\E[\| M_{i,1}\|\|M_{j,2}\|]\\
& + \|B_F\|^{2s}\sum_{j=1}^k\sum_{i=1}^k\eta_i\eta_j\epsilon^{(i)}\epsilon^{(j)}\E[\|M_{i,2}\|\|M_{j,2}\|],
\end{align*}
with $\epsilon^{(j)}=(1+\epsilon_G^2)^{k-j}-1$.
Sufficient small $\epsilon_G$ gives sufficient small ${\rm I}_\epsilon$.
\end{remark}

The next result bounds the iteration noises $N_{j,1}$ and $N_{j,2}$ under Assumptions \ref{ass:sol} and \ref{ass:stoch}; see the appendix \ref{app:lem:bound-N} for the proof.
\begin{lemma}\label{lem:bound-N}
Let Assumptions \ref{ass:sol} and \ref{ass:stoch} be fulfilled. Then for $N_{j,1}$ and $N_{j,2}$ defined in
\eqref{eqn:N} and $e_j^\delta=x_j^\delta-x^\dag$, there hold
\begin{align*}
  \E[\|N_{j,1}\|^2]^\frac12 \leq& n^\frac12 (1+c_R^2\lambda_j^\delta)\E[\|B_{F}^\frac12 e_j^\delta\|^2]^\frac12,\\
  \E[\|N_{j,2}\|^2]^\frac12 \leq& n^\frac12 \big(\tilde{C}_j\E[\|e_j^\delta\|^2]^\frac{\theta}{2}\E[\|B_F^\frac12 e_j^\delta\|^2]^\frac12+\tilde{C}^G_j\lambda_j^\delta C_{max}+\tilde{C}^F_j\delta\big),\
\end{align*}
where $\tilde{C}_j=\frac{2+\theta-\eta_F}{(1+\theta)(1-\eta_F)}c_F+ {(c_G\E[\|e_j^\delta\|^2]^\frac12+1+\tfrac{1}{1+\theta})}c_Gc_R^2\lambda_j^\delta$, $\tilde{C}^G_j=c_R(c_G\E[\|e_j^\delta\|^2]^\frac{\theta}{2}+1)$, and $\tilde{C}^F_j=(c_F+c_Gc_R\lambda_j^\delta )\E[\|e_j^\delta\|^2]^\frac{\theta}{2}+c_R\lambda_j^\delta+1$.
\end{lemma}

\begin{remark}\label{rem:basispertu5}
Without Assumption \ref{ass:sol}(v), using the estimate for $\E[\|K_G e_j^\delta\|^2]^\frac12$ in Remark \ref{rem:basispertu2}, we may bound $ \E[\|N_{j,1}\|^2]^\frac12$ and $ \E[\|N_{j,2}\|^2]^\frac12$ in Lemma \ref{lem:bound-N} by
\begin{align*}
\E[\|N_{j,1}\|^2]^\frac12 \leq& n^\frac12 (1+c_R^2\lambda_j^\delta)\E[\|B_{F}^\frac12 e_j^\delta\|^2]^\frac12+n^\frac12 c_R\lambda_j^\delta L_G \epsilon_G\E[\|e_j^\delta\|^2]^\frac12,\\
\E[\|N_{j,2}\|^2]^\frac12 \leq& n^\frac12 \big(\tilde{C}_j\E[\|e_j^\delta\|^2]^\frac{\theta}{2}\E[\|B_F^\frac12 e_j^\delta\|^2]^\frac12+\tilde{C}^G_j\lambda_j^\delta C_{max}+\tilde{C}^F_j\delta\big)\\
&+n^\frac12c_R\lambda_j^\delta {c_G(c_G\E[\|e_j^\delta\|^2]^\frac12+1+\tfrac{1}{1+\theta})}L_G \epsilon_G\E[\|e_j^\delta\|^2]^\frac{1+\theta}{2},
\end{align*}
where the additional components of the upper bounds, compared with the estimates in Lemma \ref{lem:bound-N}, tend to $0^+$ as $\epsilon_G\to 0^+$.
\end{remark}

Last, we give a bound on the weighted variance $\E[\|B_F^s(x_k^\delta- \E[x_k^\delta])\|^2]=\E[\|B_F^s(e_k^\delta- \E[e_k^\delta])\|^2]$. This
result will play an important role in deriving error estimates in Section \ref{ssec:rate}.
\begin{theorem}\label{thm:err-var}
Let Assumptions \ref{ass:sol} and \ref{ass:stoch} be fulfilled. Then for the data-driven SGD iterate error $e_{k+1}^\delta=x_{k+1}^\delta-x^\dag$,
there holds for any $s\in[0,\frac12]$,
\begin{align}\label{eqn:recur-var}
\E[\|B_F^s(e_{k+1}^\delta-\E[e_{k+1}^\delta])\|^2] \leq &n\sum_{j=1}^k\eta_j^2(\phi_j^{\tilde{s}})^2(1+c_R^2\lambda_j^\delta)^2\E[\|B_{F}^\frac12 e_j^\delta\|^2]\nonumber\\
&+n\Big(\sum_{j=1}^k\eta_j\phi_j^{\tilde{s}}\big((2+2c_R^2\lambda_j^\delta+\tilde{C}_j\E[\|e_j^\delta\|^2]^\frac{\theta}{2})\E[\|B_F^\frac12 e_j^\delta\|^2]^\frac12+\tilde{C}^G_j\lambda_j^\delta C_{max}+\tilde{C}^F_j\delta\big)\Big)\nonumber\\
&\qquad\times\Big(\sum_{j=1}^k\eta_j\phi_j^{\tilde{s}}\big(\tilde{C}_j\E[\|e_j^\delta\|^2]^\frac{\theta}{2}\E[\|B_F^\frac12 e_j^\delta\|^2]^\frac12+\tilde{C}^G_j\lambda_j^\delta C_{max}+\tilde{C}^F_j\delta\big)\Big),
\end{align}
where $\tilde{C}_j=\frac{2+\theta-\eta_F}{(1+\theta)(1-\eta_F)}c_F+{(c_G\E[\|e_j^\delta\|^2]^\frac12+1+\tfrac{1}{1+\theta})}c_Gc_R^2\lambda_j^\delta $, $\tilde{C}^G_j=c_R(c_G\E[\|e_j^\delta\|^2]^\frac{\theta}{2}+1)$, and $\tilde{C}^F_j=(c_F+c_Gc_R\lambda_j^\delta )\E[\|e_j^\delta\|^2]^\frac{\theta}{2}+c_R\lambda_j^\delta+1$.
\end{theorem}
\begin{proof}
The assertion follows directly from Lemmas \ref{lem:recursion-var} and \ref{lem:bound-N}.
\end{proof}

\begin{remark}
Under Assumptions \ref{ass:sol} and \ref{ass:stoch},
\begin{enumerate}
\item[$\rm(i)$] for linear inverse problems with linear data-driven operator $G$, the constants in the recursion \eqref{eqn:recur-var} can be simplified with $c_F=c_G=0$ to
$$\tilde{C}_j=0, \quad\tilde{C}^G_j=c_R\quad \mbox{and} \quad \tilde{C}^F_j=c_R\lambda_j^\delta+1.$$
\item[$\rm(ii)$] for nonlinear inverse problems with linear data-driven operator $G$, the constants in the recursion \eqref{eqn:recur-var} can be simplified with $c_G=0$ to
$$\tilde{C}_j=\frac{2+\theta-\eta_F}{(1+\theta)(1-\eta_F)}c_F,\quad \tilde{C}^G_j=c_R\quad \mbox{and} \quad \tilde{C}^F_j=c_F\E[\|e_j^\delta\|^2]^\frac{\theta}{2}+c_R\lambda_j^\delta+1.$$
\end{enumerate}
\end{remark}

\subsection{Convergence rates}\label{ssec:rate}
In this subsection, by using the recursions in Theorems \ref{thm:err-mean} and \ref{thm:err-var}, we derive the convergence rates of the data-driven SGD in Algorithm \ref{alg:datasgd} for exact and noisy data in Theorems \ref{thm:err-total-ex} and \ref{thm:err-total} respectively, with polynomially decaying step size and regularization parameter schedules in {Assumption
\ref{ass:stepsize}(ii)} and the source condition in Assumption \ref{ass:source}. The analysis relies heavily on the estimates listed in Appendix \ref{app:estimate}. Without loss of generality, we assume that $\|B_F\|\leq 1$ (which can be easily achieved by properly rescaling the inverse problems), $\eta_0\leq 1$,  $\max(c_R^2,c_R)\lambda_0^\delta\leq1$ and $C_{max}\lambda_0^\delta\leq\|w\|$.

Now, we analyze the case of exact data $y^\dag$, where the constants are clearer in terms of the dependence on various algorithmic parameters. First, we state some estimates on the constants defined in Theorems \ref{thm:err-mean} and \ref{thm:err-var} which is used for deriving the convergence rates; see the appendix \ref{app:lem:constants} for the proof.
\begin{lemma}\label{lem:constants}
Under the assumption $\lambda_j^\delta\leq\lambda_0^\delta\leq \min(c_R^{-2},c_R^{-1})$, for any $j\geq1$, $\theta\in(0,1]$ {and $\eta_F\in[0,1)$}, we can bound the constants $C_j, C_j^G, C_j^F, \tilde{C}_j, \tilde{C}_j^G$ and $\tilde{C}_j^F$ defined in Theorems \ref{thm:err-mean} and \ref{thm:err-var} by
\begin{align*}
&\max(C_j,\tilde{C}_j)\leq {(1+(1-\eta_F)^{-1})c_F+(2+c_G\E[\|e_j^\delta\|^2]^{\frac12})c_G},\quad
\max(C^G_j,\tilde{C}^G_j)\leq c_R (c_G(\E[\|e_j^\delta\|^2]^\frac12+1)+1),\\
\mbox{and } &\max(C^F_j,\tilde{C}^F_j)\leq(c_{F}+c_G)(\E[\|e_j^\delta\|^2]^\frac12+ 1)+2.
\end{align*}
\end{lemma}

\begin{remark}\label{rem:constant}
When the data-driven operator $G$ is linear, we can further simplify the constants above.
\begin{enumerate}
\item[$\rm(i)$] For linear inverse problems with linear data-driven operator $G$, where $c_F=c_G=0$, there hold
\begin{align*}
&\max(C_j,\tilde{C}_j)=0,\quad
\max(C^G_j,\tilde{C}^G_j)\leq c_R,\quad
\mbox{and} \quad \max(C^F_j,\tilde{C}^F_j)\leq 2.
\end{align*}
\item[$\rm(ii)$] For nonlinear inverse problems with linear data-driven operator $G$, where $c_G=0$, there hold
\begin{align*}
&\max(C_j,\tilde{C}_j)\leq {(1+(1-\eta_F)^{-1})c_F},\quad
\max(C^G_j,\tilde{C}^G_j)\leq c_R,\quad
\mbox{and} \quad \max(C^F_j,\tilde{C}^F_j)\leq c_{F}(\E[\|e_j^\delta\|^2]^\frac12+ 1)+2.
\end{align*}
\end{enumerate}
\end{remark}

We derive the convergence rates of the data-driven SGD with exact data by mathematical induction in the following theorem where the upper bounds for both the mean squared error $\E[\|e_k\|^2]$ and the mean squared residual $\E[\|B_F^\frac12e_k\|^2]$ are slightly lower than those achieved in \cite{JinZhouZou:2020}. 
\begin{theorem}{[Convergence rates for exact data]}\label{thm:err-total-ex}
Let Assumptions \ref{ass:sol}, \ref{ass:source}, \ref{ass:stepsize}(ii) and \ref{ass:stoch} be fulfilled with $\|w\|$, $\theta$, $\eta_0$ and $\lambda_0^0$ being sufficiently small, $\|B_F\|\leq1$,  $\max(c_R^2,c_R)\lambda_0^0\leq1$ and $C_{max}\lambda_0^0\leq\|w\|$. Then for the data-driven SGD iterate $x_k$ for the exact data $y^\dag$ defined in \eqref{eqn:datasgd}, the error
$e_k = x_k - x^\dag$ satisfies
\begin{align*}
  \E[\|e_k\|^2] & \leq c^* \|w\|^2 k^{-\min(2\nu(1-\alpha),\alpha)} \quad \mbox{and}\quad   \E[\|B_F^\frac12e_k\|^2] \leq c^*\|w\|^2 k^{-\min((1+2\nu)(1-\alpha),1)},
\end{align*}
where the constant $c^*$ is independent of $k$ but depends on $\alpha$, $\nu$,
$\eta_0$, $\lambda_0^\delta$, $n$, and $\theta$.
\end{theorem}
\begin{proof}
The standard bias-variance decomposition
\begin{equation*}
  \E[\|B_F^se_{k+1}\|^2] = \|B_F^s\E[e_{k+1}]\|^2 + \E[\|B_F^s(e_{k+1}-\E[e_{k+1}])\|^2],
\end{equation*}
and Theorems \ref{thm:err-mean} and \ref{thm:err-var} give the following estimate for any $s\geq0$:
\begin{align}
\E[\|B_F^se_{k+1}\|^2] &\leq \Big(\phi_0^{s+\nu} \|w\| +\sum_{j=1}^k\eta_j\phi_j ^{\tilde s}\big(C_ja_j^\frac12 b_j^\frac12+C^G_j\lambda_j^0 C_{max}\big)\Big)^2+n\sum_{j=1}^k\eta_j^2(\phi_j^{\tilde{s}})^2(1+c_R^2\lambda_j^0)^2b_j\nonumber\\
&+n\Big(\sum_{j=1}^k\eta_j\phi_j^{\tilde{s}}\big((2+2c_R^2\lambda_j^0+\tilde{C}_ja_j^\frac{\theta}{2})b_j^\frac12+\tilde{C}^G_j\lambda_j^0 C_{max}\big)\Big)\Big(\sum_{j=1}^k\eta_j\phi_j^{\tilde{s}}\big(\tilde{C}_ja_j^\frac{\theta}{2}b_j^\frac12+\tilde{C}^G_j\lambda_j^0C_{max}\big)\Big),\label{eqn:recur-err-ex}
\end{align}
where $a_j\equiv \E[\|e_j\|^2]$ and $b_j\equiv \E[\|B_F^\frac12e_j\|^2]$.
By setting $s=0$ and $s=\frac12$
in the recursion \eqref{eqn:recur-err-ex}, we can derive two coupled inequalities
\begin{align}
a_{k+1} \leq &\Big(\phi_0^{\nu} \|w\| +\sum_{j=1}^k\eta_j\phi_j ^{\frac12}\big(C_ja_j^\frac12 b_j^\frac12+C^G_j\lambda_j^0 C_{max}\big)\Big)^2+n\sum_{j=1}^k\eta_j^2(\phi_j^{\frac12})^2(1+c_R^2\lambda_j^0)^2b_j\nonumber\\
&+n\Big(\sum_{j=1}^k\eta_j\phi_j^{\frac12}\big((2+2c_R^2\lambda_j^0+\tilde{C}_ja_j^\frac{\theta}{2})b_j^\frac12+\tilde{C}^G_j\lambda_j^0 C_{max}\big)\Big)\Big(\sum_{j=1}^k\eta_j\phi_j^{\frac12}\big(\tilde{C}_ja_j^\frac{\theta}{2}b_j^\frac12+\tilde{C}^G_j\lambda_j^0C_{max}\big)\Big), \label{eqn:recur-a-ex}\\
b_{k+1} \leq &\Big(\phi_0^{\frac12+\nu} \|w\| +\sum_{j=1}^k\eta_j\phi_j ^1\big(C_ja_j^\frac12 b_j^\frac12+C^G_j\lambda_j^0 C_{max}\big)\Big)^2+n\sum_{j=1}^k\eta_j^2(\phi_j^{1})^2(1+c_R^2\lambda_j^0)^2b_j\nonumber\\
&+n\Big(\sum_{j=1}^k\eta_j\phi_j^{1}\big((2+2c_R^2\lambda_j^0+\tilde{C}_ja_j^\frac{\theta}{2})b_j^\frac12+\tilde{C}^G_j\lambda_j^0 C_{max}\big)\Big)\Big(\sum_{j=1}^k\eta_j\phi_j^{1}\big(\tilde{C}_ja_j^\frac{\theta}{2}b_j^\frac12+\tilde{C}^G_j\lambda_j^0C_{max}\big)\Big). \label{eqn:recur-b-ex}
\end{align}
First we estimate the first term $\phi_0^{s+\nu} \|w\|$ in the first bracket of both $a_{k+1}$ and $b_{k+1}$ where $s=0$ and $\frac12$ respectively.
Under Assumption \ref{ass:stepsize}(ii), for any $\nu\in(0,\frac12)$ and $s\in [0,\frac12]$, Lemma \ref{lem:estimate-B} and the inequality \eqref{eqn:bdd-phi} in Lemma \ref{lem:basicest} directly suggest that
\begin{align}
\phi_0^{s+\nu} \leq& ((s+\nu)e^{-1}(\sum_{i=1}^k\eta_i)^{-1})^{s+\nu}
\leq ((s+\nu)e^{-1}(1-2^{\alpha-1})^{-1}(1-\alpha)\eta_0^{-1}(k+1)^{-(1-\alpha)})^{s+\nu} \nonumber\\
\leq &(\frac{s+\nu}{e})^{s+\nu}(\frac{1-\alpha}{1-2^{\alpha-1}})^{s+\nu}\eta_0^{-(s+\nu)}(k+1)^{-(s+\nu)(1-\alpha)}
\leq2\eta_0^{-(s+\nu)}(k+1)^{-(s+\nu)(1-\alpha)}.\label{eqn:phi-nu}
\end{align}
The last inequality is derived by the facts that the function $(\frac{s+\nu}{e})^{s+\nu}$ is decreasing in $s+\nu$ over the interval $[0,1]$ and the function $\frac{1-\alpha}{1-2^{\alpha-1}}$ is decreasing in $\alpha$ over the interval $[0,1]$.

The rest of the proof is devoted to deriving the following bounds
\begin{equation}\label{eqn:ab}
  a_k  \leq \varrho k^{-\beta} \quad\mbox{and}\quad b_k  \leq \varrho k^{-\gamma},
\end{equation}
with $\beta=\min(2\nu(1-\alpha),\alpha)$, $\gamma=\min((1+2\nu)(1-\alpha),1)$ and $\varrho =c^*\|w\|$ for some constant $c^*$ to be specified below. The proof proceeds by mathematical induction.
For the case $k=1$, the estimates hold trivially for any sufficiently large $c^*$. Now, we assume that the bounds
hold up to the case  $k$, and prove the assertion for the case $k+1$. For any $1\leq j\leq k$, Lemma \ref{lem:constants} and the assertion $a_j\leq \varrho j^{-\beta}\leq \varrho$ directly give that
\begin{align}
\max(C_j,\tilde{C}_j)\leq& {(1+(1-\eta_F)^{-1})c_F+(2+c_G a_j^{\frac12})c_G\leq (1+(1-\eta_F)^{-1})c_F+(2+c_G \varrho^\frac12)c_G}:=c_c,\label{eqn:c_c}\\
\max(C^G_j,\tilde{C}^G_j)\leq& c_R (c_G(a_j^\frac12+1)+1)\leq c_R (c_G(\varrho^\frac12+1)+1):=c_g.\label{eqn:c_g}
\end{align}
With the conditions $\lambda_j^0\leq \lambda_0^0 j^{-\frac12\big(1-\alpha+(1+\theta)\min(2\nu(1-\alpha),\alpha)\big)}=\lambda_0^0 j^{-\frac{\gamma+\theta\beta}{2}}$ in Assumption \ref{ass:stepsize}(ii) and $c_R^2\lambda_0^0\leq1$, it follows from
\eqref{eqn:recur-a-ex}, \eqref{eqn:phi-nu} and the induction hypothesis that
\begin{align*}
a_{k+1} \leq &\Big(2\eta_0^{-\nu}\|w\| (k+1)^{-\nu(1-\alpha)}+\sum_{j=1}^k\eta_j\phi_j^{\frac12}\big(c_c \varrho j^{-\frac{\beta+\gamma}{2}}+c_g\lambda_0^0C_{max}j^{-\frac{\gamma+\theta\beta}{2}} \big)\Big)^2
+4n\sum_{j=1}^k\eta_j^2(\phi_j^{\frac12})^2 \varrho j^{-\gamma}\\
+&n\Big(\sum_{j=1}^k\eta_j\phi_j^{\frac12}\big((4+c_c \varrho^\frac{\theta}{2}j^{-\frac{\theta\beta}{2}})\varrho^\frac{1}{2}j^{-\frac{\gamma}{2}}+c_g\lambda_0^0 C_{max}j^{-\frac{\gamma+\theta\beta}{2}}\big)\Big)\Big(\sum_{j=1}^k\eta_j\phi_j^{\frac12} \big(c_c\varrho^\frac{1+\theta}{2}j^{-\frac{\gamma+\theta\beta}{2}}+c_g\lambda_0^0C_{max}j^{-\frac{\gamma+\theta\beta}{2}}\big)\Big)\\
\leq &\Big(2\eta_0^{-\nu}(k+1)^{-max(0,\nu(1-\alpha)-\frac12 \alpha)} \|w\| (k+1)^{-\frac{\beta}{2}}+\big(c_c\varrho+ c_g\lambda_0^0C_{max}\big)\sum_{j=1}^k \eta_j\phi_j^{\frac12} j^{-\frac{\gamma}{2}}\Big)^2\\
&+4n\varrho\sum_{j=1}^k\eta_j^2(\phi_j^{\frac12})^2 j^{-\gamma}
+n\big((4+c_c \varrho^\frac{\theta}{2})\varrho^\frac{1}{2}+c_g\lambda_0^0 C_{max}\big)\big(c_c\varrho^\frac{1+\theta}{2}+c_g\lambda_0^0C_{max}\big)
\Big(\sum_{j=1}^k\eta_j\phi_j^{\frac12}j^{-\frac{\gamma}{2}} \Big)^2.
\end{align*}
Next we bound the summations on the right hand side. By Proposition \ref{prop:est-ex} in the appendix, we have
\begin{align*}
  \sum_{j=1}^k\eta_j\phi_j^\frac12j^{-\frac{\gamma}{2}}\leq c_1 (k+1)^{-\frac{\beta}{2}} \quad\mbox{and}\quad
  \sum_{j=1}^k\eta_j^2(\phi_j^{\frac12})^2 j^{-\gamma}  \leq c_2(k+1)^{-\beta},
\end{align*}
with $c_1=2^{\frac{\beta}{2}-1}\eta_0^{\frac12}(B(\frac12,\zeta)+2)$ and $c_2=2^{\beta-1}\eta_0((\alpha+\beta)^{-1}+4)$, where $\zeta=1-\alpha-\frac{\gamma}{2}\geq(\frac{1}{2}-\nu)(1-\alpha)>0$ and $B(\cdot,\cdot)$ denotes the Beta function defined in \eqref{eqn:Beta}. Thus, with the notation $c_{\nu,k+1}=2\eta_0^{-\nu}(k+1)^{-max(0,\nu(1-\alpha)-\frac12 \alpha)}$ and the condition $C_{max}\lambda_0^0\leq \|w\|$, we obtain that
\begin{align}\label{eqn:bdd-a}
a_{k+1} 
\leq &\big(c_{\nu,k+1} \|w\| +c_1(c_c\varrho+ c_g\lambda_0^0C_{max})\big)^2(k+1)^{-\beta}
+ 4nc_2 \varrho(k+1)^{-\beta}\nonumber\\
&+nc_1^2\big((4+c_c \varrho^\frac{\theta}{2})\varrho^\frac{1}{2}+c_g\lambda_0^0 C_{max}\big)\big(c_c\varrho^\frac{1+\theta}{2}+c_g\lambda_0^0C_{max}\big)(k+1)^{-\beta}\nonumber\\
\leq &\Big(\big(c_{\nu,k+1} \|w\| +c_1(c_c\varrho+ c_g\|w\|)\big)^2+ 4nc_2 \varrho
+nc_1^2\big((4+c_c \varrho^\frac{\theta}{2})\varrho^\frac{1}{2}+c_g\|w\|\big)\big(c_c\varrho^\frac{1+\theta}{2}+c_g\|w\|\big)\Big)(k+1)^{-\beta}.
\end{align}

Similarly, for the term  $b_k$ and any $\theta\in(0,\frac{1-\alpha}{\beta}-1)$ (where $\frac{1-\alpha}{\beta}-1\geq \frac{1-\alpha}{2\nu(1-\alpha)}-1>0$), it follows from \eqref{eqn:recur-b-ex}, \eqref{eqn:phi-nu}, Lemma \ref{lem:constants}, the assumptions on $\lambda_j^0$ and the induction hypothesis that
\begin{align*}
b_{k+1} \leq &\Big(2\eta_0^{-(\frac12+\nu)}\|w\|(k+1)^{-(\frac12+\nu)(1-\alpha)}+\sum_{j=1}^k\eta_j\phi_j^{1}\big(c_c \varrho j^{-\frac{\beta+\gamma}{2}}+c_g\lambda_0^0C_{max}j^{-\frac{\gamma+\theta\beta}{2}} \big)\Big)^2
+4n\sum_{j=1}^k\eta_j^2(\phi_j^{1})^2 \varrho j^{-\gamma}\\
+&n\Big(\sum_{j=1}^k\eta_j\phi_j^{1}\big((4+c_c \varrho^\frac{\theta}{2}j^{-\frac{\theta\beta}{2}})\varrho^\frac{1}{2}j^{-\frac{\gamma}{2}}+c_g\lambda_0^0 C_{max}j^{-\frac{\gamma+\theta\beta}{2}}\big)\Big)\Big(\sum_{j=1}^k\eta_j\phi_j^{1} \big(c_c\varrho^\frac{1+\theta}{2}j^{-\frac{\gamma+\theta\beta}{2}}+c_g\lambda_0^0C_{max}j^{-\frac{\gamma+\theta\beta}{2}}\big)\Big)\\
\leq &\Big(c_{\frac12+\nu,k+1} \|w\| (k+1)^{-\frac{\gamma}{2}}+\big(c_c\varrho+ c_g\lambda_0^0C_{max}\big)\sum_{j=1}^k \eta_j\phi_j^{1} j^{-\frac{\gamma+\theta\beta}{2}}\Big)^2+4n\varrho\sum_{j=1}^k\eta_j^2(\phi_j^{1})^2 j^{-\gamma}\\
&+n\big((4+c_c \varrho^\frac{\theta}{2})\varrho^\frac{1}{2}+c_g\lambda_0^0 C_{max}\big)\big(c_c\varrho^\frac{1+\theta}{2}+c_g\lambda_0^0C_{max}\big)
\Big(\sum_{j=1}^k\eta_j\phi_j^{1}j^{-\frac{\gamma}{2}} \Big)\Big(\sum_{j=1}^k\eta_j\phi_j^{1}j^{-\frac{\gamma+\theta\beta}{2}} \Big),
\end{align*}
where $c_{\frac12+\nu,k+1}=2\eta_0^{-(\frac12+\nu)}(k+1)^{-max(0,\nu(1-\alpha)-\frac12 \alpha)}$
By Proposition \ref{prop:est-ex} in the appendix, there hold, with 
$\epsilon= \theta \beta$,
\begin{align*}
&\sum_{j=1}^k\eta_j\phi^1_jj^{-\frac{\gamma+\theta\beta}{2}}  \leq c_1' (k+1)^{-\frac{\epsilon}{4}-\frac{\gamma}{2}}, \quad 
\sum_{j=1}^{k} \eta_j^2(\phi_j^1)^2j^{-\gamma} \leq c_2'(k+1)^{-\gamma}\quad\mbox{and} \quad\sum_{j=1}^k\eta_j\phi_j^{1}j^{-\frac{\gamma}{2}}\leq c_3'(k+1)^{\frac{\epsilon}{4}-\frac{\gamma}{2}},
\end{align*}
with $c_1'=2^{\frac{\gamma}{2}-\frac12} \eta_0^{\frac{\theta \beta}{4(1-\alpha)}}\big(B(\frac{\theta \beta}{4(1-\alpha)},\zeta-\frac{\theta\beta}{2})+2\big)$, $c_2'=2^{\gamma+1}\eta_0^{1-\frac{\beta}{1-\alpha}}(\alpha^{-1} +1)$ and $c_3'=2^{\frac{\gamma}{2}-1} \eta_0^{\frac{\theta \beta}{4(1-\alpha)}}(B(\tfrac{\theta \beta}{4(1-\alpha)},\zeta)+2)$. Combining the preceding estimates and the condition $C_{max}\lambda_0^0\leq \|w\|$ yields
\begin{align}\label{eqn:bdd-b}
b_{k+1} \leq \Big(&\big(c_{\frac12+\nu,k+1} \|w\|+c_1'(c_c\varrho+ c_g\|w\|)\big)^2+4nc_2'\varrho\nonumber\\
&+nc_1'c_3'\big((4+c_c \varrho^\frac{\theta}{2})\varrho^\frac{1}{2}+c_g\|w\|\big)(c_c\varrho^\frac{1+\theta}{2}+c_g\|w\|)\Big) (k+1)^{-\gamma}.
\end{align}

In view of the estimates \eqref{eqn:bdd-a} and \eqref{eqn:bdd-b}, upon dividing by
$\varrho$, it suffices to prove the existence of some constant $c^*>0$ such that
\begin{align}
\big(c_{\nu,k+1}c^{*-\frac12} +c_1(c_c\varrho^\frac12+ c_gc^{*-\frac12})\big)^2
+nc_1^2(4+c_c \varrho^\frac{\theta}{2}+c_gc^{*-\frac12})(c_c\varrho^\frac{\theta}{2}+c_gc^{*-\frac12})+ 4nc_2\leq& 1,\label{con:c*}\\
\big(c_{\frac12+\nu,k+1} c^{*-\frac12}+c_1'(c_c\varrho^\frac12+ c_gc^{*-\frac12})\big)^2+nc_1'c_3'\big(4+c_c \varrho^\frac{\theta}{2}+c_gc^{*-\frac12}\big)(c_c\varrho^\frac{\theta}{2}+c_gc^{*-\frac12})+4nc_2'\leq& 1.\label{con:c**}
\end{align}
Note that for fixed $a$, both the functions $B(a,\cdot)$ and $B(\cdot,a)$ are monotonically decreasing, thus the inequalities $\frac{\theta \beta}{4(1-\alpha)}\leq\frac{1-\alpha- \beta}{4(1-\alpha)}\leq \frac12$ (derived from the condition $\theta\in(0,\frac{1-\alpha}{\beta}-1)$), $\theta\beta>0$, $\beta\leq\gamma$ and $\eta_0\leq1$ imply that $c_1\leq c_1'$ and $c_1\leq c_3'$. 
Similarly, the inequalities $0<1-\frac{\beta}{1-\alpha}\leq1$, $(\alpha+\beta)^{-1}\leq 4\alpha^{-1}$, $\beta\leq\gamma$ and $\eta_0\leq1$ suggest that $c_2\leq c_2'$ and $c_{\nu,k+1}\leq c_{\frac12+\nu,k+1}$. 
As a result, conditions \eqref{con:c*} and \eqref{con:c**} can be reduced to condition \eqref{con:c**}.
Since the constants $c_1'$, $c_1'c_3'$ and $c_2'$ are proportional to $\eta_0^{\frac{\theta \beta}{4(1-\alpha)}}$, $\eta_0^{\frac{\theta \beta}{2(1-\alpha)}}$ and $\eta_0^{1-\frac{\beta}{1-\alpha}}$ (where $1-\frac{\beta}{1-\alpha}>\frac{\theta \beta}{2(1-\alpha)}>\frac{\theta \beta}{4(1-\alpha)}>0$) respectively, for sufficiently small $\eta_0$, there hold $c_1'\leq \frac14$, $c_1'c_3'\leq (10n)^{-1}$ and $c_2'\leq (16n)^{-1}$. 
Then, for sufficiently large $c^*\geq 4\max(2c_{\frac12+\nu,k+1},c_g)^2$ (for any $k\in\mathbb{N}$) and sufficiently small {$\varrho$ such that $\varrho^{\frac12} \leq (2c_G^2)^{-1}\big(\sqrt{\tilde{c}_c^2+2c_G^2}-\tilde{c}_c\big)$ and $\varrho^\frac{\theta}{2}\leq (2(\tilde{c}_c+c_G^2))^{-1}$, where $\tilde{c}_c=(1+(1-\eta_F)^{-1})c_F+2c_G$,} with small $\|w\|=\varrho^\frac12 c^{*-\frac12}$, the above conditions hold. This completes the induction step and the proof of the theorem.
\end{proof}

\begin{remark}\label{rem:n}
We consider the condition $c^*\geq 4\max(2c_{\frac12+\nu,k+1},c_g)^2$ in the proof of above theorem, where $$c_{\frac12+\nu,k+1}=2\eta_0^{-(\frac12+\nu)}(k+1)^{-max(0,\nu(1-\alpha)-\frac12 \alpha)}.$$ 
\begin{enumerate}
\item[(i)] When $\alpha\in[\frac{2\nu}{1+2\nu},1)$, which implies that $2\nu(1-\alpha)\leq \alpha$, there holds $c_{\frac12+\nu,k+1}=2\eta_0^{-(\frac12+\nu)}$. In this case, we derive the condition $c^*\geq 4\max(4\eta_0^{-(\frac12+\nu)},c_g)^2$ which indicates that $c^*$ depends on the problem size $n$ due to the dependence of $\eta_0$ on $n$.
\item[(ii)] When $\alpha\in(0,\frac{2\nu}{1+2\nu})$, which implies that $2\nu(1-\alpha)> \alpha$, there holds $c_{\frac12+\nu,k+1}=2\eta_0^{-(\frac12+\nu)}k^{-(\nu(1-\alpha)-\frac12 \alpha)}$. There exists some $k_0\in\mathbb{N}$ such that for any $k+1\geq k_0$, $c_{\frac12+\nu,k+1}\leq\frac12$. 
For the case $k=k_0$, The estimates \eqref{eqn:ab} hold trivially for any sufficiently large $c^*$. In this case, we derive the condition $c^*\geq 4\max(2,c_g)^2$ which indicates that $c^*$ can be independent of the problem size $n$.
\end{enumerate}
\end{remark}

\begin{remark}
From Remark \ref{rem:constant}, for linear inverse problems with linear data-driven operator $G$ where the constants $c_c=\max(C_j,\tilde{C}_j)=0$ and $c_g=\max(C^G_j,\tilde{C}^G_j)\leq c_R$, the conditions \eqref{con:c*} and \eqref{con:c**} can be relaxed to 
\begin{align*}
(c_{\frac12+\nu,k+1} +c_1'c_R)^2 c^{*-1}+nc_1'c_3'(4+c_Rc^{*-\frac12})c_Rc^{*-\frac12}+4nc_2'\leq&1,
\end{align*}
which implies that there are no restrictions on $\varrho$ or $\|w\|$.
\end{remark}

\begin{remark}\label{rem:choose_alpha}
The upper bounds of the mean squared error $\E[\|e_k\|^2]$ and the mean squared residual $\E[\|B_F^\frac12 e_k\|^2]$ for the data-driven SGD with exact data derived in Theorem \ref{thm:err-total-ex} are slightly lower than that obtained in \cite{JinZhouZou:2020}. 
{ With the step size defined in Assumption \ref{ass:stepsize}(ii), the optimal convergence rate (in terms of the iteration) of $\E[\|e_k\|^2]$ is achieved at $\alpha=1-\frac{1}{1+2\nu}$.
When the decay exponent $\alpha$ is chosen close to $0$,} i.e. using an essentially constant step size, the residual $\E[\|B_F^\frac12 e_k\|^2]^\frac12 \leq c^{*\frac12}\|w\| k^{-\frac12-\nu}$, which is identical to that for the Landweber method achieved in \cite[Theorem 3.1]{HankeNeubauerScherzer:1995}.
However, when $\alpha$ approaches $0$, it may add a strict restriction on the upper bound of the error $\E[\|e_k\|^2]^\frac12$, which can not be lower than $c^{*\frac12}\|w\|k^{-\frac{\alpha}{2}}$. 
In addition, $\alpha$ also affects the constant $c^*$ through $c_i$s and $c_i'$s. In particular, it behaves like $\alpha^{-1}$ or $(1-\alpha)^{-1}$ which will explode when $\alpha(1-\alpha)$ approaches $0$. 
Therefore, careful selection of the decay exponent $\alpha$ is of great significance for the algorithm to achieve better convergence rates. This observation is also noted in \cite{JinZhouZou:2020} for the standard SGD.
\end{remark}

Last, we derive convergence rates for noisy data $y^\delta$ in Theorem \ref{thm:err-total}.
\begin{proof}[Proof of Theorem \ref{thm:err-total}]
The proof is similar to that of Theorem \ref{thm:err-total-ex}. Let $a_j\equiv\E[\|e_j^\delta\|^2]$ and
$b_j\equiv\E[\|B^\frac12e_j^\delta\|^2]$. Repeating the argument for Theorem \ref{thm:err-total-ex}, together with the assumption $c_R^2\lambda_j^\delta\leq c_R^2\lambda_0^\delta\leq1$ for any $j\geq1$, leads to the following two coupled recursions:
\begin{align*}
a_{k+1} \leq &\Big(\phi_0^{\nu} \|w\| +\sum_{j=1}^k\eta_j\phi_j ^{\frac12}\big(C_ja_j^\frac12 b_j^\frac12+C^G_j\lambda_j^\delta C_{max}+C^F_j \delta\big)\Big)^2+4n\sum_{j=1}^k\eta_j^2(\phi_j^{\frac12})^2 b_j\\
&+n\Big(\sum_{j=1}^k\eta_j\phi_j^{\frac12}\big((4+\tilde{C}_ja_j^\frac{\theta}{2})b_j^\frac12+\tilde{C}^G_j\lambda_j^\delta C_{max}+\tilde{C}^F_j\delta\big)\Big)\Big(\sum_{j=1}^k\eta_j\phi_j^{\frac12}\big(\tilde{C}_ja_j^\frac{\theta}{2}b_j^\frac12+\tilde{C}^G_j\lambda_j^\delta C_{max}+\tilde{C}^F_j\delta\big)\Big),\\
b_{k+1} \leq &\Big(\phi_0^{\frac12+\nu} \|w\| +\sum_{j=1}^k\eta_j\phi_j ^1\big(C_ja_j^\frac12 b_j^\frac12+C^G_j\lambda_j^\delta C_{max}+C^F_j\delta\big)\Big)^2+4n\sum_{j=1}^k\eta_j^2(\phi_j^{1})^2 b_j\\
&+n\Big(\sum_{j=1}^k\eta_j\phi_j^{1}\big((4+\tilde{C}_ja_j^\frac{\theta}{2})b_j^\frac12+\tilde{C}^G_j\lambda_j^\delta C_{max}+\tilde{C}^F_j\delta\big)\Big)\Big(\sum_{j=1}^k\eta_j\phi_j^{1}\big(\tilde{C}_ja_j^\frac{\theta}{2}b_j^\frac12+\tilde{C}^G_j\lambda_j^\delta C_{max}+\tilde{C}^F_j\delta\big)\Big). 
\end{align*}
Next we prove the following bounds
\begin{equation*}
  a_k  \leq \varrho k^{-\beta}\quad \mbox{and}\quad b_k \leq \varrho k^{-\gamma},
\end{equation*}
for all $k\leq k^*=[(\frac{\delta}{\|w\|})^{-\frac{2}{\gamma+\epsilon}}]$, with $\beta=\min(2\nu(1-\alpha),\alpha)$, $\gamma=\min((1+2\nu)(1-\alpha),1)$, $\epsilon\in(0,2\theta\beta)$ (where $\theta\in(0,\frac{1-\alpha}{\beta}-1)$) and $\varrho =c^*\|w\|^2$ for some constant $c^*$ to be specified below.
Similar to Theorem \ref{thm:err-total-ex}, the proof
proceeds by mathematical induction. The assertion holds trivially for the case $k=1$. Now assume that the bounds hold up to some
$k<k^*$, and we prove the assertion for the case $k+1\leq k^*$. For any $1\leq j\leq k$, Lemma \ref{lem:constants} and the assertion $a_j\leq \varrho j^{-\beta}\leq \varrho$ directly give { the estimates \eqref{eqn:c_c} and \eqref{eqn:c_g} and} that
\begin{align*}
\max(C^F_j,\tilde{C}^F_j)\leq(c_{F}+c_G)(a_j^\frac12+ 1)+2\leq(c_{F}+c_G)(\varrho^\frac12+ 1)+2:=c_f.
\end{align*}
Upon substituting the induction hypothesis and the condition $\lambda_j^\delta=\lambda_0^\delta j^{-\frac{\gamma+\theta\beta}{2}}$ in Assumption \ref{ass:stepsize}(ii), we obtain that
\begin{align*}
a_{k+1} \leq &\Big(2\eta_0^{-\nu}\|w\| (k+1)^{-\nu(1-\alpha)}+\sum_{j=1}^k\eta_j\phi_j^{\frac12}\big(c_c \varrho j^{-\frac{\beta+\gamma}{2}}+c_g\lambda_0^0C_{max}j^{-\frac{\gamma+\theta\beta}{2}}+c_f \delta \big)\Big)^2
+4n\sum_{j=1}^k\eta_j^2(\phi_j^{\frac12})^2 \varrho j^{-\gamma}\\
&+n\Big(\sum_{j=1}^k\eta_j\phi_j^{\frac12}\big((4+c_c \varrho^\frac{\theta}{2}j^{-\frac{\theta\beta}{2}})\varrho^\frac{1}{2}j^{-\frac{\gamma}{2}}+c_g\lambda_0^0 C_{max}j^{-\frac{\gamma+\theta\beta}{2}}+c_f \delta\big)\Big)\\
&\quad\;\;\times\Big(\sum_{j=1}^k\eta_j\phi_j^{\frac12} \big(c_c\varrho^\frac{1+\theta}{2}j^{-\frac{\gamma+\theta\beta}{2}}+c_g\lambda_0^0C_{max}j^{-\frac{\gamma+\theta\beta}{2}}+c_f \delta\big)\Big)\\
\leq &\Big(c_{\nu,k+1} \|w\| (k+1)^{-\frac{\beta}{2}}+\big(c_c\varrho+ c_g\lambda_0^0C_{max}\big)\sum_{j=1}^k \eta_j\phi_j^{\frac12} j^{-\frac{\gamma}{2}}+c_f\sum_{j=1}^k \eta_j\phi_j^{\frac12}\delta\Big)^2
+4n\varrho\sum_{j=1}^k\eta_j^2(\phi_j^{\frac12})^2 j^{-\gamma}\\
&+n\Big(\big((4+c_c \varrho^\frac{\theta}{2})\varrho^\frac{1}{2}+c_g\lambda_0^0 C_{max}\big)\sum_{j=1}^k\eta_j\phi_j^{\frac12}j^{-\frac{\gamma}{2}} +c_f\sum_{j=1}^k \eta_j\phi_j^{\frac12}\delta\Big) \\
&\quad\;\;\times\Big(\big(c_c\varrho^\frac{1+\theta}{2}+c_g\lambda_0^0C_{max}\big)
\sum_{j=1}^k\eta_j\phi_j^{\frac12}j^{-\frac{\gamma}{2}}+c_f\sum_{j=1}^k \eta_j\phi_j^{\frac12}\delta \Big).
\end{align*}
Further, using the estimates in Proposition \ref{prop:est-ex} in the appendix that 
\begin{align*}
\sum_{j=1}^k\eta_j\phi_j^\frac12j^{-\frac{\gamma}{2}}\leq c_1 (k+1)^{-\frac{\beta}{2}}, \quad \sum_{j=1}^k\eta_j^2(\phi_j^{\frac12})^2 j^{-\gamma}  \leq c_2(k+1)^{-\beta},\quad\mbox{and}\quad  
\sum_{j=1}^{k}\eta_j \phi_j^\frac12 \leq c_3 (k+1)^{\frac{1-\alpha}{2}}
\end{align*}
with $c_1=2^{\frac{\beta}{2}-1}\eta_0^{\frac12}(B(\frac12,\zeta)+2)$, $c_2=2^{\beta-1}\eta_0((\alpha+\beta)^{-1}+4)$, and $c_3=2^{-1}\eta_0^\frac12( B(\tfrac12,1-\alpha)+2)$, where $\zeta=1-\alpha-\frac{\gamma}{2}\geq(\frac{1}{2}-\nu)(1-\alpha)>0$, we can bound the right hand side by
\begin{align*}
a_{k+1}\leq &\Big(c_{\nu,k+1} \|w\| (k+1)^{-\frac{\beta}{2}}+c_1\big(c_c\varrho+ c_g\lambda_0^0C_{max}\big) (k+1)^{-\frac{\beta}{2}}+c_3c_f(k+1)^{\frac{1-\alpha}{2}}\delta\Big)^2
+4nc_2\varrho (k+1)^{-\beta}\\
&+n\Big(c_1\big((4+c_c \varrho^\frac{\theta}{2})\varrho^\frac{1}{2}+c_g\lambda_0^0 C_{max}\big)(k+1)^{-\frac{\beta}{2}} +c_3c_f(k+1)^{\frac{1-\alpha}{2}}\delta\Big) \\
&\quad\;\;\times\Big(c_1\big(c_c\varrho^\frac{1+\theta}{2}+c_g\lambda_0^0C_{max}\big)
(k+1)^{-\frac{\beta}{2}}+c_3c_f(k+1)^{\frac{1-\alpha}{2}}\delta \Big).
\end{align*}
Finally, by the choice of $k^*$, for any $k\leq k^*-1$, there holds
\begin{equation}\label{eqn:k-delta}
  (k+1)^{\frac{1-\alpha}{2}}\delta \leq (k+1)^{-\frac{\gamma-1+\alpha+\epsilon}{2}}\|w\|= (k+1)^{-\frac{\beta+\epsilon}{2}}\|w\|,
\end{equation} 
and thus
\begin{align}
a_{k+1}
\leq &\Big(\big(c_{\nu,k+1} \|w\| +c_1(c_c\varrho+ c_g\lambda_0^0C_{max}) +c_3c_f\|w\|\big)^2
+4nc_2\varrho \nonumber\\
&+n\big(c_1((4+c_c \varrho^\frac{\theta}{2})\varrho^\frac{1}{2}+c_g\lambda_0^0 C_{max}) +c_3c_f\|w\|\big)\big(c_1(c_c\varrho^\frac{1+\theta}{2}+c_g\lambda_0^0C_{max})
+c_3c_f\|w\| \big)\Big)(k+1)^{-\beta}.\label{eqn:bdd-a-noise}
\end{align}

For the term $b_{k+1}$, 
it follows from the same steps for bounding $a_{k+1}$ that
\begin{align*}
b_{k+1} \leq &\Big(2\eta_0^{-(\frac12+\nu)}\|w\| (k+1)^{-(\frac12+\nu)(1-\alpha)}+\sum_{j=1}^k\eta_j\phi_j^{1}\big(c_c \varrho j^{-\frac{\beta+\gamma}{2}}+c_g\lambda_0^0C_{max}j^{-\frac{\gamma+\theta\beta}{2}}+c_f \delta \big)\Big)^2
\\
&+4n\sum_{j=1}^k\eta_j^2(\phi_j^{1})^2 \varrho j^{-\gamma}+n\Big(\sum_{j=1}^k\eta_j\phi_j^{1}\big((4+c_c \varrho^\frac{\theta}{2}j^{-\frac{\theta\beta}{2}})\varrho^\frac{1}{2}j^{-\frac{\gamma}{2}}+c_g\lambda_0^0 C_{max}j^{-\frac{\gamma+\theta\beta}{2}}+c_f \delta\big)\Big)\\
&\qquad\qquad\qquad\qquad\qquad\quad\times\Big(\sum_{j=1}^k\eta_j\phi_j^{1} \big(c_c\varrho^\frac{1+\theta}{2}j^{-\frac{\gamma+\theta\beta}{2}}+c_g\lambda_0^0C_{max}j^{-\frac{\gamma+\theta\beta}{2}}+c_f \delta\big)\Big)\\
\leq &\Big(c_{\frac12+\nu,k+1} \|w\| (k+1)^{-\frac{\gamma}{2}}+\big(c_c\varrho+ c_g\lambda_0^0C_{max}\big)\sum_{j=1}^k \eta_j\phi_j^{1} j^{-\frac{\gamma+\theta\beta}{2}}+c_f\sum_{j=1}^k \eta_j\phi_j^{1}\delta\Big)^2\\
&+4n\varrho\sum_{j=1}^k\eta_j^2(\phi_j^{1})^2 j^{-\gamma}+n\Big(\big((4+c_c \varrho^\frac{\theta}{2})\varrho^\frac{1}{2}+c_g\lambda_0^0 C_{max}\big)\sum_{j=1}^k\eta_j\phi_j^{1}j^{-\frac{\gamma}{2}} +c_f\sum_{j=1}^k \eta_j\phi_j^{1}\delta\Big) \\
&\qquad\qquad\qquad\qquad\qquad\quad\times\Big(\big(c_c\varrho^\frac{1+\theta}{2}+c_g\lambda_0^0C_{max}\big)
\sum_{j=1}^k\eta_j\phi_j^{1}j^{-\frac{\gamma+\theta\beta}{2}}+c_f\sum_{j=1}^k \eta_j\phi_j^{1}\delta \Big).
\end{align*}
And further, with the estimates in Proposition \ref{prop:est-ex} in the appendix that
\begin{align*}
&\sum_{j=1}^k\eta_j\phi^1_j j^{-\frac{\gamma+\theta\beta}{2}}
\leq c_1' (k+1)^{-\frac{\epsilon}{4}-\frac{\gamma}{2}},\quad \sum_{j=1}^{k} \eta_j^2(\phi_j^1)^2j^{-\gamma} \leq c_2'(k+1)^{-\gamma},\\
&\sum_{j=1}^k\eta_j\phi^1_j j^{-\frac{\gamma}{2}} 
\leq  c_3'(k+1)^{\frac{\epsilon}{4}-\frac{\gamma}{2}}
\quad \mbox{and}\quad
\sum_{j=1}^{k}\eta_j\phi_j^1  \leq  c_4'(k+1)^{\frac{\epsilon}{4}}, 
\end{align*}
with 
\begin{align*}
c_1'=&2^{\frac{\gamma}{2}-\frac12} \eta_0^{\frac{2\theta \beta-\epsilon}{4(1-\alpha)}}\big(B(\frac{2\theta \beta-\epsilon}{4(1-\alpha)},\zeta-\tfrac{\theta\beta}{2})+2\big), \quad c_2'=2^{\gamma+1}\eta_0^{1-\frac{\beta}{1-\alpha}}(\alpha^{-1} +1),\\
c_3'=&2^{\frac{\gamma}{2}-1} \eta_0^{\frac{\epsilon}{4(1-\alpha)}}(B(\tfrac{\epsilon}{4(1-\alpha)},\zeta)+2), \quad  \mbox{and} \quad c_4'=2^{-1} \eta_0^{\frac{\epsilon}{4(1-\alpha)}}(B(\tfrac{\epsilon}{4(1-\alpha)},1-\alpha)+2),
\end{align*}
we obtain that
\begin{align*}
b_{k+1}\leq &\Big(c_{\frac12+\nu,k+1} \|w\| (k+1)^{-\frac{\gamma}{2}}+c_1'\big(c_c\varrho+ c_g\lambda_0^0C_{max}\big) (k+1)^{-\frac{\epsilon}{4}-\frac{\gamma}{2}}+c_4'c_f(k+1)^{\frac{\epsilon}{4}}\delta\Big)^2\\
&+4nc_2'\varrho(k+1)^{-\gamma}+n\Big(c_3'\big((4+c_c \varrho^\frac{\theta}{2})\varrho^\frac{1}{2}+c_g\lambda_0^0 C_{max}\big)(k+1)^{\frac{\epsilon}{4}-\frac{\gamma}{2}} +c_4'c_f(k+1)^{\frac{\epsilon}{4}}\delta\Big) \\
&\qquad\qquad\qquad\qquad\qquad\times\Big(c_1'\big(c_c\varrho^\frac{1+\theta}{2}+c_g\lambda_0^0C_{max}\big)
 (k+1)^{-\frac{\epsilon}{4}-\frac{\gamma}{2}}+c_4'c_f(k+1)^{\frac{\epsilon}{4}}\delta \Big).
\end{align*}
By the choice of $k^*$, for any $k\leq k^*-1$, there holds $\delta\leq (k+1)^{-\frac{\gamma+\epsilon}{2}}\|w\|$. 
Finally, we can bound $b_{k+1}$ by
\begin{align}
b_{k+1} \leq &\big(c_{\frac12+\nu,k+1} \|w\| +c_1'(c_c\varrho+ c_g\lambda_0^0C_{max})+c_4'c_f\|w\|\big)^2(k+1)^{-\gamma}+4nc_2'\varrho(k+1)^{-\gamma}\nonumber\\
&+n\Big(c_3'\big((4+c_c \varrho^\frac{\theta}{2})\varrho^\frac{1}{2}+c_g\lambda_0^0 C_{max}\big) +c_4'c_f\|w\|\Big)(k+1)^{\frac{\epsilon}{4}-\frac{\gamma}{2}} \nonumber\\
&\qquad\times\Big(c_1'(c_c\varrho^\frac{1+\theta}{2}+c_g\lambda_0^0C_{max})
+c_4'c_f\|w\| \Big)(k+1)^{-\frac{\epsilon}{4}-\frac{\gamma}{2}}\nonumber\\
\leq &\Big(\big(c_{\frac12+\nu,k+1} \|w\| +c_1'(c_c\varrho+ c_g\lambda_0^0C_{max})+c_4'c_f\|w\|\big)^2+4nc_2'\varrho\nonumber\\
&+n\Big(c_3'\big((4+c_c \varrho^\frac{\theta}{2})\varrho^\frac{1}{2}+c_g\lambda_0^0 C_{max}\big) +c_4'c_f\|w\|\Big)\Big(c_1'(c_c\varrho^\frac{1+\theta}{2}+c_g\lambda_0^0C_{max})
+c_4'c_f\|w\| \Big)\Big)(k+1)^{-\gamma}.\label{eqn:bdd-b-noise}
\end{align}
Note that for fixed $a$, the Beta function $B(a,\cdot)$ is monotonically decreasing, thus the inequality $\zeta=1-\alpha-\frac{\gamma}{2}\leq 1-\alpha$ implies that $c_3\leq c_1$ and $c_4'\leq c_3'$.
Then in view of the bounds
\eqref{eqn:bdd-a-noise} and \eqref{eqn:bdd-b-noise}, upon dividing by
$\varrho$, with the condition $C_{max}\lambda_0^\delta\leq\|w\|$, it suffices to prove the existence of some constant $c^*>0$ such that
\begin{align}
\big(c_{\nu,k+1} c^{*-\frac12} +c_1(c_c\varrho^\frac12+ (c_g+c_f)c^{*-\frac12})\big)^2
+nc_1^2c_fc^{*-\frac12}\big(4+c_c \varrho^\frac{\theta}{2}+(c_g+c_f)c^{*-\frac12}\big)&\nonumber\\
+nc_1^2\big(4+c_c \varrho^\frac{\theta}{2}+(c_g+c_f)c^{*-\frac12}\big)\big(c_c\varrho^\frac{\theta}{2}+c_gc^{*-\frac12} \big)+4nc_2&\leq 1,\label{con:c***}\\
\big(c_{\frac12+\nu,k+1} c^{*-\frac12} +c_1'(c_c\varrho^\frac12+ c_gc^{*-\frac12})+c_3'c_fc^{*-\frac12}\big)^2
+nc_3'^2c_fc^{*-\frac12} \big(4+c_c \varrho^\frac{\theta}{2}+(c_g+c_f)c^{*-\frac12}\big)&\nonumber\\
+nc_1'c_3'\big(4+c_c \varrho^\frac{\theta}{2}+(c_g+c_f)c^{*-\frac12}\big)(c_c\varrho^\frac{\theta}{2}+c_gc^{*-\frac12}) +4nc_2'&\leq 1.\label{con:c****}
\end{align}
Following the analysis on the constants in the proof of Theorem \ref{thm:conv-exact}, we have $c_1\leq c_1'$, $c_1\leq c_3'$, $c_2\leq c_2'$ and $c_{\nu,k+1}\leq c_{\frac12+\nu,k+1}$ which imply that conditions \eqref{con:c***} and \eqref{con:c****} can be reduced to condition \eqref{con:c****}.
Since the constants $c_1'$, $c_2'$ and $c_3'$ are proportional to $\eta_0^{\frac{2\theta \beta-\epsilon}{4(1-\alpha)}}$, $\eta_0^{1-\frac{\beta}{1-\alpha}}$ and $\eta_0^{\frac{\epsilon}{4(1-\alpha)}}$ respectively, for sufficiently small $\eta_0$, there hold 
$\max(2c_1',c_3')\leq \min((11n)^{-\frac12},\frac14)$ and $c_2'\leq (16n)^{-1}$.
Then, for sufficiently large $c^*\geq 4\max(2c_{\frac12+\nu,k+1},c_g,c_f)^2$ (for any $k\in\mathbb{N}$) and sufficiently small {$\varrho$ such that $\varrho^{\frac12} \leq (2c_G^2)^{-1}\big(\sqrt{\tilde{c}_c^2+2c_G^2}-\tilde{c}_c\big)$ and $\varrho^\frac{\theta}{2}\leq (2(\tilde{c}_c+c_G^2))^{-1}$, where $\tilde{c}_c=(1+(1-\eta_F)^{-1})c_F+2c_G$,} with small $\|w\|=\varrho^\frac12 c^{*-\frac12}$, the above conditions hold. This completes the induction step and the proof of the theorem.
\end{proof}

\begin{remark}\label{rem:linear_err_total}
From Remark \ref{rem:constant}, for linear inverse problems with linear data-driven operator $G$ where the constants $c_c=\max(C_j,\tilde{C}_j)=0$, $c_g=\max(C^G_j,\tilde{C}^G_j)\leq c_R$ and $c_f=\max(C^F_j,\tilde{C}^F_j)\leq 2$, the condition \eqref{con:c****} can be relaxed to 
\begin{align*}
\big(c_{\frac12+\nu,k+1}  +c_1'c_R+2c_3'\big)^2 c^{*-1}
+nc_3'(c_1'c_R+2c_3') c^{*-\frac12}\big(4+(c_R+2)c^{*-\frac12}\big)+4nc_2'&\leq 1,
\end{align*}
which implies that there are no restrictions on $\varrho$ or $\|w\|$.
\end{remark}

{
\begin{remark}\label{rem:alpha}
By the stopping index $k^*=[(\frac{\delta}{\|w\|})^{-\frac{2}{\min((1+2\nu)(1-\alpha),1)+\epsilon}}]$ provided in Theorem 2.2, when $\epsilon$ is close to 0, the convergence rate (in terms of the noise level) is given by $$\E[\|e_{k^*}^\delta\|^2] \leq c^* \|w\|^{2-2\min(\frac{2\nu}{1+2\nu},\alpha)}\delta^{2\min(\frac{2\nu}{1+2\nu},\alpha)}.$$
To achieve the optimal convergence rate, the decay exponent $\alpha$ of the step size should be greater than $\frac{2\nu}{1+2\nu}$. 
When $\alpha \geq \frac{2\nu}{1+2\nu}$, the impact of the constant $c^*$ on the convergence behavior increases, potentially affecting the convergence rate either positively or negatively, as discussed in Remark 4.12. 
Furthermore, the stopping index $k^*$ will increase as $\alpha$ grows. 
Therefore, to ensure the accuracy and efficiency of the method, a suitable decay exponent $\alpha$ is necessary.
\end{remark}}

\section{Numerical experiments}\label{sec:numer}

In this section, we provide numerical experiments {for both linear and nonlinear inverse problems} to complement the analysis. 

{At the beginning, we shall describe the general idea for constructing the data-driven operator $G$. 
In light of Assumption \ref{ass:sol}(v) for deriving the convergence rate in Section \ref{sec:rate}, we design a neural network with an autoencoder architecture \cite{Goodfellow2016} to approximate the forward operator $F$ by capturing its principal features.
In particular, we consider a class of problems with the forward operator $F=f\circ A$, where $A$ is a compact linear operator and $f$ is a nonlinear operator. 
One can either train a nonlinear autoencoder neural network to simulate the entire operator $F$ or train a linear autoencoder neural network to extract the principal features of $A$, followed by a fully connected or convolutional neural network for approximating $f$.
In this work, we adopted the latter structure, where we can use exact operators to serve the role of well-trained neural networks in order to avoid the influence of the capacity of varying neural networks and the optimization error of training, which are not the focus of our study.
Specifically, we generate several approximate matrices $\tilde{A}$ of $A$ via truncated singular value decomposition, which retain different numbers of principal singular values, to serve as the linear autoencoder architecture. 
We denote the matrix retaining the $N$ principal singular values of $A=\sum_{j=1}^\infty \sigma_j\langle\varphi_j,\cdot\rangle\psi_j$ by $\tilde{A}_{N}=\sum_{j=1}^N \sigma_j\langle\varphi_j,\cdot\rangle\psi_j$, where $\{\varphi_j\}_{j=1}^N$ acts as the encoder and $\{\psi_j\}_{j=1}^N$ as the decoder.
Then, we define the data-driven operator as $G=f\circ \tilde{A}_N$.
} 

\subsection{Linear inverse problems}\label{sec:num_linear}
We {first} focus on the linear inverse problem rather than the nonlinear case discussed in theoretical analysis to observe more transparent dependencies of algorithms on parameters. 
To this end, we employ three examples, denoted by \texttt{phillips} (mildly ill-posed), \texttt{gravity} (moderately ill-posed) and \texttt{shaw} (severely ill-posed) in the public MATLAB package Regutools \cite{P.C.Hansen2007} (available at \url{http://people.compute.dtu.dk/pcha/Regutools/}, last accessed on August 20, 2020). 
These examples are Fredholm/Volterra integral equations of the first kind, discretized using either Galerkin approximation with piecewise constant basis functions or quadrature rules, and all discretized into a linear system of size $n=1000$ with the forward operator $A_{n\times n}$. {The data-driven operator $G$ is chosen as the truncated singular value decomposition $\tilde{A}_{N}$ of $A$, retaining $N$ principal singular values. 
In this setting, Assumption \ref{ass:sol} holds with constants $L_G\leq L_F=\max_i\|a_i\|$, $\eta_F=0$, $c_F=c_G=0$ and $c_R=1$, and Assumption \ref{ass:stoch} holds with any $\theta\in(0,1)$.}
We first normalize the exact solution $x_e$ provided by the package to the reference solution $x^\dag:= x_e/\| x_e\|_{\ell^\infty}$ with $\|\cdot\|_{\ell^\infty}$ denoting the maximum norm of vectors. 
Then, we generate the exact data $y^\dag:=A x^\dag$ and the noisy data $y^\delta:=y^\dag+\delta_0\|y^\dag\|_{\ell^\infty}\xi$, where $\delta_0>0$ represents the relative noise level and each component of $\xi$ follows the standard Gaussian distribution. 

{Now, we shall briefly describe the algorithmic parameters used in the experiment. Both the step sizes and the regularization parameters are chosen as either constant or polynomially decaying schedules, as given in Assumption \ref{ass:stepsize}(ii), which are commonly used in SGD to ensure the convergence. 
The step size is defined as $\eta_k:=\eta_0 k^{-\alpha}$, where the initial step size $\eta_0=c_0/(2\max_i(\|a_i\|^2))$ (with $c_0$ taken from the set $\{1,2\}$) and the decay exponent $\alpha$ is chosen from the set $\{0,0.1,0.3\}$, while the regularization parameter is defined as $\lambda_k^\delta=\lambda_0^\delta k^{-\alpha'}$ where the initial index $\lambda_0^\delta=1$ 
and the decay exponent $\alpha'$ is chosen from the set $\{0,0.1,0.3,0.5\}$.
For the convergence of data-driven SGD (see Theorem \ref{thm:conv-noisy}), the condition $L_F^2 \eta_k <1-\eta_F=1$ and $\sum_{k=1}^\infty \eta_k =\infty$ in Assumption \ref{ass:stepsize}(i) are (almost) satisfied with $c_0=1,2$ and $\alpha=0,0.1,0.3$, while the condition on the regularization parameter fails to hold under our setting.
This inconsistency is due to the limitations of the theoretical analysis and the fact that the convergence behavior is proven for a general data-driven operator, such that $C_{min}\leq\|G(x^\dag)-y^\dag\|\leq C_{max}$, which may not be an appropriate approximation of the forward operator. 
When $C_{max}$ is very small, a constant $\lambda_k^\delta$ (i.e., $\alpha'=0$) can also guarantee the convergence of DSGD.
For deriving certain convergence rates (see Theorem \ref{thm:err-total} and Remark \ref{rem:linear_err_total}), Assumption \ref{ass:stepsize}(ii) (under Assumptions \ref{ass:source} and \ref{ass:stoch}) holds with $c_0=1$ and $\alpha'=0.5$ or $\alpha'\geq (0.5+\nu)(1-\alpha)$, while the smallness condition on $\eta_0$ and $\lambda_0^\delta$ fails to hold.
One may design novel step size and regularization parameter schedules instead of the polynomially decaying type to improve the algorithm; we leave this to future research.

In order to indicate the advantage of the data-driven SGD over the standard SGD, we compare these two methods with the same type of step size schedules. 
The parameter $c_0$ is taken from the set $\{1,2\}$ so that $\eta_0$ satisfies the condition for the convergence of data-driven SGD (see Assumption \ref{ass:stepsize}) and SGD (see \cite[Assumption 2.2]{JinZhouZou:2020}), and is chosen to optimize the average performance of SGD on the specific problem across different noise levels.}
Furthermore, to show the order optimality of these methods with particular step size schedules, we evaluate it against the Landweber method (with a constant step size $1/\|A\|_F^2$) which is proven to be an order optimal regularization method \cite{EnglHankeNeubauer:1996}.
Each method is initialized with $x_1 = 0$, and the maximum number of epochs is fixed at $1$e6 for Landweber method and $1$e5 for (data-driven) SGD, where one epoch refers to $1$ Landweber iteration and $n$ (data-driven) SGD iterations, { with $n=1000$ being the problem size.}  
{The results for \texttt{shaw} with the relative noise level $\delta_0=$1e-3 and the decay exponent $\alpha=0.3$ (where the step size is too small, resulting in the required iterations exceeding 1e5 epochs) although presented in Table \ref{tab:shaw} (and also in Tables \ref{tab:shaw2} and \ref{tab:shaw3}) are not taken into consideration in this work.}
All statistical quantities presented below are computed from 10 independent runs.

\subsubsection{Order optimality of data-driven SGD}\label{sec:num_opt}

In Sections \ref{sec:num_opt} and \ref{sec:num_lambda}, we adopt the data-driven matrix $\tilde{A}_{10}$ (which retains approximately 98\% of the principal components of $A$) for \texttt{phillips} and \texttt{gravity}, and $\tilde{A}_{6}$ (which retains approximately 99\% of the principal components of $A$) for \texttt{shaw}. 
For data-driven SGD (DSGD), SGD, and Landweber method (LM), the stopping indices (counted in epoch) $k_{\rm dsgd}$, $k_{\rm sgd}$ and $k_{\rm lm}$ are taken such that the corresponding mean squared errors $e_{\rm dsgd}=\E[\|x_{k_{\rm dsgd}}^\delta-x^\dag\|^2]$, $e_{\rm sgd}=\E[\|x_{k_{\rm sgd}}^\delta-x^\dag\|^2]$ and $e_{\rm lm}=\E[\|x_{k_{\rm lm}}^\delta-x^\dag\|^2]$ are the smallest along the iteration trajectories. 
{This choice of the stopping index is motivated by the lack of provably order-optimal \textit{a posteriori} stopping rules for DSGD.}
{The numerical results for the three examples -- \texttt{phillips}, \texttt{gravity}, and \texttt{shaw} -- are presented in Tables \ref{tab:phil}, \ref{tab:gravity}, and \ref{tab:shaw}, respectively.}

\begin{table}[htp!]
  \centering\small
  \begin{threeparttable}
  \caption{Comparison of DSGD (with $\tilde{A}_{10}$), SGD and LM for \texttt{phillips}.\label{tab:phil}}
    \begin{tabular}{cccccccccccc}
    \toprule
    \multicolumn{2}{c}{Method}&
    \multicolumn{2}{c}{ DSGD ($c_0=1$, $\alpha'=0$)}&\multicolumn{2}{c}{ SGD ($c_0=1$)}&\multicolumn{2}{c}{LM}\\
    \cmidrule(lr){3-4} \cmidrule(lr){5-6} \cmidrule(lr){7-8}
    $\delta_0$ &$\alpha$ &$e_{\rm dsgd}$&$k_{\rm dsgd}$&$e_{\rm sgd}$&$k_{\rm sgd}$&$e_{\rm lm}$&$k_{\rm lm}$\\
    \midrule
    1e-3&$0$ &1.62e-2 &38.21 &1.87e-2 &39.31 &  1.65e-2 &  5851 \cr
       &$0.1$ &1.50e-2 &85.96 &1.80e-2 &128.37   \cr
       &$0.3$ &1.36e-2 &1517.88 &1.70e-2 &2300.83\cr
    \hline
  5e-3&$0$ &1.29e-1 &10.01 &1.27e-1 &11.58 & 9.28e-2 &  1036  \cr
       &$0.1$ &1.21e-1 &33.65 &1.25e-1 &33.66 \cr
       &$0.3$ &1.09e-1 &340.10 &1.14e-1 &273.10\cr
       \hline
  1e-2&$0$ &3.79e-1 &5.45 &2.40e-1 &2.64  & 1.28e-1 &  249   \cr
       &$0.1$ &2.60e-1 &9.65 &1.98e-1 &9.66 \cr
       &$0.3$ &2.26e-1 &39.49 &1.73e-1 &46.75\cr
       \hline
  5e-2&$0$ &3.54e0 &0.33 &1.54e0 &0.57 & 5.34e-1 &  136 \cr
       &$0.1$ &1.61e0 &1.53 &9.75e-1 &1.84\cr
       &$0.3$  &7.60e-1 &5.07 &5.88e-1 &10.62 \cr
    \bottomrule
    \end{tabular}
    \end{threeparttable}
\end{table}
\begin{table}[htp!]
  \centering\small
  \begin{threeparttable}
  \caption{Comparison of DSGD (with $\tilde{A}_{10}$), SGD and LM for \texttt{gravity}.\label{tab:gravity}}
    \begin{tabular}{cccccccccccc}
    \toprule
    \multicolumn{2}{c}{Method}&
    \multicolumn{2}{c}{ DSGD ($c_0=1$, $\alpha'=0$)}&\multicolumn{2}{c}{ SGD ($c_0=1$)}&\multicolumn{2}{c}{LM}\\
    \cmidrule(lr){3-4} \cmidrule(lr){5-6} \cmidrule(lr){7-8}
    $\delta_0$ &$\alpha$ &$e_{\rm dsgd}$&$k_{\rm dsgd}$&$e_{\rm sgd}$&$k_{\rm sgd}$&$e_{\rm lm}$&$k_{\rm lm}$\\
    \midrule
    1e-3&$0$ & 8.62e-2 &59.21 & 9.81e-2 & 128.37 &  9.39e-2 &  27201 \cr
       &$0.1$ &  8.23e-2 &257.48 &  9.45e-2 & 267.65  \cr
       &$0.3$ &8.36e-2 &5103.99 &9.58e-2 &7429.32\cr
    \hline
  5e-3&$0$ &3.16e-1 &4.75 &3.08e-1 &11.58      & 3.27e-1 &  2515  \cr
       &$0.1$ &2.82e-1 &11.58 &3.24e-1 &18.75    \cr
       &$0.3$ &3.02e-1 &126.54 &3.18e-1 &266.76 \cr
       \hline
  1e-2&$0$ &7.01e-1 &3.99  &6.09e-1 &4.97 & 5.73e-1 &  793   \cr
       &$0.1$ &5.57e-1 &10.59  &5.67e-1 &11.21  \cr
       &$0.3$ &5.64e-1 &49.63  &6.07e-1 &49.66 \cr
       \hline
  5e-2&$0$ &5.41e0 &0.36 &2.83e0 &0.57& 2.07e0 &  149 \cr
       &$0.1$  &3.16e0 &0.57 &2.50e0 &0.57\cr
       &$0.3$   &2.67e0 &1.62 &2.30e0 &5.24\cr
    \bottomrule
    \end{tabular}
    \end{threeparttable}
\end{table}

Observed from the results for all three examples (which have different degrees of
ill-posedness), both DSGD (with the constant regularization parameter $\lambda_k^\delta$, where $\alpha'=0$, which is more relaxed than the assumptions in the theoretical analysis in Theorems \ref{thm:conv-noisy} and \ref{thm:err-total}) and SGD can achieve an accuracy (with much fewer iterations) comparable with that for the optimal Landweber method, which indicates that both DSGD and SGD are optimal methods when combined with suitable step size schedules.  
It is also observed that smaller decay exponents $\alpha$ (with a fixed suitable initial step size) enable DSGD and SGD to achieve comparable accuracy with fewer iterations. 
{However, the accuracy can still be improved by increasing $\alpha$, which shortens the step size and consequently increases the number of iterations. This aligns with the condition for the stopping index, i.e., $k^*=[(\frac{\delta}{\|w\|})^{-\frac{2}{\min((1+2\nu)(1-\alpha),1)+\epsilon}}]$, as given in Theorem \ref{thm:err-total}.}
The best accuracy of numerical results is usually obtained at the intermediate value $\alpha=0.1$ (optimal decay exponent), which is consistent with the analysis in {Remarks \ref{rem:choose_alpha} and \ref{rem:alpha}}. {And the higher the noise level is, the larger the optimal decay exponent $\alpha$ is required.} 
{It is worth noting that, for \texttt{shaw} (where the regularity index $\nu$ is very low), a larger step size schedule (e.g., $c_0=3$, which is outside the range specified by either Assumption \ref{ass:stepsize} or \cite[Assumption 2.2]{JinZhouZou:2020}) also allows DSGD (with decaying step sizes or regularization parameters) to achieve comparable accuracy to LM. 
However, using larger constant step sizes and regularization parameters leads to divergence from the very first few iterations. 
Therefore, the numerical results for this case are not presented in this work.
Similar observations for SGD are given in \cite{JinZhouZou:2021,JinZhouZou:2022}, which generally concludes that the larger the regularity index $\nu$ is, the smaller the value of $c_0$ should be to fully realize the benefit of the smoothness for initial errors and achieve the optimal accuracy. 
In practice, since the regularity index $\nu$ and the relative noise level $\delta_0$ are unknown, we should use a step size that satisfies Assumption \ref{ass:stepsize} to guarantee desirable accuracy of DSGD. 
However, if $\nu$ or $\delta_0$ are known, we can further optimize the efficiency of DSGD (and SGD) by designing better step sizes based on that knowledge.}

\begin{table}[htp!]
  \centering\small
  \begin{threeparttable}
  \caption{Comparison of DSGD (with $\tilde{A}_{6}$), SGD and LM for \texttt{shaw}.\label{tab:shaw}}
    \begin{tabular}{cccccccccccc}
    \toprule
    \multicolumn{2}{c}{Method}&
    \multicolumn{2}{c}{ DSGD ($c_0=2$, $\alpha'=0$)}
    &\multicolumn{2}{c}{ SGD ($c_0=2$)}&\multicolumn{2}{c}{LM}\\
    \cmidrule(lr){3-4} \cmidrule(lr){5-6} \cmidrule(lr){7-8}
    $\delta_0$ &$\alpha$ 
    &$e_{\rm dsgd}$&$k_{\rm dsgd}$&$e_{\rm sgd}$&$k_{\rm sgd}$&$e_{\rm lm}$&$k_{\rm lm}$\\
    \midrule
    1e-3&$0$ &2.82e-1 &2893.54 &2.81e-1 &2649.27  &  2.81e-1 &  760983 \cr
       &$0.1$ &2.81e-1 &12405.07 &2.81e-1 &12405.08  \cr
       &{$0.3$} & {4.50e-1} &{99998.96} &{4.50e-1} &{99999.32}  \cr
    \hline
  5e-3&$0$ &5.33e-1 &58.75 &5.42e-1 &65.07 &5.25e-1 &18588  \cr
       &$0.1$ &5.01e-1 &186.54 &5.28e-1 &203.01\cr
       &$0.3$ &4.98e-1 &4203.87 &5.28e-1 &4693.20 \cr
       \hline
  1e-2&$0$ &6.31e-1 &38.19 &6.90e-1 &41.67 & 6.67e-1 &  12385   \cr
       &$0.1$ &5.60e-1 &106.06 &6.99e-1 &134.69 \cr
       &$0.3$ &5.36e-1 &2190.53 &6.70e-1 &2623.51\cr
       \hline
  5e-2&$0$ &4.38e0 &14.32  &3.22e0 &11.14  &2.91e0 &3392\cr
       &$0.1$ &2.33e0 &30.69 &2.84e0 &30.69 \cr
       &$0.3$ &2.24e0 &397.04 &2.93e0 &394.07\cr
    \bottomrule
    \end{tabular}
    \end{threeparttable}
\end{table}
Now, we compare the results of DSGD with SGD. 
We discuss the results for the examples \texttt{phillips}, \texttt{gravity} and \texttt{shaw} separately.
{For \texttt{phillips} (mildly ill-posed, as shown in Table \ref{tab:phil}) and \texttt{gravity} (moderately ill-posed, as shown in Table \ref{tab:gravity}),} when the noise level $\delta_0$ is relatively low, DSGD can provide higher accuracy with fewer iterations than SGD, which represents a surprising advantage of DSGD over SGD.
However, when the noise level increases, the accuracy of DSGD may be lower than that of SGD for two possible reasons: (i) the regularization term in DSGD introduces additional noisy data errors at each iteration (see Algorithm \ref{alg:datasgd}), which affects the attainable accuracy of DSGD; (ii) the regularization term algorithmically enlarges the step size of the gradient descent concerning all components (which may include relatively high-frequency components) captured by the data-driven matrix $\tilde{A}_{10}$ (see Algorithm \ref{alg:datasgd} and Assumption \ref{ass:sol}(v)), which makes the step size too large to achieve higher accuracy than SGD. 
Moreover, large noise can be mistaken for these relatively high-frequency components, causing damage to the algorithm if not handled properly.
The additional data error can be reduced by using smaller step sizes and regularization parameters (see Section \ref{sec:num_lambda}), and the issues concerning relatively high-frequency components can be avoided by removing these components from the data-driven matrix (see Section \ref{sec:num_G}).

\begin{figure}
  \centering
  \setlength{\tabcolsep}{4pt}
  \begin{tabular}{ccc}
\includegraphics[width=0.31\textwidth,trim={1.5cm 0 0.5cm 0.5cm}]{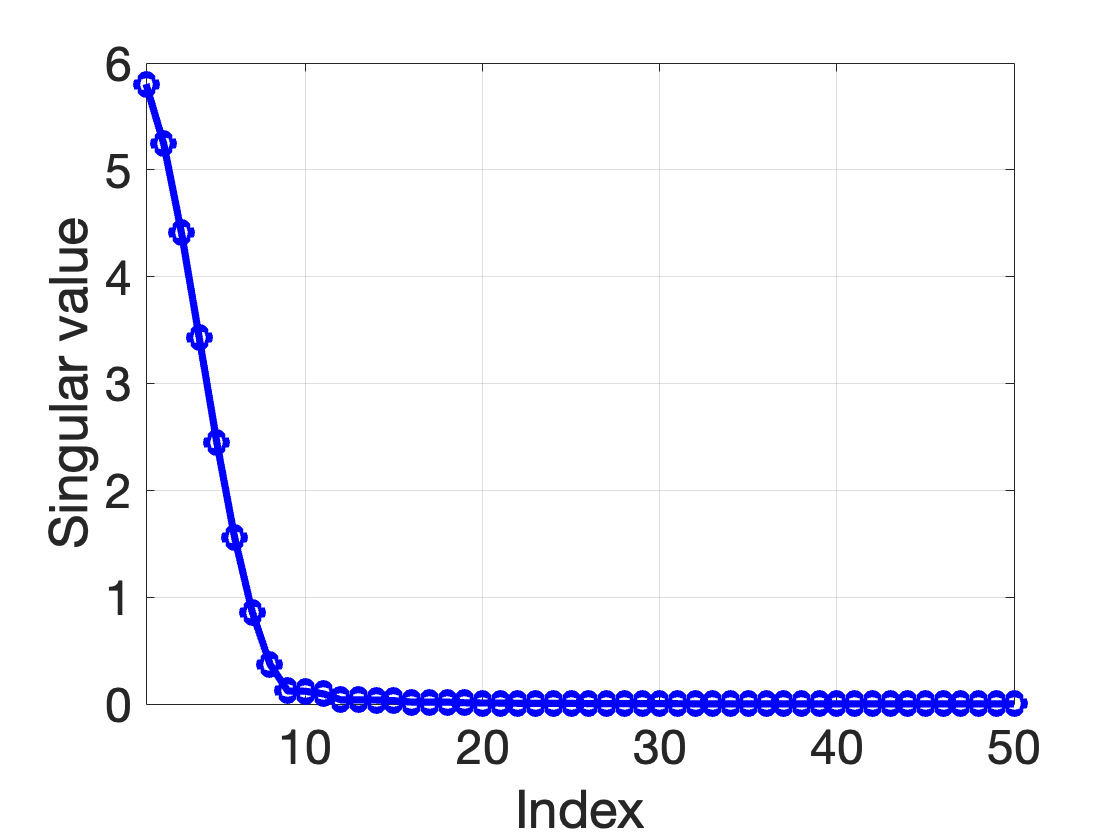}&
\includegraphics[width=0.31\textwidth,trim={1.5cm 0 0.5cm 0.5cm}]{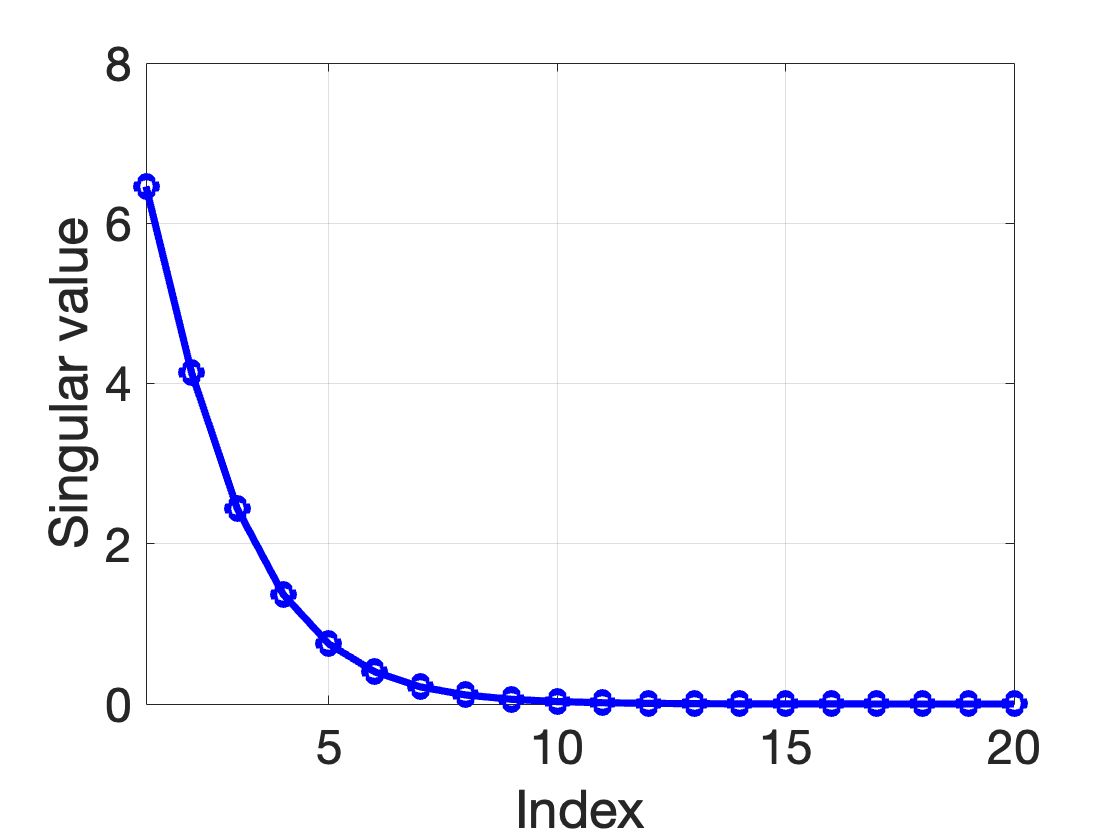}&
\includegraphics[width=0.31\textwidth,trim={1.5cm 0 0.5cm 0.5cm}]{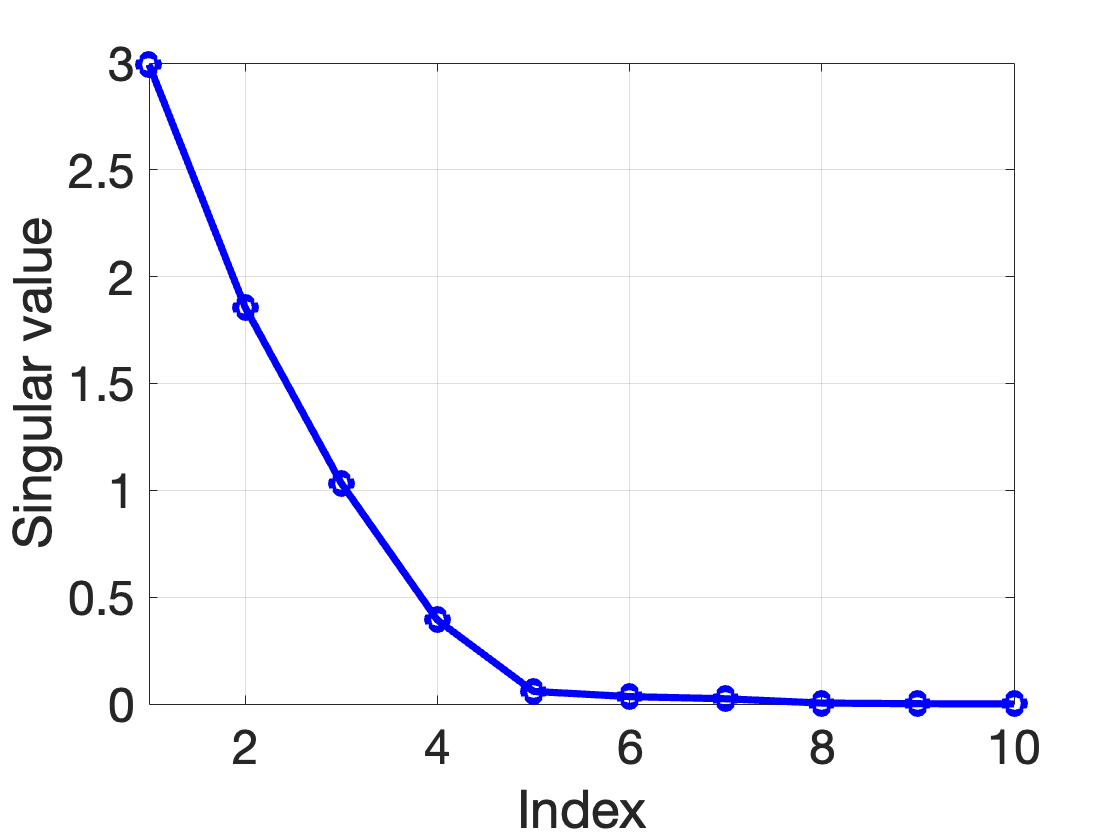}\\
\texttt{phillips}& \texttt{gravity} & \texttt{shaw}
\end{tabular}
\caption{Singular Value Spectrum\label{fig:eig}}
\end{figure}

On the contrary, in the severely ill-posed example \texttt{shaw}, {as shown in Table \ref{tab:shaw},} DSGD provides higher accuracy than SGD for noisier rather than less noisy problems. 
This observation can be explained by the singular value spectrum of $A$ in Figure \ref{fig:eig}. 
The data-driven matrix $\tilde{A}_{6}$ misses several principal components of $A$ that are useful for less noisy problems. 
However, as we discussed before, when the noise is relatively large, these components need to be removed; see Section \ref{sec:num_G} for details.

\subsubsection{Dependence on the regularization parameter}\label{sec:num_lambda}

In order to investigate the impact of the regularization parameter $\lambda_k^\delta=\lambda_0^\delta k^{-\alpha'}$ on DSGD, we present the numerical results of this algorithm with different decay exponent $\alpha'\in\{0,0.1,0.3,0.5\}$ {for the three examples -- \texttt{phillips}, \texttt{gravity}, and \texttt{shaw} -- in Tables \ref{tab:phil2}, \ref{tab:gravity2}, and \ref{tab:shaw2}, respectively.}

\begin{table}[htp!]
  \centering\small
  \begin{threeparttable}
  \caption{Comparison of DSGD ($c_0=1$) with different $\lambda_k^\delta$ and SGD ($c_0=1$) for \texttt{phillips}.\label{tab:phil2}}
    \begin{tabular}{cccccccccccccccc}
    \toprule
    \multicolumn{2}{c}{Method}&
    \multicolumn{2}{c}{DSGD ($\alpha'=0$)}&\multicolumn{2}{c}{ DSGD ($\alpha'=0.1$)}&\multicolumn{2}{c}{ DSGD ($\alpha'=0.3$)}&\multicolumn{2}{c}{ DSGD ($\alpha'=0.5$)}&\multicolumn{2}{c}{SGD}\\
    \cmidrule(lr){3-4} \cmidrule(lr){5-6} \cmidrule(lr){7-8}
    \cmidrule(lr){9-10}
    \cmidrule(lr){11-12}
     $\delta_0$&    $\alpha$ &$e_{\rm dsgd}$&$k_{\rm dsgd}$&$e_{\rm dsgd}$&$k_{\rm dsgd}$&$e_{\rm dsgd}$&$k_{\rm dsgd}$&$e_{\rm dsgd}$&$k_{\rm dsgd}$&$e_{\rm sgd}$&$k_{\rm sgd}$\\
    \midrule
    1e-3&$0$ &1.62e-2 &38.21 &1.67e-2 &38.21 &1.82e-2 &39.31 &1.85e-2 &39.31 &1.87e-2 &39.31
\cr
       &$0.1$ &1.50e-2 &85.96 &1.63e-2 &102.90 &1.75e-2 &128.37 &1.76e-2 &128.37 &1.80e-2 &128.37

\cr
       &$0.3$ &1.36e-2 &1517.88 &1.54e-2 &2008.65 &1.67e-2 &2379.25 &1.68e-2 &2379.25  &1.70e-2 &2300.83 
\cr
    \hline
  5e-3&$0$ &1.29e-1 &10.01 &1.33e-1 &11.58 &1.38e-1 &11.58 &1.38e-1 &11.58 &1.27e-1 &11.58
\cr
       &$0.1$ &1.21e-1 &33.65 &1.28e-1 &33.66 &1.35e-1 &33.66 &1.35e-1 &33.66 &1.25e-1 &33.66
\cr
       &$0.3$ &1.09e-1 &340.10 &1.20e-1 &340.10 &1.25e-1 &340.10 &1.24e-1 &340.10 &1.14e-1 &273.10
\cr
       \hline
  1e-2&$0$ &3.79e-1 &5.45 &3.28e-1 &4.40 &2.90e-1 &4.40 &2.78e-1 &4.40 &2.40e-1 &2.64  
\cr
       &$0.1$  &2.60e-1 &9.65 &2.45e-1 &9.65 &2.31e-1 &9.66 &2.24e-1 &9.66 &1.98e-1 &9.66 
\cr
       &$0.3$ &2.26e-1 &39.49 &2.16e-1 &48.34 &2.05e-1 &48.34 &1.99e-1 &48.34 &1.73e-1 &46.75
\cr       
       \hline
  5e-2&$0$ &3.54e0 &0.33 &2.36e0 &0.44 &1.79e0 &0.44 &1.64e0 &0.57 &1.54e0 &0.57 
\cr
       &$0.1$ &1.61e0 &1.53 &1.30e0 &1.53 &1.09e0 &1.84 &1.04e0 &1.84 &9.75e-1 &1.84
\cr
       &$0.3$ &7.60e-1 &5.07 &6.77e-1 &10.62 &6.39e-1 &10.62 &6.30e-1 &10.62 &5.88e-1 &10.62 
\cr
    \bottomrule
    \end{tabular}
    \end{threeparttable}
\end{table}
\begin{table}[htp!]
  \centering\small
  \begin{threeparttable}
  \caption{Comparison of DSGD ($c_0=1$) with different $\lambda_k^\delta$ and SGD ($c_0=1$) for \texttt{gravity}.\label{tab:gravity2}}
    \begin{tabular}{cccccccccccccccc}
    \toprule
    \multicolumn{2}{c}{Method}&
    \multicolumn{2}{c}{DSGD ($\alpha'=0$)}&\multicolumn{2}{c}{ DSGD ($\alpha'=0.1$)}&\multicolumn{2}{c}{ DSGD ($\alpha'=0.3$)}&\multicolumn{2}{c}{ DSGD ($\alpha'=0.5$)}&\multicolumn{2}{c}{SGD}\\
    \cmidrule(lr){3-4} \cmidrule(lr){5-6} \cmidrule(lr){7-8}
    \cmidrule(lr){9-10}
    \cmidrule(lr){11-12}
     $\delta_0$&    $\alpha$ &$e_{\rm dsgd}$&$k_{\rm dsgd}$&$e_{\rm dsgd}$&$k_{\rm dsgd}$&$e_{\rm dsgd}$&$k_{\rm dsgd}$&$e_{\rm dsgd}$&$k_{\rm dsgd}$&$e_{\rm sgd}$&$k_{\rm sgd}$\\
    \midrule
    1e-3&$0$ &8.62e-2 &59.21 &9.15e-2 &59.21 &9.69e-2 &59.21 &9.78e-2 &128.37 &9.81e-2 &128.37
\cr
       &$0.1$ &8.23e-2 &257.48 &8.86e-2 &267.65 &9.34e-2 &267.65 &9.40e-2 &267.65 &9.45e-2 &267.65
\cr
       &$0.3$ &8.36e-2 &5103.99 &9.15e-2 &6320.27 &9.54e-2 &7429.32 &9.57e-2 &7429.32 &9.58e-2 &7429.32
\cr
    \hline
  5e-3&$0$ &3.16e-1 &4.75 &2.99e-1 &10.59 &2.90e-1 &10.59 &2.92e-1 &10.59 &3.08e-1 &11.58
\cr
       &$0.1$ &2.82e-1 &11.58 &2.96e-1 &13.37 &3.07e-1 &18.75 &3.11e-1 &19.51 &3.24e-1 &18.75
\cr
       &$0.3$ &3.02e-1 &126.54 &3.04e-1 &266.76 &3.06e-1 &266.76 &3.08e-1 &341.71 &3.18e-1 &266.76
\cr
       \hline
  1e-2&$0$ &7.01e-1 &3.99 &6.19e-1 &4.97 &5.94e-1 &4.97 &5.95e-1 &4.97 &6.09e-1 &4.97 
\cr
       &$0.1$ &5.57e-1 &10.59 &5.46e-1 &10.60 &5.52e-1 &11.21 &5.56e-1 &11.21 &5.67e-1 &11.21 
\cr
       &$0.3$ &5.64e-1 &49.63 &5.89e-1 &49.64 &6.20e-1 &49.66 &6.27e-1 &49.84 &6.07e-1 &49.66
\cr       
       \hline
  5e-2&$0$ &5.41e0 &0.36 &4.01e0 &0.57 &3.27e0 &0.57 &3.11e0 &0.57 &2.83e0 &0.57
\cr
       &$0.1$ &3.16e0 &0.57 &2.92e0 &0.57 &2.81e0 &0.57 &2.81e0 &0.57 &2.50e0 &0.57
\cr
       &$0.3$  &2.67e0 &1.62 &2.57e0 &5.24 &2.52e0 &5.24 &2.51e0 &5.24 &2.30e0 &5.24
\cr
    \bottomrule
    \end{tabular}
    \end{threeparttable}
\end{table}
In the examples {\texttt{phillips} (as shown in Table \ref{tab:phil2}) and \texttt{gravity} (as shown in Table \ref{tab:gravity2}),} DSGD, with any regularization parameters, enjoys better accuracy for the problems with relatively low noise levels and stops no later than SGD; while for the cases with high noise levels, DSGD gives lower accuracy than SGD, due to the large step size and data errors, which is also observed in Section \ref{sec:num_opt}. 
For problems with large noise, larger step size decay exponents $\alpha$ or regularization parameters decay exponents $\alpha'$ allow DSGD to improve the attainable accuracy. 
However, in \texttt{shaw} {(as shown in Table \ref{tab:shaw2})}, the observations are opposite to that in \texttt{phillips} or \texttt{gravity}.
For all cases, the behavior of DSGD tends to that of SGD as the regularization parameter becomes smaller and smaller, which makes the data-driven regularization term negligible.

\begin{table}[htp!]
  \centering\small
  \begin{threeparttable}
  \caption{Comparison of DSGD ($c_0=2$) with different $\lambda_k^\delta$ and SGD ($c_0=2$) for \texttt{shaw}.\label{tab:shaw2}}
    \begin{tabular}{cccccccccccccccc}
    \toprule
    \multicolumn{2}{c}{Method}&
    \multicolumn{2}{c}{DSGD ($\alpha'=0$)}&\multicolumn{2}{c}{ DSGD ($\alpha'=0.1$)}&\multicolumn{2}{c}{ DSGD ($\alpha'=0.3$)}&\multicolumn{2}{c}{ DSGD ($\alpha'=0.5$)}&\multicolumn{2}{c}{SGD}\\
    \cmidrule(lr){3-4} \cmidrule(lr){5-6} \cmidrule(lr){7-8}
    \cmidrule(lr){9-10}
    \cmidrule(lr){11-12}
     $\delta_0$&    $\alpha$ &$e_{\rm dsgd}$&$k_{\rm dsgd}$&$e_{\rm dsgd}$&$k_{\rm dsgd}$&$e_{\rm dsgd}$&$k_{\rm dsgd}$&$e_{\rm dsgd}$&$k_{\rm dsgd}$&$e_{\rm sgd}$&$k_{\rm sgd}$\\
    \midrule
    1e-3&$0$ &2.82e-1 &2893.54 &2.81e-1 &2649.27 &2.81e-1 &2649.27 &2.81e-1 &2649.27 &2.81e-1 &2649.27 
\cr
       &$0.1$ &2.81e-1 &12405.07 &2.81e-1 &12405.08 &2.81e-1 &12405.08 &2.81e-1 &12405.08 &2.81e-1 &12405.08 
\cr
       &{$0.3$} & {4.50e-1} &{99998.96}  & {4.50e-1} &{99999.32}& {4.50e-1} &{99999.32}& {4.50e-1} &{99999.32} &{4.50e-1} &{99999.32} 
\cr
    \hline
  5e-3&$0$ &5.33e-1 &58.75 &5.23e-1 &65.07 &5.37e-1 &65.07 &5.41e-1 &65.07 &5.42e-1 &65.07
\cr
       &$0.1$ &5.01e-1 &186.54 &5.08e-1 &195.67 &5.24e-1 &200.87 &5.27e-1 &203.01 &5.28e-1 &203.01
\cr
       &$0.3$ &4.98e-1 &4203.87 &5.10e-1 &4461.03 &5.26e-1 &4693.20 &5.28e-1 &4693.20 &5.28e-1 &4693.20 
\cr
       \hline
  1e-2&$0$ &6.31e-1 &38.19 &6.12e-1 &40.32 &6.70e-1 &41.67 &6.85e-1 &41.67 &6.90e-1 &41.67 
\cr
       &$0.1$ &5.60e-1 &106.06 &6.19e-1 &115.72 &6.84e-1 &128.40 &6.96e-1 &134.69 &6.99e-1 &134.69  
\cr
       &$0.3$ &5.36e-1 &2190.53 &6.00e-1 &2409.55 &6.61e-1 &2623.51 &6.67e-1 &2623.51 &6.70e-1 &2623.51
\cr       
       \hline
  5e-2&$0$ &4.38e0 &14.32 &3.30e0 &11.14 &3.19e0 &11.14 &3.22e0 &11.14 &3.22e0 &11.14 
\cr
       &$0.1$ &2.33e0 &30.69 &2.55e0 &30.69 &2.80e0 &30.69 &2.85e0 &30.69 &2.84e0 &30.69  
\cr
       &$0.3$ &2.24e0 &397.04 &2.65e0 &397.08 &2.91e0 &394.07 &2.94e0 &394.07 &2.93e0 &394.07 
\cr
    \bottomrule
    \end{tabular}
    \end{threeparttable}
\end{table}

There is no doubt that DSGD, with its optimal attainable accuracy and excellent speed, is a better choice than SGD (and LM) when solving relatively mildly ill-posed inverse problems with low noise levels or relatively severely ill-posed inverse problems with high noise levels.
For the mildly or moderately ill-posed problems with high noise levels, DSGD also shows great potential for achieving higher accuracy than SGD when combined with sufficiently small step size and regularization parameter schedules.
However, in practice, we prefer larger step size schedules, which have lower computational complexity, for achieving some desirable (may not be the highest) accuracy. In this case, SGD is more efficient.

\subsubsection{Dependence on the data-driven model}\label{sec:num_G}

Intuitively, when using the exact matrix $A$ as the data-driven matrix in the regularization term, DSGD can be viewed as the standard SGD with a larger step size schedule, which may prevent the algorithm from achieving optimal accuracy.
Meanwhile, from the observation in Sections \ref{sec:num_opt} and \ref{sec:num_lambda}, the regularization term with data-driven matrix $\tilde{A}_{10}$ for \texttt{phillips} and \texttt{gravity}, and $\tilde{A}_{6}$ for \texttt{shaw} improve the accuracy of SGD. 
To study the impact of the proportion of principal features of $A$ captured by the data-driven matrix on DSGD, we present the numerical results of DSGD with the constant regularization parameter and different $\tilde{A}_{N}$ (with $\tilde{A}_{N}$ denoting the matrix retains $N$ principal singular values of $A$) {for the three examples -- \texttt{phillips}, \texttt{gravity}, and \texttt{shaw} -- in Tables \ref{tab:phil3}, \ref{tab:gravity3}, and \ref{tab:shaw3}, respectively.}  

\begin{table}[htp!]
  \centering\small
  \begin{threeparttable}
  \caption{Comparison of DSGD ($c_0=1$, $\alpha'=0$) with different $\tilde{A}_{N}$ and SGD ($c_0=1$) for \texttt{phillips}.\label{tab:phil3}}
    \begin{tabular}{cccccccccccccccc}
    \toprule
    \multicolumn{2}{c}{Method}&
    \multicolumn{2}{c}{DSGD ($N=3$)}&\multicolumn{2}{c}{ DSGD ($N=5$)}&\multicolumn{2}{c}{ DSGD ($N=10$)}&\multicolumn{2}{c}{ DSGD ($N=1000$)}&\multicolumn{2}{c}{SGD}\\
    \cmidrule(lr){3-4} \cmidrule(lr){5-6} \cmidrule(lr){7-8}
    \cmidrule(lr){9-10}
    \cmidrule(lr){11-12}
     $\delta_0$&    $\alpha$ &$e_{\rm dsgd}$&$k_{\rm dsgd}$&$e_{\rm dsgd}$&$k_{\rm dsgd}$&$e_{\rm dsgd}$&$k_{\rm dsgd}$&$e_{\rm dsgd}$&$k_{\rm dsgd}$&$e_{\rm sgd}$&$k_{\rm sgd}$\\
    \midrule
    1e-3&$0$ &5.92e-1 &11.5 &3.41e-2 &49.45 &1.62e-2 &38.21 &2.39e-2 &25.73 &1.87e-2 &39.31
\cr
       &$0.1$ &8.06e-2 &179.17 &1.87e-2 &129.41 &1.50e-2 &85.96 &2.10e-2 &59.19 &1.80e-2 &128.37

\cr
       &$0.3$ &1.76e-2 &2366.35 &1.65e-2 &2313.59 &1.36e-2 &1517.88 &1.83e-2 &942.05 &1.70e-2 &2300.83 
\cr
    \hline
  5e-3&$0$ &6.51e-1 &10.91 &1.40e-1 &10.01 &1.29e-1 &10.01 &1.74e-1 &10.01 &1.27e-1 &11.58
\cr
       &$0.1$ &1.95e-1 &36.11 &1.25e-1 &24.52 &1.21e-1 &33.65 &1.48e-1 &13.59 &1.25e-1 &33.66
\cr
       &$0.3$ &1.12e-1 &241.80 &1.10e-1 &272.96  &1.09e-1 &340.10 &1.32e-1 &184.75 &1.14e-1 &273.10
\cr
       \hline
  1e-2&$0$ &7.35e-1 &1.85 &2.70e-1 &1.78 &3.79e-1 &5.45 &3.91e-1 &1.52 &2.40e-1 &2.64  
\cr
       &$0.1$  &2.80e-1 &10.31 &2.00e-1 &10.2 &2.60e-1 &9.65 &2.87e-1 &4.40 &1.98e-1 &9.66 
\cr
       &$0.3$ &1.72e-1 &46.29 &1.67e-1 &40.18 &2.26e-1 &39.49 &2.27e-1 &35.37 &1.73e-1 &46.75
\cr       
\hline
  5e-2&$0$ &2.38e0 &0.44 &2.40e0 &1.52 &3.54e0 &0.33 &3.58e0 &0.33 &1.54e0 &0.57 
\cr
       &$0.1$ &1.23e0 &1.53 &1.19e0 &1.52 &1.61e0 &1.53 &1.68e0 &1.53 &9.75e-1 &1.84
\cr
       &$0.3$ &5.98e-1 &10.62 &5.83e-1 &10.62 &7.60e-1 &5.07 &7.69e-1 &4.40 &5.88e-1 &10.62 
\cr
    \bottomrule
    \end{tabular}
    \end{threeparttable}
\end{table}

In {\texttt{phillips} (as shown in Table \ref{tab:phil3}) and \texttt{gravity} (as shown in Table \ref{tab:gravity3}),} 
the data-driven matrices $\tilde{A}_{3}$, $\tilde{A}_{5}$, $\tilde{A}_{10}$ and $\tilde{A}_{1000}$ retain approximately 50\%, 90\%, 98\% and 100\% of the principal components of $A$ respectively. 
Clearly, DSGD combined with suitable step size schedules and parameters $N$ has the capability to provide better accuracy than SGD. 
In general, the higher the noise level is, the smaller the value of $N$ needs to be taken, which means that fewer and lower-frequency components of $A$ will be captured by the data-driven matrix. Otherwise, large noise may be incorrectly identified as relatively high-frequency components, which can prevent the iteration from achieving optimal accuracy.
Similar behavior for DSGD with different $N$ is observed from the results of \texttt{shaw} {(as shown in Table \ref{tab:shaw3})}, where the data-driven matrices $\tilde{A}_{3}$, $\tilde{A}_{4}$, $\tilde{A}_{6}$ and $\tilde{A}_{1000}$ retain approximately 90\%, 98\%, 99\% and 100\% of the principal components of $A$ respectively. 
The difference is that, when the noise level is sufficiently low, SGD with a larger step size schedule (i.e., DSGD with $N=n=1000$) is more efficient than DSGD as smaller $N$ will not improve the accuracy but will increase the computational complexity. 

\begin{table}[htp!]
  \centering\small
  \begin{threeparttable}
  \caption{Comparison of DSGD ($c_0=1$, $\alpha'=0$) with different $\tilde{A}_{N}$ and SGD ($c_0=1$) for \texttt{gravity}.\label{tab:gravity3}}
    \begin{tabular}{cccccccccccccccc}
    \toprule
    \multicolumn{2}{c}{Method}&
    \multicolumn{2}{c}{DSGD ($N=3$)}&\multicolumn{2}{c}{ DSGD ($N=5$)}&\multicolumn{2}{c}{ DSGD ($N=10$)}&\multicolumn{2}{c}{ DSGD ($N=1000$)}&\multicolumn{2}{c}{SGD}\\
    \cmidrule(lr){3-4} \cmidrule(lr){5-6} \cmidrule(lr){7-8}
    \cmidrule(lr){9-10}
    \cmidrule(lr){11-12}
     $\delta_0$&    $\alpha$ &$e_{\rm dsgd}$&$k_{\rm dsgd}$&$e_{\rm dsgd}$&$k_{\rm dsgd}$&$e_{\rm dsgd}$&$k_{\rm dsgd}$&$e_{\rm dsgd}$&$k_{\rm dsgd}$&$e_{\rm sgd}$&$k_{\rm sgd}$\\
    \midrule
    1e-3&$0$ &2.00e-1 &35.89 &9.76e-2 &128.37 &8.62e-2 &59.21 &9.71e-2 &39.70 &9.81e-2 &128.37
\cr
       &$0.1$ &1.03e-1 &267.75 &9.39e-2 &261.92 &8.23e-2 &257.48 &9.79e-2 &157.32 &9.45e-2 &267.65
\cr
       &$0.3$ &9.60e-2 &7451.99 &9.54e-2 &7704.13 &8.36e-2 &5103.99 &9.80e-2 &2614.97 &9.58e-2 &7429.32
\cr
    \hline
  5e-3&$0$ &4.22e-1 &11.27 &3.19e-1 &11.53 &3.16e-1 &4.75 &3.27e-1 &4.75 &3.08e-1 &11.58
\cr
       &$0.1$ &3.37e-1 &18.97 &3.16e-1 &18.71 &2.82e-1 &11.58 &2.92e-1 &11.58 &3.24e-1 &18.75
\cr
       &$0.3$ &3.21e-1 &198.21 &3.13e-1 &262.69 &3.02e-1 &126.54 &3.12e-1 &126.54 &3.18e-1 &266.76
\cr
       \hline
  1e-2&$0$ &7.45e-1 &2.54 &6.58e-1 &4.97 &7.01e-1 &3.99 &7.16e-1 &1.75 &6.09e-1 &4.97 
\cr
       &$0.1$ &5.74e-1 &11.22 &5.52e-1 &10.59 &5.57e-1 &10.59 &5.90e-1 &10.59 &5.67e-1 &11.21 
\cr
       &$0.3$ &5.86e-1 &55.28 &5.83e-1 &49.63 &5.64e-1 &49.63 &5.75e-1 &49.63 &6.07e-1 &49.66
\cr       
       \hline
  5e-2&$0$ &3.72e0 &0.36 &4.47e0 &0.57 &5.41e0 &0.36 &5.41e0 &0.36 &2.83e0 &0.57
\cr
       &$0.1$ &2.65e0 &0.57 &2.65e0 &0.57 &3.16e0 &0.57 &3.17e0 &0.57 &2.50e0 &0.57
\cr
       &$0.3$ &2.29e0 &5.24 &2.27e0 &3.98 &2.67e0 &1.62 &2.68e0 &2.62 &2.30e0 &5.24
\cr
    \bottomrule
    \end{tabular}
    \end{threeparttable}
\end{table}

\begin{table}[htp!]
  \centering\small
  \begin{threeparttable}
  \caption{Comparison of DSGD ($c_0=2$, $\alpha'=0$) with different $\tilde{A}_{N}$ and SGD ($c_0=2$) for \texttt{shaw}.\label{tab:shaw3}}
    \begin{tabular}{cccccccccccccccc}
    \toprule
    \multicolumn{2}{c}{Method}&
    \multicolumn{2}{c}{DSGD ($N=3$)}&\multicolumn{2}{c}{ DSGD ($N=4$)}&\multicolumn{2}{c}{ DSGD ($N=6$)}&\multicolumn{2}{c}{ DSGD ($N=1000$)}&\multicolumn{2}{c}{SGD}\\
    \cmidrule(lr){3-4} \cmidrule(lr){5-6} \cmidrule(lr){7-8}
    \cmidrule(lr){9-10}
    \cmidrule(lr){11-12}
     $\delta_0$&    $\alpha$ &$e_{\rm dsgd}$&$k_{\rm dsgd}$&$e_{\rm dsgd}$&$k_{\rm dsgd}$&$e_{\rm dsgd}$&$k_{\rm dsgd}$&$e_{\rm dsgd}$&$k_{\rm dsgd}$&$e_{\rm sgd}$&$k_{\rm sgd}$\\
    \midrule
    1e-3&$0$ &3.38e-1 &2487.09 &2.82e-1 &2894.04 &2.82e-1 &2893.54 &2.80e-1 &1345.96 &2.81e-1 &2649.27 
\cr
       &$0.1$ &2.85e-1 &13158.28 &2.81e-1 &12405.08 &2.81e-1 &12405.07 &2.80e-1 &5917.86 &2.81e-1 &12405.08  
\cr
       &{$0.3$} & {4.50e-1} &{99996.42} & {4.50e-1} &{99999.45}& {4.50e-1} &{99998.96}& {3.69e-1} &{99999.32} &{4.50e-1} &{99999.32}
\cr
    \hline
  5e-3&$0$ &6.21e-1 &65.17 &5.65e-1 &66.02 &5.33e-1 &58.75 &5.64e-1 &30.65 &5.42e-1 &65.07
\cr
       &$0.1$ &5.33e-1 &204.33 &5.29e-1 &200.82 &5.01e-1 &186.54 &5.39e-1 &96.71 &5.28e-1 &203.01
\cr
       &$0.3$  &5.28e-1 &4708.58 &5.28e-1 &4692.52 &4.98e-1 &4203.87 &5.29e-1 &1770.67 &5.28e-1 &4693.20 
\cr
       \hline
  1e-2&$0$ &8.17e-1 &40.67 &7.66e-1 &42.29 &6.31e-1 &38.19 &7.96e-1 &24.55 &6.90e-1 &41.67 
\cr
       &$0.1$ &7.10e-1 &130.58 &7.03e-1 &136.26 &5.60e-1 &106.06 &7.17e-1 &58.67 &6.99e-1 &134.69  
\cr
       &$0.3$ &6.68e-1 &2613.80 &6.69e-1 &2623.02 &5.36e-1 &2190.53 &6.74e-1 &979.01 &6.70e-1 &2623.51
\cr  
\hline
  5e-2&$0$ &4.66e0 &8.30 &4.43e0 &10.60 &4.38e0 &14.32 &4.80e0 &6.02 &3.22e0 &11.14 
\cr
       &$0.1$ &3.01e0 &30.24 &3.00e0 &30.69 &2.33e0 &30.69 &3.04e0 &16.86 &2.84e0 &30.69  
\cr
       &$0.3$ &2.93e0 &396.91 &2.92e0 &397.04 &2.24e0 &397.04 &3.00e0 &164.06 &2.93e0 &394.07 
\cr
    \bottomrule
    \end{tabular}   \end{threeparttable}
\end{table}

Based on these observations, we arrive at a similar conclusion to the discussion in section \ref{sec:num_lambda}:
DSGD, when combined with appropriate step sizes and data-driven matrices, is more efficient than SGD (and LM) in solving relatively mildly ill-posed inverse problems with any noise level or relatively severely ill-posed inverse problems with high noise levels.
However, SGD is more efficient when solving inverse problems that are less noisy and severely ill-posed. 
{In practice, since the levels of noise and ill-posedness are unknown, DSGD is an excellent choice when combined with a suitable data-driven operator, as it performs better than SGD in most cases and does not compromise the accuracy of SGD in other cases.}

{\subsection{Nonlinear inverse problems}\label{sec:num_nonlinear}

In this section, we consider two simple nonlinear inverse problems derived from the linear problems \texttt{phillips} and \texttt{shaw} by defining $y^\dag=F(x^\dag):=(Ax^\dag)^2$, where $A$ and $x^\dag$ are given in Section \ref{sec:num_linear}, and $(\cdot)^2$ is applied component-wise. 
These two problems are named \texttt{squared-phillips} and \texttt{squared-shaw}, respectively. 
The Jacobian matrix of $F$ at the point $x$ is given by $F'(x)=2{\rm diag}(Ax)A$, where ${\rm diag}(Ax)$ is the diagonal matrix with the components of $Ax$ on the diagonal, and the gradient of $F_i$ at $x$ is given by $F_i'(x)=2\langle a_i, x\rangle a_i^t$, where $a_i^t$ is the $i$th row of $A$. 
We define the data-driven operator $G=(\tilde{A}_N)^2$ and adopt $\tilde{A}_{10}$ for \texttt{squared-phillips} and $\tilde{A}_{6}$ for \texttt{squared-shaw}, where $N$ is selected as the best choice for the corresponding linear problems, as observed from Tables \ref{tab:phil3} and \ref{tab:shaw3}, respectively.
The data-driven SGD (DSGD) for \texttt{squared-phillips} and \texttt{squared-shaw} is updated by
\begin{align*}
x_{k+1}^\delta = x_k^\delta - 2\eta_k \big(\langle a_{i_k}, x_k^\delta\rangle(\langle a_{i_k}, x_k^\delta\rangle^2-y_{i_k}^\delta)a_{i_k}+\lambda_k^\delta \langle \tilde{a}_{N,i_k}, x_k^\delta\rangle(\langle \tilde{a}_{N,i_k} , x_k^\delta\rangle^2-y_{i_k}^\delta)\tilde{a}_{N,i_k}\big)
\end{align*}
with $\tilde{a}_{N,i_k}^t$ being the $i_k$th row of $\tilde{A}_N$.
We set $\eta_k=c_0/\big(2\max_i(\|F'_{i}(x^\dag)\|^2)\big)$ and $\lambda_k^\delta=1$, and compare the convergence behavior of DSGD with that of SGD, LM, and the data-driven LM (DLM) using the same data-driven operators and regularization parameters as DSGD. 
For \texttt{squared-phillips}, the constant step size is chosen as $1/{\|F'(x^\dag)\|_F^2}$ for LM and $1/{2\|F'(x^\dag)\|_F^2}$ for DLM; while for \texttt{squared-shaw}, the constant step size is chosen as $2/{3\|F'(x^\dag)\|_F^2}$ for LM and $1/(3\|F'(x^\dag)\|_F^2)$ for DLM to ensure the convergence of the algorithms.
For \texttt{squared-phillips}, we set $c_0=2$ for SGD and $c_0=1$ for DSGD; while for \texttt{squared-shaw}, we set $c_0=4/3$ for SGD and $c_0=2/3$ for DSGD. 
We present the results for \texttt{squared-phillips} and \texttt{squared-shaw} in Figures \ref{fig:err_non_phil} and \ref{fig:err_non_shaw}, respectively.

\begin{figure}[hbt!]
\centering
\setlength{\tabcolsep}{4pt}
\begin{tabular}{cc}
\includegraphics[width=0.31\textwidth,trim={1.5cm 0 0.5cm 0.5cm}]{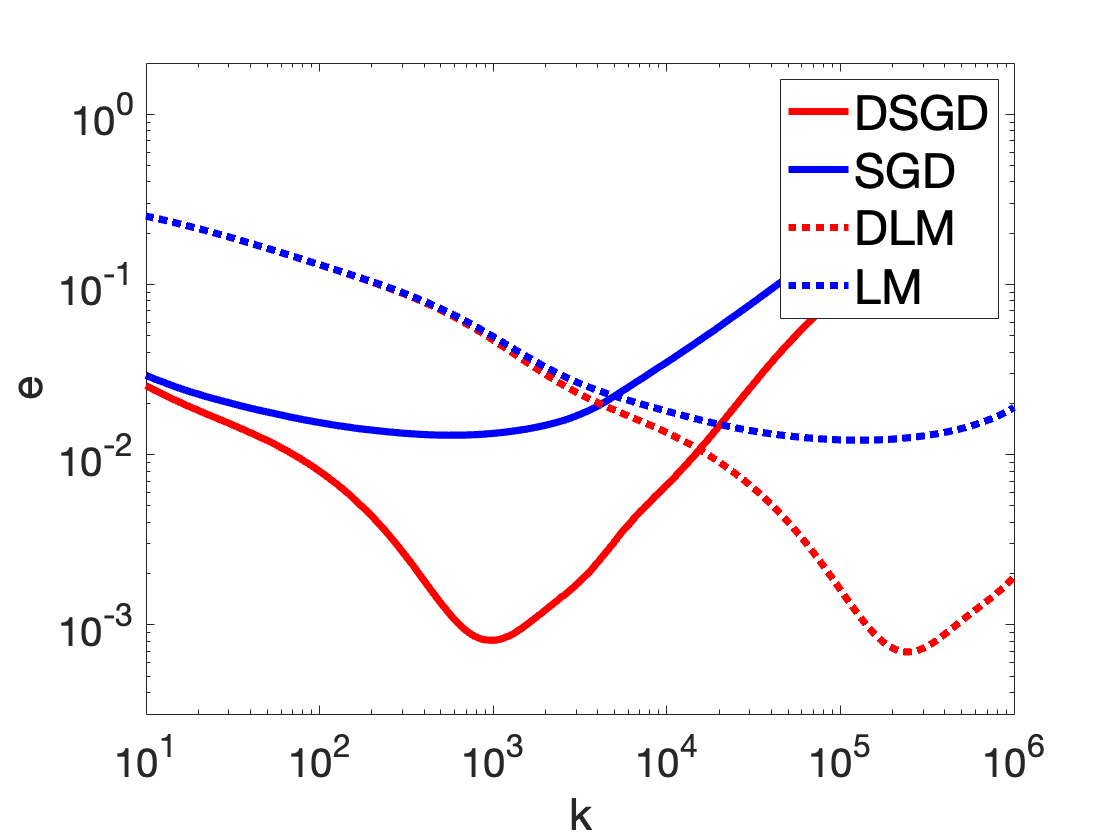}&
\includegraphics[width=0.31\textwidth,trim={1.5cm 0 0.5cm 0.5cm}]{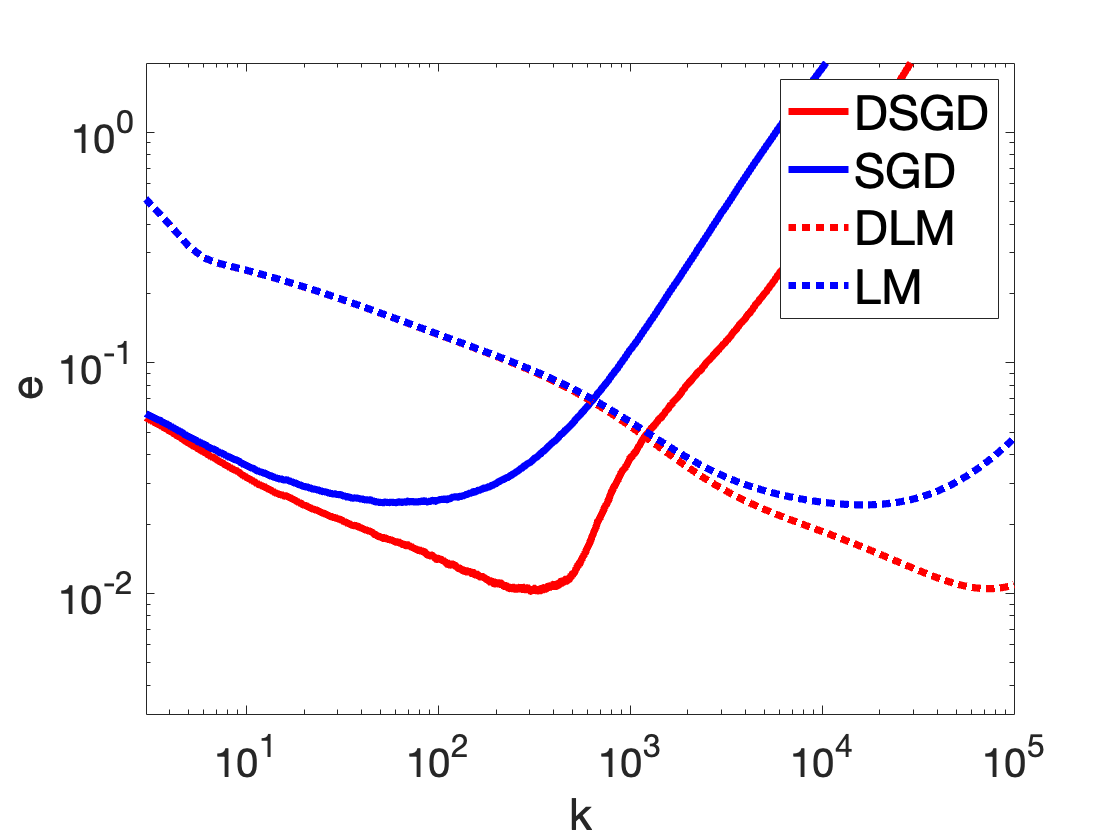}\\
$\delta_0=$ 1e-3 & $\delta_0=$1e-2
\end{tabular}
\caption{The convergence of relative mean squared errors $e=\frac{\E[\|x_k^\delta-x^\dag\|^2]}{\|x^\dag\|^2}$ of four methods for \texttt{squared-phillips}. \label{fig:err_non_phil}}
\end{figure}

\begin{figure}[hbt!]
\centering
\setlength{\tabcolsep}{4pt}
\begin{tabular}{cc}
\includegraphics[width=0.31\textwidth,trim={1.5cm 0 0.5cm 0.5cm}]{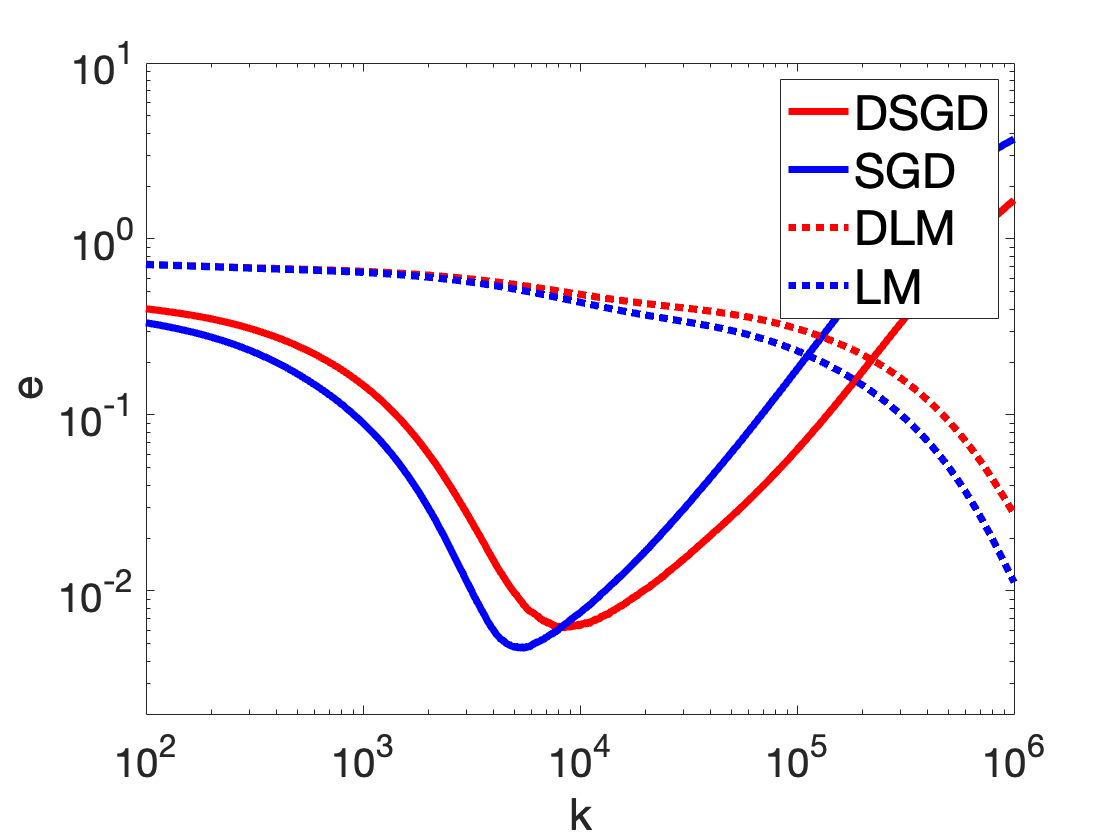}&
\includegraphics[width=0.31\textwidth,trim={1.5cm 0 0.5cm 0.5cm}]{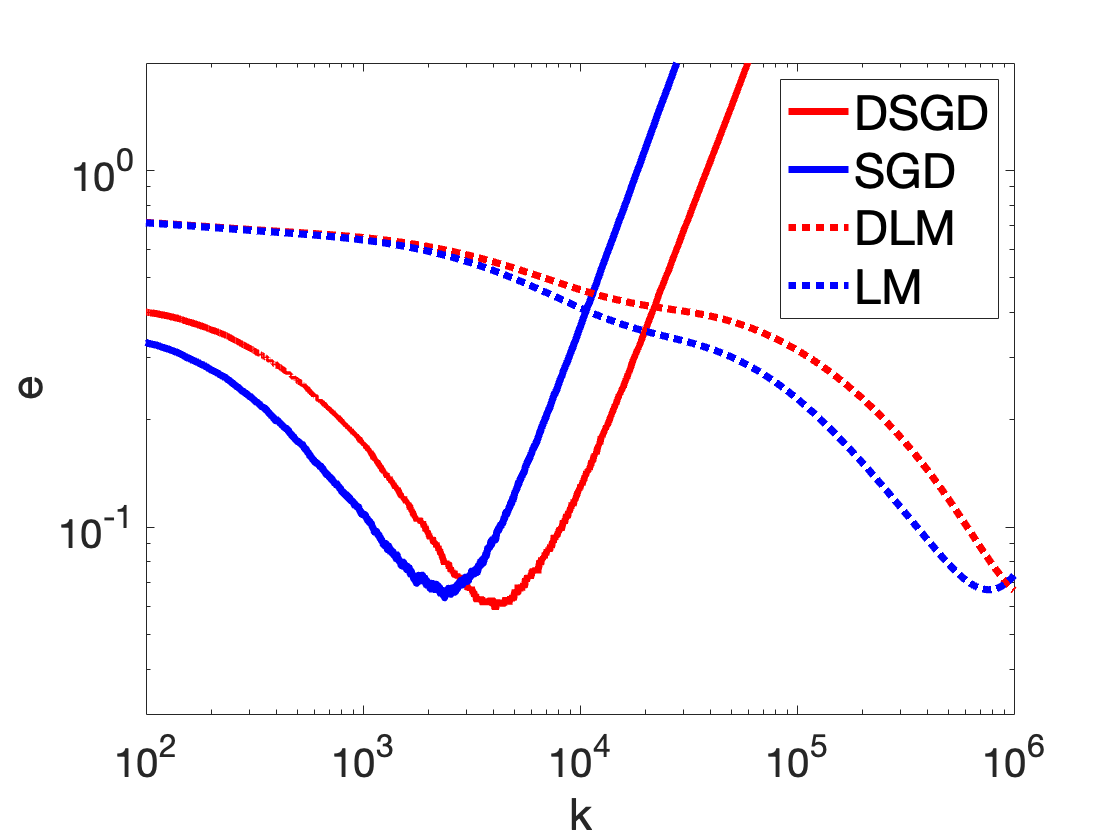}\\
$\delta_0=$ 5e-3 & $\delta_0=$5e-2
\end{tabular}
\caption{The convergence of relative mean squared errors $e=\frac{\E[\|x_k^\delta-x^\dag\|^2]}{\|x^\dag\|^2}$ of four methods for \texttt{squared-shaw}.\label{fig:err_non_shaw}}
\end{figure}

The results indicate that for both nonlinear problems, the stochastic methods, i.e., SGD and DSGD, are significantly more efficient than the corresponding deterministic methods, i.e., LM and DLM.
For \texttt{squared-phillips}, DSGD performs much better than SGD, while for \texttt{squared-shaw}, DSGD can achieve better accuracy than SGD with more iterations in the noisy case (when $\delta_0=$5e-2). 
However, in the less noisy case, DSGD performs slightly worse than SGD. 
These observations are mostly consistent with those for the linear problems in Section \ref{sec:num_linear}.
Thus, it is promising to improve the convergence behavior of DSGD by using decaying step size and regularization parameter schedules, as well as more suitable data-driven operators. 
We shall address this interesting topic in future work.
}

\section{Concluding remarks} \label{sec:conc}

In this work, we first established the regularizing property of a new data-driven regularized stochastic gradient descent (with a data-driven operator that can only partially explain the model for the true data) for a class of nonlinear inverse problems, under the tangential cone condition and \textit{a priori} rules on the parameter (step size, regularization parameter, and stopping index) choice. 
Then, we derived the convergence rates of this algorithm with polynomially decaying step size and regularization parameter schedules under the additional source condition, range invariance condition, and its stochastic variant. 
The analysis is motivated by both data-driven iteratively regularized Landweber iteration and the standard stochastic gradient descent for solving nonlinear inverse problems, and the results extend the existing works in \cite{AspriBanert:2020} and \cite{JinZhouZou:2020}. 
Finally, we present several numerical experiments on {both linear and nonlinear} inverse problems, demonstrating the advantages of the data-driven SGD over the standard SGD and Landweber method.

The algorithm proposed in this work combines the standard stochastic gradient descent method with a data-driven model introduced in the regularization term. 
It is known that training data can be used to increase the possibility of selecting better initial guesses which provide greater regularity indexes in the source condition and thus allow the algorithm to achieve higher convergence rates. 
Choosing appropriate initial guesses based on data-driven models to improve the convergence rates and providing theoretical support for it is an important topic that desires to be investigated.
We leave this interesting question to future works.

\section{Acknowledgments}
The author thanks Professor Fioralba Cakoni for the discussion on this work {and is also grateful to the three anonymous referees whose constructive comments have led to an improved presentation of the paper.}
Part of the work was completed during the visit to Hong Kong and Beijing supported by AMS-Simons Travel Grant.

\appendix
\section{Auxiliary estimates}\label{app:estimate}
In this appendix, we collect a set of supplementary estimates and lengthy technical proofs of several results. We begin with the proofs for analyzing the regularizing property of the data-driven SGD in Section \ref{sec:conv}.

\subsection{Proof of Proposition \ref{prop:mono-exact}}\label{app:prop:mono-exact}
{ Using a similar technique to that in \cite[Lemma 2.2]{AspriBanert:2020}, we first bound the mean squared residual of the data-driven model $G$, i.e., $\E[\|G(x_k^\delta)-y^\delta\|^2]$ in the following lemma which is used in Propositions \ref{prop:mono-exact} and \ref{prop:Cauchy}.}

\begin{lemma}\label{lem:resG}
Let Assumption \ref{ass:sol}(i) be fulfilled. Then for any data-driven SGD iterate $x_k^\delta\in\mathcal{B}_\rho(x^\dag)$ in \eqref{eqn:datasgd} and the error $e_k^\delta=x_k^\delta-x^\dag$,
there holds
\begin{equation*}
{\|G(x_k^\delta)-y^\delta\|\leq L_G\|e_k^\delta\|+\|G(x^\dag)-y^\delta\| \quad \mbox{and}\quad }\E[\|G(x_k^\delta)-y^\delta\|^2]^\frac12\leq L_G\E[\|e_k^\delta\|^2]^\frac12+\|G(x^\dag)-y^\delta\|.
\end{equation*}
Further, if Assumption \ref{ass:sol}(iii) is fulfilled, then there holds
\begin{equation*}
{\|G(x_k^\delta)-y^\delta\|\leq L_G\|e_k^\delta\|+C_{max}+\delta \quad \mbox{and}\quad }\E[\|G(x_k^\delta)-y^\delta\|^2]^\frac12\leq L_G\E[\|e_k^\delta\|^2]^\frac12+C_{max}+\delta.
\end{equation*}
\end{lemma}
\begin{proof}
By the triangle inequality and Assumption \ref{ass:sol}(i), there holds
\begin{align*}
    \|G(x_k^\delta)-y^\delta\|
    \leq& \|G(x_k^\delta)-G(x^\dag)\|+\|G(x^\dag)-y^\delta\|
    \leq \|\int_{0}^1 G'\big(x^\dag+t(x_k^\delta-x^\dag)\big)(x_k^\delta-x^\dag){\rm d} t\|+\|G(x^\dag)-y^\delta\| \\
    \leq& L_G\|x_k^\delta-x^\dag\|+\|G(x^\dag)-y^\delta\|.
\end{align*}
If Assumption \ref{ass:sol}(iii) holds, 
then we have
\begin{align*}
\|G(x^\dag)-y^\delta\|\leq\|G(x^\dag)-y^\dag\|+\|y^\dag-y^\delta\|\leq C_{max}+\delta.
\end{align*}
Finally, by taking full expectation, we obtain the desired assertion.
\end{proof}

Now, we give the proof of Proposition \ref{prop:est-ex}.
\begin{proof}
We define the inner product denoted by $\langle\cdot,\cdot\rangle$. With the definition of $x_{k}^\delta$ in \eqref{eqn:datasgd}, completing the square gives
\begin{align*}
    \|e_{k+1}^\delta\|^2 - \|e_{k}^\delta\|^2  
    \leq&2\eta_k\Big(\eta_k\|F'_{i_k}(x_k^\delta)^*(F_{i_k}(x_k^\delta)-y_{i_k}^\delta)\|^2-\langle e_k^\delta,  F_{i_k}'(x_k^\delta)^*(F_{i_k}(x_k^\delta)-y_{i_k}^\delta)\rangle\Big)\\   &+2\eta_k\lambda_k^\delta\Big(\eta_k\lambda_k^\delta\|G'_{i_k}(x_k^\delta)^*(G_{i_k}(x_k^\delta)-y_{i_k}^\delta)\|^2-\langle e_k^\delta,  G_{i_k}'(x_k^\delta)^*(G_{i_k}(x_k^\delta)-y_{i_k}^\delta)\rangle\Big)\\
    = & 2\eta_k\Big(\eta_k\|F'_{i_k}(x_k^\delta)^*(F_{i_k}(x_k^\delta)-y_{i_k}^\delta)\|^2-\langle F_{i_k}'(x_k^\delta)e_k^\delta,  F_{i_k}(x_k^\delta)-y_{i_k}^\delta\rangle\Big)\\   &+2\eta_k\lambda_k^\delta\Big(\eta_k\lambda_k^\delta\|G'_{i_k}(x_k^\delta)^*(G_{i_k}(x_k^\delta)-y_{i_k}^\delta)\|^2-\langle G_{i_k}'(x_k^\delta)e_k^\delta, G_{i_k}(x_k^\delta)-y_{i_k}^\delta\rangle\Big)\\
    :=& 2{\rm I}_1+2{\rm I}_2.
\end{align*}

Now, we bound ${\rm I}_1$ and ${\rm I}_2$ one by one. First, for ${\rm I}_1$, we split the factor $F'_{i_k}(x_k^\delta)e_k^\delta$ into three terms,
\begin{align*}
F_{i_k}'(x_k^\delta)e_k^\delta=&  (F_{i_k}(x_k^\delta)-y_{i_k}^\delta)+(y_{i_k}^\delta-F_{i_k}(x^\dag)) + (F_{i_k}(x^\dag)-F_{i_k}(x_k^\delta)+F_{i_k}'(x_k^\delta)e_k^\delta)\\
=& (F_{i_k}(x_k^\delta)-y_{i_k}^\delta)+\xi_{i_k} + (F_{i_k}(x^\dag)-F_{i_k}(x_k^\delta)-F_{i_k}'(x_k^\delta)(x^\dag - x_k^\delta)).
\end{align*}
Together with the inequality, derived directly from {Assumption \ref{ass:sol}(i), that $$\eta_k\|F'_{i_k}(x_k^\delta)^*(F_{i_k}(x_k^\delta)-y_{i_k}^\delta)\|^2\leq \eta_k L_F^2\|F_{i_k}(x_k^\delta)-y_{i_k}^\delta\|^2,$$}
we can bound ${\rm I}_1$ by
\begin{align}
&{\rm I}_1=\eta_k\Big(\eta_k\|F'_{i_k}(x_k^\delta)^*(F_{i_k}(x_k^\delta)-y_{i_k}^\delta)\|^2-\langle F_{i_k}'(x_k^\delta)e_k^\delta,  F_{i_k}(x_k^\delta)-y_{i_k}^\delta\rangle\Big)\nonumber\\
\leq&\eta_k\Big({(L_F^2\eta_k-1)}\|F_{i_k}(x_k^\delta)-y_{i_k}^\delta\|^2-\langle \xi_{i_k},  F_{i_k}(x_k^\delta)-y_{i_k}^\delta\rangle-\langle F_{i_k}(x^\dag)-F_{i_k}(x_k^\delta)-F_{i_k}'(x_k^\delta)(x^\dag - x_k^\delta),  F_{i_k}(x_k^\delta)-y_{i_k}^\delta\rangle\Big).\label{eqn:I1}
\end{align}
{
Then, under Assumption \ref{ass:sol}(ii), the Cauchy-Schwarz inequality and the triangle inequality $\| F_{i_k}(x_k^\delta)-y^\dag_{i_k}\|\leq\| F_{i_k}(x_k^\delta)-y^\delta_{i_k}\|+ \|\xi_{i_k}\|$ suggest that
\begin{align*}
{\rm I}_1
\leq & \eta_k\|F_{i_k}(x_k^\delta)-y_{i_k}^\delta\|\Big((L_F^2\eta_k-1)\|F_{i_k}(x_k^\delta)-y_{i_k}^\delta\|+\|\xi_{i_k}\| +\| F_{i_k}(x^\dag)-F_{i_k}(x_k^\delta)-F_{i_k}'(x_k^\delta)(x^\dag - x_k^\delta)\|\Big)\\
\leq & \eta_k\|F_{i_k}(x_k^\delta)-y_{i_k}^\delta\|\Big((L_F^2\eta_k-1)\|F_{i_k}(x_k^\delta)-y_{i_k}^\delta\|+\|\xi_{i_k}\| +\eta_F\| F_{i_k}(x_k^\delta)-y^\dag_{i_k}\|\Big)\\
\leq & \eta_k\|F_{i_k}(x_k^\delta)-y_{i_k}^\delta\|\Big((L_F^2\eta_k-1)\|F_{i_k}(x_k^\delta)-y_{i_k}^\delta\|+\|\xi_{i_k}\| +\eta_F(\| F_{i_k}(x_k^\delta)-y^\delta_{i_k}\|+ \|\xi_{i_k}\|)\Big)\\
\leq & \eta_k\|F_{i_k}(x_k^\delta)-y_{i_k}^\delta\|\Big((L_F^2\eta_k+\eta_F-1)\|F_{i_k}(x_k^\delta)-y_{i_k}^\delta\|+(1+\eta_F)\|\xi_{i_k}\|\Big).
\end{align*}
The identity $\delta^2\geq\|\xi\|^2=\frac1n \sum_{i=1}^n \|\xi_{i}\|^2$ implies that $\|\xi_{i}\|\leq \sqrt{n}\delta$ for any $i=1,\dots,n$, which yields that
\begin{align*}
{\rm I}_1
\leq & -(1-L_F^2\eta_k-\eta_F)\eta_k\|F_{i_k}(x_k^\delta)-y_{i_k}^\delta\|^2+\sqrt{n}(1+\eta_F)\eta_k \delta\|F_{i_k}(x_k^\delta)-y_{i_k}^\delta\|.
\end{align*}
Further, Young's inequality $2ab\leq c a^2 + c^{-1} b^2$,
with the choice $a=\|F_{i_k}(x_k^\delta)-y_{i_k}^\delta\|$, $b=\tfrac12\sqrt{n}(1+\eta_F)\delta$ and $c= (1-L_F^2\eta_k-\eta_F)>0$ gives that
\begin{align*}
{\rm I}_1
\leq & -c\eta_k a^2+2\eta_k ab\leq -c\eta_k a^2+\eta_k(c a^2 + c^{-1} b^2)=c^{-1} \eta_k b^2=\frac{n(1+\eta_F)^2}{4(1-L_F^2\eta_k-\eta_F)}\eta_k\delta^2.
\end{align*}

Similarly, for ${\rm I}_2$, we derive that following estimate with Assumption \ref{ass:sol}(i) and the Cauchy-Schwarz inequality: 
\begin{align*}
{\rm I}_2\leq&\eta_k\lambda_k^\delta\big(\eta_k\lambda_k^\delta L_G^2\|G_{i_k}(x_k^\delta)-y_{i_k}^\delta\|^2-\langle G_{i_k}'(x_k^\delta)e_k^\delta, G_{i_k}(x_k^\delta)-y_{i_k}^\delta\rangle\big)\\
\leq&n\eta_k\lambda_k^\delta\big(\eta_k\lambda_k^\delta L_G^2\frac1n\sum_{i=1}^n\|G_{i}(x_k^\delta)-y_{i}^\delta\|^2+L_G \|e_k^\delta\|\frac1n\sum_{i=1}^n\|G_{i}(x_k^\delta)-y_{i}^\delta\|\big)\\
\leq&n\eta_k\lambda_k^\delta\big(\eta_k\lambda_k^\delta L_G^2\|G(x_k^\delta)-y^\delta\|^2+L_G \|e_k^\delta\| \|G(x_k^\delta)-y^\delta\|\big):=n {\rm I}_3.
\end{align*}
Further, by Lemma \ref{lem:resG} (under Assumptions \ref{ass:sol}(i) and (iii)), we have
\begin{align*}
{\rm I}_3
\leq&\eta_k\lambda_k^\delta\big(\eta_k\lambda_k^\delta L_G^2(L_G\|e_k^\delta\|+C_{max}+\delta)^2+L_G \|e_k^\delta\| (L_G\|e_k^\delta\|+C_{max}+\delta)\big)\\
\leq&\eta_k\lambda_k^\delta\big((1+\eta_k\lambda_k^\delta L_G^2)L_G^2\|e_k^\delta\|^2+(1+2\eta_k\lambda_k^\delta L_G^2)L_G\|e_k^\delta\|(C_{max}+\delta)+\eta_k\lambda_k^\delta L_G^2(C_{max}+\delta)^2\big).
\end{align*}
The inequality $L_G\|e_k^\delta\|(C_{max}+\delta)\leq \tfrac12 (L_G^2\|e_k^\delta\|^2+(C_{max}+\delta)^2)$ implies that
\begin{align*}
{\rm I}_3
\leq&\eta_k\lambda_k^\delta L_G^2(\tfrac32+2\eta_k\lambda_k^\delta L_G^2)\|e_k^\delta\|^2+\eta_k\lambda_k^\delta(\tfrac12+2\eta_k\lambda_k^\delta L_G^2)(C_{max}+\delta)^2.
\end{align*}
Combining the above two estimates of ${\rm I}_1$ and ${\rm I}_2$ gives that
\begin{align*}
&\|e_{k+1}^\delta\|^2 \leq 2{\rm I}_1+2{\rm I}_2+\|e_{k}^\delta\|^2\leq2{\rm I}_1+2n{\rm I}_3+\|e_{k}^\delta\|^2\\
\leq &\big(1+2n\eta_k\lambda_k^\delta L_G^2(\tfrac32+2\eta_k\lambda_k^\delta L_G^2)\big)\|e_k^\delta\|^2+\frac{n(1+\eta_F)^2}{2(1-L_F^2\eta_k-\eta_F)}\eta_k\delta^2+2n\eta_k\lambda_k^\delta(\tfrac12+2\eta_k\lambda_k^\delta L_G^2)(C_{max}+\delta)^2.
\end{align*}}

{Next, we bound $\E[{\rm I}_1]$ and $\E[{\rm I}_2]$ using a similar strategy to that used for estimating ${\rm I}_1$ and ${\rm I}_2$.
Under Assumption \ref{ass:sol}(ii), by the measurability of the iterate $x_k^\delta$ with respect to the filtration
$\mathcal{F}_{k}$, we derive from \eqref{eqn:I1} that
\begin{align*}
\E[{\rm I}_1|\mathcal{F}_{k}]
\leq &{(L_F^2\eta_k-1)\frac{\eta_k}{n}}\sum_{i=1}^n\|F_{i}(x_k^\delta)-y_{i}^\delta\|^2-\frac{\eta_k}{n}\sum_{i=1}^n\langle \xi_{i},  F_{i}(x_k^\delta)-y_{i}^\delta\rangle\\
&-\frac{\eta_k}{n}\sum_{i=1}^n\langle F_{i}(x^\dag)-F_{i}(x_k^\delta)-F_{i}'(x_k^\delta)(x^\dag - x_k^\delta),  F_{i}(x_k^\delta)-y_{i}^\delta\rangle\\
= &{(L_F^2\eta_k-1)}\eta_k\|F(x_k^\delta)-y^\delta\|^2-\eta_k\langle \xi,  F(x_k^\delta)-y^\delta\rangle-\eta_k\langle F(x^\dag)-F(x_k^\delta)-F'(x_k^\delta)(x^\dag - x_k^\delta),  F(x_k^\delta)-y^\delta\rangle\\
\leq & {(L_F^2\eta_k-1)}\eta_k\|F(x_k^\delta)-y^\delta\|^2 + \eta_k\delta\|F(x_k^\delta)-y^\delta\| + \eta_k\eta_F\|F(x_k^\delta)-y^\dag\|\|F(x_k^\delta)-y^\delta\|\\
\leq & \eta_k\|F(x_k^\delta)-y^\delta\|\Big({(L_F^2\eta_k-1)}\|F(x_k^\delta)-y^\delta\| + \delta + \eta_F(\|F(x_k^\delta)-y^\delta\|+\delta)\Big)\\
\leq & \eta_k\|F(x_k^\delta)-y^\delta\|\Big({(L_F^2\eta_k+\eta_F-1)}\|F(x_k^\delta)-y^\delta\| +(1+\eta_F)\delta\Big).
\end{align*}}
{Then, under Assumption \ref{ass:sol}(i), we derive from the definition of ${\rm I}_2$ that 
\begin{align*}
\E[{\rm I}_2|\mathcal{F}_{k}]\leq&\E\big[\eta_k\lambda_k^\delta\big(\eta_k\lambda_k^\delta L_G^2\|G_{i_k}(x_k^\delta)-y_{i_k}^\delta\|^2-\langle G_{i_k}'(x_k^\delta)e_k^\delta, G_{i_k}(x_k^\delta)-y_{i_k}^\delta\rangle\big)\big|\mathcal{F}_{k}\big]\\
=&\eta_k\lambda_k^\delta\big(\eta_k\lambda_k^\delta L_G^2\|G(x_k^\delta)-y^\delta\|^2-\langle G'(x_k^\delta)e_k^\delta, G(x_k^\delta)-y^\delta\rangle\big)\\
=&\eta_k\lambda_k^\delta\big(\eta_k\lambda_k^\delta L_G^2\|G(x_k^\delta)-y^\delta\|^2+L_G \|e_k^\delta\| \|G(x_k^\delta)-y^\delta\|\big)={\rm I}_3.
\end{align*}
By taking full conditional of the inequality and using the triangle inequality,  we obtain that
\begin{align*}
\E[{\rm I}_1]
\leq & -{(1-L_F^2\eta_k-\eta_F)}\eta_k\E[\|F(x_k^\delta)-y^\delta\|^2] +(1+\eta_F)\eta_k\delta\E[\|F(x_k^\delta)-y^\delta\|^2]^\frac12\\
\mbox{and}\quad \E[{\rm I}_2]=\E[{\rm I}_3]
\leq&\eta_k\lambda_k^\delta L_G^2(\tfrac32+2\eta_k\lambda_k^\delta L_G^2)\E[\|e_k^\delta\|^2]+\eta_k\lambda_k^\delta(\tfrac12+2\eta_k\lambda_k^\delta L_G^2)(C_{max}+\delta)^2,
\end{align*}
which implies that 
\begin{align*}
&{\E[\|e_{k+1}^\delta\|^2]  \leq 2\E[{\rm I}_1]+2\E[{\rm I}_2]+\E[\|e_{k}^\delta\|^2]}\\
\leq & \big(1+2\eta_k\lambda_k^\delta L_G^2(\tfrac32+2\eta_k\lambda_k^\delta L_G^2)\big)\E[\|e_k^\delta\|^2]+2\eta_k\lambda_k^\delta(\tfrac12+2\eta_k\lambda_k^\delta L_G^2)(C_{max}+\delta)^2\\
&-2(1-L_F^2\eta_k-\eta_F)\eta_k\E[\|F(x_k^\delta)-y^\delta\|^2] +2(1+\eta_F)\eta_k\delta\E[\|F(x_k^\delta)-y^\delta\|^2]^\frac12.
\end{align*}
Finally, by Young's inequality $2ab\leq c a^2 + c^{-1} b^2$,
with the choice $a=\E[\|F(x_k^\delta)-y^\delta\|^2]^\frac12$, $b=(1+\eta_F)\delta$ and $c= 2(1-L_F^2\eta_k-\eta_F)>0$, we estimate the last two terms of the above upper bound of $\E[\|e_{k+1}^\delta\|^2]$ by
\begin{align*}
-2(1-L_F^2\eta_k-\eta_F)\eta_k\E[\|F(x_k^\delta)-y^\delta\|^2] +2(1+\eta_F)\eta_k\delta\E[\|F(x_k^\delta)-y^\delta\|^2]^\frac12\leq \frac{(1+\eta_F)^2}{2(1-L_F^2\eta_k-\eta_F)}\eta_k\delta^2.
\end{align*}
}
This completes the proof of the proposition.
\end{proof}

\subsection{Proof of Proposition \ref{prop:Cauchy}}\label{app:prop:Cauchy}
To prove Proposition \ref{prop:Cauchy}, we first collect a preliminary result from \cite{HankeNeubauerScherzer:1995} which is used in Proposition \ref{prop:Cauchy}. This
result is a useful characterization of all possible solutions $x^*$ of problem \eqref{eqn:nonlin} \cite[Proposition 2.1]{HankeNeubauerScherzer:1995}.
\begin{lemma}\label{lem:linear}
Let Assumptions \ref{ass:sol}(i) and (ii) be fulfilled. 
\begin{itemize}
\item[$\rm(i)$] {The following inequalities hold for any $x,\tilde x\in  \mathcal{B}_\rho(x^\dag)$:
\begin{equation*}
(1+\eta_F)^{-1}\|F'(x)(x-\tilde x)\|\leq \|F(x)-F(\tilde x)\|\leq (1-\eta_F)^{-1}\|F'(x)(x-\tilde x)\|.
\end{equation*}}
\item[$\rm(ii)$]If $x^*\in  \mathcal{B}_\rho(x^\dag)$ is a solution of \eqref{eqn:nonlin}, then any other solution $\tilde x^*\in  \mathcal{B}_\rho(x^\dag)$ satisfies
$x^*-\tilde x^* \in \mathcal{N}(F'(x^*))$, and vice versa.
\end{itemize} 
\end{lemma}

Now, we give the proof of Proposition \ref{prop:Cauchy}.
\begin{proof}
The argument below follows closely 
\cite[Theorem 2.5]{AspriBanert:2020} and \cite[Lemma 3.3]{JinZhouZou:2020}, which can be
traced back to \cite{McCormickRodrigue:1975}. For the convenience of readers, we state similar results to those in \cite[Lemma 3.3]{JinZhouZou:2020} first. For any $j\geq k$, choose an index $\ell$
with $j\geq \ell\geq k$ such that
\begin{equation}\label{eqn:minimal-res}
  \E[\|F(x_\ell)-y^\dag\|^2]\leq \E[\| F(x_i)-y^\dag\|^2], \quad \forall k\leq i \leq j.
\end{equation}
We claim that $\lim_{j\geq k, k\to\infty}\E[\|e_j-e_k\|^2]=0$ which implies that the sequence $\{x_k\}_{k\ge 1}$ is actually a Cauchy sequence. In fact, we can bound $\E[\|e_j-e_k\|^2]^\frac12$ with the triangle inequality 
\begin{equation*}
  \E[\|e_j-e_k\|^2]^\frac{1}{2} \leq \E[\|e_j-e_\ell\|^2]^\frac{1}{2} + \E[\|e_\ell-e_k\|^2]^\frac{1}{2},
\end{equation*}
where
\begin{equation}\label{eqn:err-exact}
\begin{aligned}
  \E[\|e_j-e_\ell\|^2] & = 2\E[\langle e_\ell-e_j,e_\ell\rangle] + \E[\|e_j\|^2] - \E[\|e_\ell\|^2],\\
  \E[\|e_\ell-e_k\|^2] & = 2\E[\langle e_\ell-e_k,e_\ell\rangle] + \E[\|e_k\|^2] - \E[\|e_\ell\|^2].
\end{aligned}
\end{equation}
By Corollary \ref{cor:mono}, $\{x_k\}_{k\geq 1}\subset\mathcal{B}_\rho(x^\dag)$ and $\{\E[\|e_{k}\|^2]\}_{k\geq 1}$ is a Cauchy sequence which implies that $$\lim_{j\geq \ell, \ell\to\infty}(\E[\|e_j\|^2] - \E[\|e_\ell\|^2])=0\; \mbox{  and } \lim_{\ell\geq k, k\to\infty}(\E[\|e_k\|^2] - \E[\|e_\ell\|^2])=0.$$
Now, we show that $\lim_{k\to \infty}\E[\langle e_\ell-e_k,e_\ell\rangle]=0$ and $\lim_{\ell\to \infty}\E[\langle e_\ell-e_j,e_\ell\rangle]=0$.
By the definition of the data-driven SGD iterate $x_k$ in \eqref{eqn:datasgd}, we have
\begin{equation*}
  e_\ell - e_k = \sum_{i=k}^{\ell-1}(e_{i+1}-e_i) = \sum_{i=k}^{\ell-1}\eta_i\Big(F_{i_i}'(x_i)^*(y_{i_i}^\dag-F_{i_i}(x_i))+\lambda_i^0 G_{i_i}'(x_i)^*(y_{i_i}^\dag-G_{i_i}(x_i))\Big).
\end{equation*}
Then we can bound $\E[\langle e_\ell-e_k,e_\ell\rangle]$, using the triangle inequality, by
\begin{align*}
  |\E[\langle e_\ell-e_k,e_\ell\rangle]| & = |\E[\sum_{i=k}^{\ell-1}\langle \eta_i \Big(F_{i_i}'(x_i)^*(y_{i_i}^\dag-F_{i_i}(x_i))+\lambda_i^0 G_{i_i}'(x_i)^*(y_{i_i}^\dag-G_{i_i}(x_i))\Big),e_\ell\rangle]|\\
  & \leq \sum_{i=k}^{\ell-1}\eta_i|\E[\langle F_{i_i}'(x_i)^*(y_{i_i}^\dag-F_{i_i}(x_i)),e_\ell\rangle]|+\sum_{i=k}^{\ell-1}\eta_i \lambda_i^0|\E [\langle G_{i_i}'(x_i)^*(y_{i_i}^\dag-G_{i_i}(x_i)),e_\ell\rangle]|\\
  & =\sum_{i=k}^{\ell-1}\eta_i|\E[\langle y_{i_i}^\dag-F_{i_i}(x_i),F_{i_i}'(x_i)(x_\ell-x^\dag)\rangle]|+\sum_{i=k}^{\ell-1}\eta_i \lambda_i^0|\E [\langle y_{i_i}^\dag-G_{i_i}(x_i),G_{i_i}'(x_i)e_\ell\rangle]|\\
  & := {\rm I}_1+{\rm I}_2.
\end{align*}
Next, we estimate ${\rm I}_1$ and ${\rm I}_2$ one by one. 
By taking the conditional expectation, together with the Cauchy-Schwarz inequality, we have
\begin{align*}
{\rm I}_1
  & = \sum_{i=k}^{\ell-1}\eta_i|\E[\langle y_{i_i}^\dag-F_{i_i}(x_i),F_{i_i}'(x_i)(x_\ell-x^\dag)\rangle]| =\sum_{i=k}^{\ell-1}\eta_i|\E\big[\E[\langle y_{i_i}^\dag-F_{i_i}(x_i),F_{i_i}'(x_i)(x_\ell-x^\dag)\rangle|\mathcal{F}_i]\big]|\\
  & =\sum_{i=k}^{\ell-1}\eta_i|\E[\langle y^\dag-F(x_i),F'(x_i)(x_\ell-x^\dag)\rangle]| \leq \sum_{i=k}^{\ell-1}\eta_i\E[\|y^\dag-F(x_i)\|^2]^\frac{1}{2}\E[\|F'(x_i)(x_\ell-x^\dag)\|^2]^\frac{1}{2}.
\end{align*}
By the decomposition $x_\ell-x^\dag=(x_\ell-x_i)+(x_i-x^\dag)$ and the triangle inequality, there holds
\begin{align*}
{\rm I}_1
  & \leq \sum_{i=k}^{\ell-1}\eta_i\E[\|y^\dag-F(x_i)\|^2]^\frac{1}{2}\E[\|F'(x_i)\big((x_\ell-x_i)+(x_i-x^\dag)\big)\|^2]^\frac{1}{2}\\
  & \leq \sum_{i=k}^{\ell-1}\eta_i\E[\|y^\dag-F(x_i)\|^2]^\frac{1}{2}\Big(\E[\|F'(x_i)(x_\ell-x_i)\|^2]^\frac{1}{2}+\E[\|F'(x_i)(x_i-x^\dag)\|^2]^\frac{1}{2}\Big).
\end{align*}
By Assumption \ref{ass:sol}(ii) and Lemma \ref{lem:linear}(i), we have
\begin{equation*}
\|F'(x_i)(x_i-x)\|\leq (1+\eta_F)\|F(x_i)-F(x)\|,
\end{equation*}
where $x=x^\dag$ or $x_\ell$ with the index $\ell$ satisfying the inequality \eqref{eqn:minimal-res}, which implies
\begin{align*}
  {\rm I}_1 
   & \leq (1+\eta_F)\sum_{i=k}^{\ell-1}\eta_i\E[\|y^\dag-F(x_i)\|^2]^\frac{1}{2}\Big(\E[\|F(x_i)-F(x_\ell)\|^2]^\frac{1}{2}+\E[\|F(x_i)-F(x^\dag)\|^2]^\frac{1}{2}\Big)\\
   & \leq (1+\eta_F)\sum_{i=k}^{\ell-1}\eta_i\E[\|y^\dag-F(x_i)\|^2]^\frac{1}{2}\Big(\E[\|F(x_i)-y^\dag+y^\dag-F(x_\ell)\|^2]^\frac{1}{2}+\E[\|F(x_i)-y^\dag\|^2]^\frac{1}{2}\Big)\\
   & \leq (1+\eta_F)\sum_{i=k}^{\ell-1}\eta_i\E[\|y^\dag-F(x_i)\|^2]^\frac{1}{2}\Big(2\E[\|F(x_i)-y^\dag\|^2]^\frac{1}{2}+\E[\|F(x_\ell)-y^\dag\|^2]^\frac{1}{2}\Big)\\
   & \leq 3(1+\eta_F)\sum_{i=k}^{\ell-1}\eta_i\E[\|F(x_i)-y^\dag\|^2].
\end{align*}

For ${\rm I}_2$,  the Cauchy-Schwarz inequality gives that
\begin{align*}
 {\rm I}_2 
 & =\sum_{i=k}^{\ell-1}\eta_i\lambda_i^0|\E[\langle y^\dag-G(x_i),G'(x_i)e_\ell\rangle]| \leq \sum_{i=k}^{\ell-1}\eta_i\lambda_i^0\E[\|y^\dag-G(x_i)\|^2]^\frac{1}{2}\E[\|G'(x_i)e_\ell\|^2]^\frac{1}{2}.
\end{align*}
By Assumption \ref{ass:sol}(i) and Lemma \ref{lem:resG}, there holds
\begin{align*}
 {\rm I}_2 
  & \leq \sum_{i=k}^{\ell-1}\eta_i\lambda_i^0(L_G\E[\|e_i\|^2]^\frac12+C_{max})L_G\E[\|e_\ell\|^2]^\frac{1}{2}.
\end{align*}
Then, with the fact that $\lim_{k\to \infty}\E[\|e_{k}\|^2]=C_e$ obtained from Corollary \ref{cor:mono}, there exists some $k_0\in\mathbb{N}$ such that for any $k\geq k_0$, $\E[\|e_k\|^2]\leq 2C_e$. Thus, for any $k\geq k_0$, we have
\begin{align*}
{\rm I}_2 
  & \leq (L_G(2C_e)^\frac12+C_{max})L_G(2C_e)^\frac{1}{2}\sum_{i=k}^{\ell-1}\eta_i\lambda_i^0.
\end{align*}

Combining the above two estimates of ${\rm I}_1$ and ${\rm I}_2$ gives that, for any $k>k_0$,
\begin{align*}
  |\E[\langle e_\ell-e_k,e_\ell\rangle]| 
  & \leq {\rm I}_1+{\rm I}_2\leq 3(1+\eta_F)\sum_{i=k}^{\ell-1}\eta_i\E[\|F(x_i)-y^\dag\|^2]+(L_G(2C_e)^\frac12+C_{max})L_G(2C_e)^\frac{1}{2}\sum_{i=k}^{\ell-1}\eta_i\lambda_i^0.
\end{align*}
Similarly, we can deduce for any $\ell\geq k_0$
\begin{equation*}
  |\E[\langle e_j-e_\ell,e_\ell\rangle]| \leq 3(1+\eta_F)\sum_{i=\ell}^{j-1}\eta_i\E[\|F(x_i)-y^\dag\|^2]+(L_G(2C_e)^\frac12+C_{max})L_G(2C_e)^\frac{1}{2}\sum_{i=\ell}^{j-1}\eta_i\lambda_i^0.
\end{equation*}
Under 
Assumption \ref{ass:stepsize}(i), these two estimates and Corollary \ref{cor:mono} imply $\lim_{k\to \infty}\E[\langle e_\ell-e_k,e_\ell\rangle]=0$ and $\lim_{\ell\to \infty}\E[\langle e_\ell-e_j,e_\ell\rangle]=0$. Thus, the sequence $\{e_k\}_{k\geq 1}$ and $\{x_k\}_{k\geq1}$ are Cauchy sequences.
\end{proof}

\subsection{Proof of Lemma \ref{lem:conti-noise}}\label{app:lem:conti-noise}
By Corollary \ref{cor:mono_delta}, for any $k\leq k(\delta)$, we have $x_k,x_k^\delta\in\mathcal{B}_\rho(x^\dag)$. Now, we prove the assertion by mathematical induction. The assertion holds
trivially for $k=1$, since $x_1^\delta - x_1=0$. Now suppose that it
holds for all indices up to $k$ and any path $(i_1,\ldots,i_{k-1})\in\mathcal{F}_k$.
Next, by the definitions of the data-driven SGD iterates $x_k$ and $x_k^\delta$ defined by \eqref{eqn:datasgd}:
\begin{align*}
   x_{k+1} = x_k - \eta_k \big(F_{i_k}'(x_k)^*(F_{i_k}(x_k)-y_{i_k}^\dag)+\lambda_k^0 G_{i_k}'(x_k)^*(G_{i_k}(x_k)-y_{i_k}^\dag)\big),\\
   x_{k+1}^\delta = x_k^\delta - \eta_k \big(F_{i_k}'(x_k^\delta)^*(F_{i_k}(x_k^\delta)-y_{i_k}^\delta)+\lambda_k^\delta G_{i_k}'(x_k^\delta)^*(G_{i_k}(x_k^\delta)-y_{i_k}^\delta)\big).
\end{align*}
Therefore, for any fixed path $(i_1,\ldots,i_k)$, there holds
\begin{align*}
   &x_{k+1}^\delta -x_{k+1} \\=& (x_k^\delta-x_k) - \eta_k \left(F_{i_k}'(x_k^\delta)^*(F_{i_k}(x_k^\delta)-y_{i_k}^\delta)-F_{i_k}'(x_k)^*(F_{i_k}(x_k)-y_{i_k}^\dag)\right)\\
   &- \eta_k \left(\lambda_k^\delta G_{i_k}'(x_k^\delta)^*(G_{i_k}(x_k^\delta)-y_{i_k}^\delta)-\lambda_k^0 G_{i_k}'(x_k)^*(G_{i_k}(x_k)-y_{i_k}^\dag)\right)\\
   =& (x_k^\delta-x_k) - \eta_k \Big(F_{i_k}'(x_k^\delta)^*\big((F_{i_k}(x_k^\delta)-y_{i_k}^\delta)-(F_{i_k}(x_k)-y_{i_k}^\dag)\big)+(F_{i_k}'(x_k^\delta)^*-F_{i_k}'(x_k)^*)(F_{i_k}(x_k)-y_{i_k}^\dag)\Big)\\
   &- \eta_k \Big(\lambda_k^\delta G_{i_k}'(x_k^\delta)^*\big((G_{i_k}(x_k^\delta)-y_{i_k}^\delta)-(G_{i_k}(x_k)-y_{i_k}^\dag)\big)+\lambda_k^\delta (G_{i_k}'(x_k^\delta)^*-G_{i_k}'(x_k)^*)(G_{i_k}(x_k)-y_{i_k}^\dag)\\
   &\qquad\quad+(\lambda_k^\delta-\lambda_k^0) G_{i_k}'(x_k)^*(G_{i_k}(x_k)-y_{i_k}^\dag)\Big)\\
   =& (x_k^\delta-x_k) - \eta_k \Big(F_{i_k}'(x_k^\delta)^*(F_{i_k}(x_k^\delta)-F_{i_k}(x_k)-\xi_{i_k})+(F_{i_k}'(x_k^\delta)^*-F_{i_k}'(x_k)^*)(F_{i_k}(x_k)-y_{i_k}^\dag)\Big)\\
   &- \eta_k \Big(\lambda_k^\delta G_{i_k}'(x_k^\delta)^*(G_{i_k}(x_k^\delta)-G_{i_k}(x_k)-\xi_{i_k})+\lambda_k^\delta (G_{i_k}'(x_k^\delta)^*-G_{i_k}'(x_k)^*)(G_{i_k}(x_k)-y_{i_k}^\dag)\\
   &\qquad\quad+(\lambda_k^\delta-\lambda_k^0) G_{i_k}'(x_k)^*(G_{i_k}(x_k)-y_{i_k}^\dag)\Big).
\end{align*}
Together with the triangle inequality, we have
\begin{align*}
 & \|x_{k+1}^\delta-x_{k+1}\|\\ \leq& \|x_k^\delta -x_k\|+ \eta_k \Big(\|F_{i_k}'(x_k^\delta)^*(F_{i_k}(x_k^\delta)-F_{i_k}(x_k)-\xi_{i_k})\|+\|(F_{i_k}'(x_k^\delta)^*-F_{i_k}'(x_k)^*)(F_{i_k}(x_k)-y_{i_k}^\dag)\|\Big)\\
   &+ \eta_k \Big(\lambda_k^\delta\| G_{i_k}'(x_k^\delta)^*(G_{i_k}(x_k^\delta)-G_{i_k}(x_k)-\xi_{i_k})\|+\lambda_k^\delta\| (G_{i_k}'(x_k^\delta)^*-G_{i_k}'(x_k)^*)(G_{i_k}(x_k)-y_{i_k}^\dag)\|\\
   &\qquad\quad+(\lambda_k^\delta-\lambda_k^0)\| G_{i_k}'(x_k)^*(G_{i_k}(x_k)-y_{i_k}^\dag)\|\Big)\\
\leq&\|x_k^\delta -x_k\|+\eta_k({\rm{I}_1}+{\rm{I}_2}),
\end{align*}
where
\begin{align*}
{\rm{I}_1} =&\|F_{i_k}'(x_k^\delta)^*\|\big(\|F_{i_k}(x_k^\delta)-F_{i_k}(x_k)\|+\|\xi_{i_k}\|\big)+\|F_{i_k}'(x_k^\delta)^*-F_{i_k}'(x_k)^*\|\|F_{i_k}(x_k)-y_{i_k}^\dag\|,\\
{\rm{I}_2} =&\lambda_k^\delta\| G_{i_k}'(x_k^\delta)^*\|\big(\|G_{i_k}(x_k^\delta)-G_{i_k}(x_k)\|+\|\xi_{i_k}\|\big)+\lambda_k^\delta\| G_{i_k}'(x_k^\delta)^*-G_{i_k}'(x_k)^*\|\|G_{i_k}(x_k)-y_{i_k}^\dag\|\\
&+(\lambda_k^\delta-\lambda_k^0)\| G_{i_k}'(x_k)^*\|\|G_{i_k}(x_k)-y_{i_k}^\dag\|.
\end{align*}
Then, by Assumption \ref{ass:sol}(i), we can bound ${\rm{I}_1}$ and ${\rm{I}_2}$ by
\begin{align*}
{\rm{I}_1}\leq& L_F\big(\|F_{i_k}(x_k^\delta)-F_{i_k}(x_k)\|+\delta\big)+\|F_{i_k}'(x_k^\delta)^*-F_{i_k}'(x_k)^*\|\|F_{i_k}(x_k)-y_{i_k}^\dag\|,\\{\rm{I}_2}\leq& \lambda_k^\delta L_G\big(\|G_{i_k}(x_k^\delta)-G_{i_k}(x_k)\|+\delta\big)+\big(\lambda_k^\delta\| G_{i_k}'(x_k^\delta)^*-G_{i_k}'(x_k)^*\|+(\lambda_k^\delta-\lambda_k^0)L_G\big)\|G_{i_k}(x_k)-y_{i_k}^\dag\|.
\end{align*}
Finally, by the induction hypothesis that $\lim_{\delta\to 0}\|x_k^\delta-x_k\|=0$, the continuity of $F_{i_k}$, $F'_{i_k}$, $G_{i_k}$ and $G'_{i_k}$, and the fact $\lim_{\delta \to 0}\lambda_k^\delta=\lambda_k^0$, we can derive that, for any path
$(i_1,\ldots,i_{k})\in\mathcal{F}_{k+1}$,
\begin{align*}
\lim_{\delta\to 0}  \|x_{k+1}^\delta-x_{k+1}\|=0,
\end{align*}
which implies $\lim_{\delta\to 0^+}\E[\|x_{k+1}^\delta-x_{k+1}\|^2]^\frac12=0.$
This completes the proof.

\subsection{Proof of Lemma \ref{lem:recursion-mean}}\label{app:lem:recursion-mean}
We first collect the following elementary bound on the linearization error $ \|H(x)-H(x^\dag)-K_H(x-x^\dag)\|$ for $H=F$ or $G$ from \cite{JinZhouZou:2020}.
\begin{lemma}\label{lem:R}
Under Assumption \ref{ass:sol}(iv), for $H=F$ or $G$ and any $x\in  \mathcal{B}_\rho(x^\dag)$, there holds
\begin{equation*}
  \|H(x)-H(x^\dag)-K_H(x-x^\dag)\| \leq \frac{c_H}{2}\|K_H(x-x^\dag)\|\|x-x^\dag\|.
\end{equation*}
Further, under Assumption \ref{ass:stoch}, there holds
\begin{equation*}
  \E[\|H(x_k^\delta)-H(x^\dag)-K_H(x_k^\delta-x^\dag)\|^2]^\frac12 \leq \frac{c_H}{1+\theta}\E[\|K_H(x_k^\delta-x^\dag)\|^2]^\frac12\E[\|x_k^\delta-x^\dag\|^2]^\frac{\theta}{2}.
\end{equation*}
\end{lemma}
Now, we give the proof of Lemma \ref{lem:recursion-mean}.
\begin{proof}
By the definition of the data-driven SGD iterate $x_k^\delta$ in \eqref{eqn:datasgd} and Assumption \ref{ass:sol}(iv), there holds
\begin{align*}
  e_{k+1}^\delta & = e_k^\delta - \eta_k \big(F_{i_k}'(x_k^\delta)^*(F_{i_k}(x_k^\delta)-y_{i_k}^\delta)+\lambda_k^\delta G_{i_k}'(x_k^\delta)^*(G_{i_k}(x_k^\delta)-y_{i_k}^\delta)\big)\\
  & = e_k^\delta - \eta_k \big(K_{F,i_k}^* R_{F,x_k^\delta}^{i_k *}(F_{i_k}(x_k^\delta)-y_{i_k}^\delta)+\lambda_k^\delta K_{G,i_k}^*R_{G,x_k^\delta}^{i_k *}(G_{i_k}(x_k^\delta)-y_{i_k}^\delta)\big): = e_k^\delta - \eta_k \big({\rm I}_{F,k,i_k}+\lambda_k^\delta {\rm I}_{G,k,i_k}\big).
\end{align*}
Then we decompose ${\rm I}_{H,k,i_k}$ for $H=F$ or $G$ into
\begin{align*}
{\rm I}_{H,k,i_k} =&K_{H,i_k}^* R_{H,x_k^\delta}^{i_k *}(H_{i_k}(x_k^\delta)-y_{i_k}^\delta)
=K_{H,i_k}^*( R_{H,x_k^\delta}^{i_k *}-I )(H_{i_k}(x_k^\delta)-y_{i_k}^\delta)+K_{H,i_k}^* (H_{i_k}(x_k^\delta)-y_{i_k}^\delta)\\
=&K_{H,i_k}^*( R_{H,x_k^\delta}^{i_k *}-I )(H_{i_k}(x_k^\delta)-y_{i_k}^\delta)+K_{H,i_k}^* K_{H,i_k}(x_k^\delta-x^\dag)\\
&+K_{H,i_k}^* (H_{i_k}(x_k^\delta)-H_{i_k}(x^\dag)-K_{H,i_k}(x_k^\delta-x^\dag)+H_{i_k}(x^\dag)-y_{i_k}^\dag-\xi_{i_k})\\
:=&K_{H,i_k}^* K_{H,i_k}e_k^\delta+K_{H,i_k}^*v_{H,k,i_k},
\end{align*}
{where the random variables $v_{H,k,i_k}$ and $v_{G,k,i_k}$ are defined in \eqref{eqn:w_Fik} and \eqref{eqn:w_Gik} respectively.}
Thus, by the measurability of the iterate $x_k^\delta$ (and thus $e_k^\delta$) with respect to the filtration
$\mathcal{F}_{k}$, the conditional expectation $\E[e_{k+1}^\delta|\mathcal{F}_{k}]$ is given by
\begin{align*}
 &\E[e_{k+1}^\delta|\mathcal{F}_{k}]
   = e_k^\delta - \frac{\eta_k}{n}\sum_{i=1}^n \big({\rm I}_{F,k,i}+\lambda_k^\delta {\rm I}_{G,k,i}\big)\\
   = \;&e_k^\delta - \frac{\eta_k}{n}\sum_{i=1}^n \big(K_{F,i}^* K_{F,i}e_k^\delta+K_{F,i}^*v_{F,k,i}\big)- \frac{\eta_k\lambda_k^\delta }{n}\sum_{i=1}^n \big(K_{G,i}^* K_{G,i}e_k^\delta+K_{G,i}^*v_{G,k,i}\big)\\
   = \;&e_k^\delta - \eta_k\big(K_{F}^* K_{F}e_k^\delta+K_{F}^*v_{F,k}\big)- \eta_k\lambda_k^\delta \big(K_{G}^* K_{G}e_k^\delta+K_{G}^*v_{G,k}\big)\\
   = \;&\big(I-\eta_k(K_{F}^* K_{F}+\lambda_k^\delta K_{G}^* K_{G})\big)e_k^\delta - \eta_k K_{F}^*v_{F,k}- \eta_k\lambda_k^\delta K_{G}^*v_{G,k},
\end{align*}
where the random variables $v_{F,k}$ and $v_{G,k}$ are defined in \eqref{eqn:w_F} and \eqref{eqn:w_G}.
Then taking full conditional, with $B_H=K_H^* K_H$ for $H=F$ or $G$, there holds
\begin{align*}
  \E[e_{k+1}^\delta]   = \big(I-\eta_k(B_{F}+\lambda_k^\delta B_{G})\big)\E[e_k^\delta] - \eta_k (K_{F}^*\E[v_{F,k}]+ \lambda_k^\delta K_{G}^*\E[v_{G,k}]).
\end{align*}
Thus, with the notation $\Pi_j^k(B)$ from \eqref{eqn:Pi-B}, applying the recursion repeatedly yields
\begin{equation*}
  \E[e_{k+1}^\delta] = \Pi_{1}^k(B)e_1^\delta -\sum_{j=1}^{k}\eta_j\Pi_{j+1}^{k}(B) (K_{F}^*\E[v_{F,j}]+ \lambda_j^\delta K_{G}^*\E[v_{G,j}]).
\end{equation*}
This completes the proof of the lemma.
\end{proof}

\subsection{Proof of Lemma \ref{lem:recursion-var}}\label{app:lem:recursion-var}
Collected from the proof of Lemma \ref{lem:recursion-mean}, we rewrite the error $e_{k+1}^\delta=x_{k+1}^\delta-x^\dag$ and the mean error $\E[e_{k+1}^\delta]$ as
\begin{align*}
  e_{k+1}^\delta  =& e_k^\delta - \eta_k \Big(K_{F,i_k}^* K_{F,i_k}e_k^\delta+K_{F,i_k}^*v_{F,k,i_k}+\lambda_k^\delta \big(K_{G,i_k}^* K_{G,i_k}e_k^\delta+K_{G,i_k}^*v_{G,k,i_k}\big)\Big)\\
  =& \big(I-\eta_k(K_{F,i_k}^* K_{F,i_k}+\lambda_k^\delta K_{G,i_k}^* K_{G,i_k})\big)e_k^\delta - \eta_k (K_{F,i_k}^*v_{F,k,i_k}+\lambda_k^\delta K_{G,i_k}^*v_{G,k,i_k}),\\
  \E[e_{k+1}^\delta]   =& \big(I-\eta_k(B_{F}+\lambda_k^\delta B_{G})\big)\E[e_k^\delta] - \eta_k (K_{F}^*\E[v_{F,k}]+ \lambda_k^\delta K_{G}^*\E[v_{G,k}]).
\end{align*}
where the random variables $v_{F,k,i_k}$, $v_{G,k,i_k}$, $v_{F,k}$ and $v_{G,k}$ are defined in \eqref{eqn:w_Fik}, \eqref{eqn:w_Gik}, \eqref{eqn:w_F} and \eqref{eqn:w_G} respectively.
Then, subtracting the recursion for $\E[e_{k+1}^\delta]$ from that for $e_{k+1}^\delta$ indicates that the random variable\\ $z_{k+1}:=e_{k+1}^\delta-\E[e_{k+1}^\delta]$ satisfies
\begin{align}
z_{k+1}  =& \big(I-\eta_k(B_F+\lambda_k^\delta B_G)\big)e_k^\delta +\eta_k\big(B_F-K_{F,i_k}^* K_{F,i_k}+\lambda_k^\delta (B_G- K_{G,i_k}^* K_{G,i_k})\big)e_k^\delta \nonumber\\
&- \eta_k (K_{F,i_k}^*v_{F,k,i_k}+\lambda_k^\delta K_{G,i_k}^*v_{G,k,i_k})-\big(I-\eta_k(B_{F}+\lambda_k^\delta B_{G})\big)\E[e_k^\delta] + \eta_k (K_{F}^*\E[v_{F,k}]+ \lambda_k^\delta K_{G}^*\E[v_{G,k}])\nonumber\\
= &\big(I-\eta_k(B_F+\lambda_k^\delta B_G)\big)z_k +\eta_k\big(B_F-K_{F,i_k}^* K_{F,i_k}+\lambda_k^\delta (B_G- K_{G,i_k}^* K_{G,i_k})\big)e_k^\delta \nonumber\\
&+ \eta_k \big(K_{F}^*\E[v_{F,k}]-K_{F,i_k}^*v_{F,k,i_k}+ \lambda_k^\delta (K_{G}^*\E[v_{G,k}]-K_{G,i_k}^*v_{G,k,i_k})\big)\nonumber\\
= &\big(I-\eta_k(B_F+\lambda_k^\delta B_G)\big)z_k +\eta_k M_{k,1}+\eta_k M_{k,2},\label{eqn:recurs-z}
\end{align}
with the random variables $M_{j,1}$ and $M_{j,2}$ given by
\begin{align*}
   M_{j,1} & =  \big(B_F-K_{F,i_j}^* K_{F,i_j}+\lambda_j^\delta (B_G- K_{G,i_j}^* K_{G,i_j})\big)e_j^\delta ,\\
   M_{j,2} & = K_{F}^*\E[v_{F,j}]-K_{F,i_j}^*v_{F,j,i_j}+ \lambda_j^\delta (K_{G}^*\E[v_{G,j}]-K_{G,i_j}^*v_{G,k,i_j}).
\end{align*}
With the initial condition $z_1=0$ (since $x_1^\delta$ is deterministic), we repeatedly apply the recursion \eqref{eqn:recurs-z} and obtain a formula for $z_{k+1}$ that
\begin{equation*}
  z_{k+1} = \sum_{j=1}^{k} \eta_j \Pi_{j+1}^k(B) (M_{j,1}+M_{j,2}).
\end{equation*}
The random variables $M_{j,1}$ (the conditionally independent factor) and $M_{j,2}$ (the conditionally dependent factor) represent the iteration noise, due to the random choice of the index $i_j$.
In fact, for any $i>j$, by the measurability of $x_i^\delta$ and $x_j^\delta$ with respect to the filtration $\mathcal{F}_{i}$, we derive that
\begin{equation*}
  \E[\langle M_{i,1},M_{j,1}\rangle] =\E[\E[\langle M_{i,1},M_{j,1}\rangle|\mathcal{F}_{j}]]=\E[\langle \E[M_{i,1}|\mathcal{F}_{j}],M_{j,1}\rangle]=\E[\langle 0,M_{j,1}\rangle]=0,
\end{equation*}
which directly implies the conditional independence.
Further, a similar argument yields $\E[\langle M_{i,1},M_{j,2}\rangle] = 0$, for any $i> j$. 
Then we can decompose the weighted computational variance $\E[\|B^sz_{k+1}\|^2]$ as
\begin{align*}
\E[\|B_F^sz_{k+1}\|^2] &= \sum_{i=1}^k\sum_{j=1}^k\eta_i\eta_j\E[\langle B_F^s\Pi_{i+1}^k(B)(M_{i,1}+M_{i,2}),B_F^s\Pi_{j+1}^k(B)(M_{j,2}+M_{j,2})\rangle]\\
&= \sum_{i=1}^k\sum_{j=1}^k\eta_i\eta_j\E[\langle B_F^s\Pi_{i+1}^k(B)M_{i,1},B_F^s\Pi_{j+1}^k(B)M_{j,1}\rangle]\\
& \quad + 2\sum_{i=1}^k\sum_{j=1}^k\eta_i\eta_j\E[\langle B_F^s\Pi_{i+1}^k(B)M_{i,1},B_F^s\Pi_{j+1}^k(B)M_{j,2}\rangle]\\
& \quad + \sum_{i=1}^k\sum_{j=1}^k\eta_i\eta_j\E[\langle B_F^s\Pi_{i+1}^k(B)M_{i,2},B_F^s\Pi_{j+1}^k(B)M_{j,2}\rangle] \\
&= \sum_{j=1}^k \eta_j^2\E[\|B_F^s\Pi_{j+1}^k(B)M_{j,1}\|^2]+ 2\sum_{j=1}^k\sum_{i=1}^j \eta_i\eta_j\E[\langle B_F^s\Pi_{i+1}^kM_{i,1}, B_F^s\Pi_{j+1}^kM_{j,2}\rangle]\\
& \quad + \sum_{j=1}^k\sum_{i=1}^k\eta_i\eta_j\E[\langle B_F^s\Pi_{i+1}^k(B)M_{i,2},B_F^s\Pi_{j+1}^k(B)M_{j,2}\rangle].
\end{align*}
With the notation $\varphi_{i}$ that denotes the $i$th Cartesian basis vector in $\mathbb{R}^n$ scaled by $n^\frac12$, we can rewrite the random variables $M_{j,1}$ and $M_{j,2}$ as
\begin{align*}
   M_{j,1} & =  (K_F^*K_F+\lambda_j^\delta K_G^*K_G)e_j^\delta-
   (K_{F}^* K_{F,i_j}+\lambda_j^\delta K_{G}^* K_{G,i_j})e_j^\delta\varphi_{i_j} ,\\
   M_{j,2} & = K_{F}^*\E[v_{F,j}]-K_{F}^*v_{F,j,i_j}\varphi_{i_j}+ \lambda_j^\delta (K_{G}^*\E[v_{G,j}]-K_{G}^*v_{G,k,i_j}\varphi_{i_j}).
\end{align*}
Further, under Assumption \ref{ass:sol}(v), there holds
\begin{align*}
   M_{j,1} & =  (K_F^*K_F+\lambda_j^\delta K_F^*R^*K_G)e_j^\delta-
   (K_{F}^* K_{F,i_j}+\lambda_j^\delta K_F^*R^* K_{G,i_j})e_j^\delta\varphi_{i_j}\\
   & =  K_F^*\big((K_F+\lambda_j^\delta R^*K_G)e_j^\delta-
   ( K_{F,i_j}+\lambda_j^\delta R^* K_{G,i_j})e_j^\delta\varphi_{i_j} \big):=K_F^*N_{j,1},\\
   M_{j,2}  & = K_{F}^*\E[v_{F,j}]-K_{F}^*v_{F,j,i_j}\varphi_{i_j}+ \lambda_j^\delta (K_{F}^*R^*\E[v_{G,j}]-K_{F}^*R^*v_{G,k,i_j}\varphi_{i_j})\\
   & = K_{F}^*\big(\E[v_{F,j}]-v_{F,j,i_j}\varphi_{i_j}+ \lambda_j^\delta R^*(\E[v_{G,j}]-v_{G,k,i_j}\varphi_{i_j})\big):=K_F^*N_{j,2}.
\end{align*}
Thus, using the identity $\|B_F^s\Pi_{j+1}^k(B)K_F^*\|^2=\|B_F^{s+\frac12}\Pi_{j+1}^k(B)\|^2=\|B_F^{\tilde{s}}\Pi_{j+1}^k(B)\|^2=(\phi_j^{\tilde{s}})^2$ and the Cauchy-Schwarz inequality, we can rewrite the decomposition of the weighted computational variance $\E[\|B^sz_{k+1}\|^2]$ as
\begin{align*}
\E[\|B_F^sz_{k+1}\|^2] 
=& \sum_{j=1}^k \eta_j^2\E[\|B_F^s\Pi_{j+1}^k(B)K_F^*N_{j,1}\|^2]+ 2\sum_{j=1}^k\sum_{i=1}^j\eta_i\eta_j\E[\langle B_F^s\Pi_{i+1}^k(B)K_F^*N_{i,1}, B_F^s\Pi_{j+1}^k(B)K_F^*N_{j,2}\rangle]\\
& + \sum_{j=1}^k\sum_{i=1}^k\eta_i\eta_j\E[\langle B_F^s\Pi_{i+1}^k(B)K_F^*N_{i,2},B_F^s\Pi_{j+1}^k(B)K_F^*N_{j,2}\rangle]\\
\leq& \sum_{j=1}^k \eta_j^2(\phi_j^{\tilde{s}})^2\E[\|N_{j,1}\|^2]+ 2\sum_{j=1}^k\sum_{i=1}^k\eta_i\eta_j\phi_i^{\tilde{s}}\phi_j^{\tilde{s}}\E[\|N_{i,1}\| \|N_{j,2}\|]+ \sum_{i=1}^k\sum_{j=1}^k\eta_i\eta_j\phi_i^{\tilde{s}}\phi_j^{\tilde{s}}\E[\|N_{i,2}\| \|N_{j,2}\|]\\
\leq&\sum_{j=1}^k\eta_j^2(\phi_j^{\tilde{s}})^2\E[\|N_{j,1}\|^2]+\sum_{j=1}^k\sum_{i=1}^k\eta_i\eta_j\phi_i^{\tilde{s}}\phi_j^{\tilde{s}}(2\E[\|N_{i,1}\|^2]^\frac12+\E[\|N_{i,2}\|^2]^\frac12)\E[\|N_{j,2}\|^2]^\frac12.
\end{align*}
Finally, the equation 
\begin{align*}
&\sum_{j=1}^k\sum_{i=1}^k\eta_i\eta_j\phi_i^{\tilde{s}}\phi_j^{\tilde{s}}(2\E[\|N_{i,1}\|^2]^\frac12+\E[\|N_{i,2}\|^2]^\frac12)\E[\|N_{j,2}\|^2]^\frac12\\
=&\big(\sum_{j=1}^k\eta_j\phi_j^{\tilde{s}}(2\E[\|N_{j,1}\|^2]^\frac12+\E[\|N_{j,2}\|^2]^\frac12)\big)\big(\sum_{j=1}^k\eta_j\phi_j^{\tilde{s}}\E[\|N_{j,2}\|^2]^\frac12\big)    
\end{align*} 
completes the proof of the lemma.

\subsection{Proof of Lemma \ref{lem:bound-N}}\label{app:lem:bound-N}

{To prove Lemma \ref{lem:bound-N}, we first derive a refined estimate for the residual $\E[\|G(x_k^\delta)-y^\delta\|^2]^\frac12$, under Assumptions \ref{ass:sol}(i)(iii)(iv), which is also used in the proof of Lemma \ref{lem:bound-w}. 
\begin{lemma}\label{lem:resG_refined}
Let Assumptions \ref{ass:sol}(i)(iv) be fulfilled. Then there holds
\begin{equation*}
\E[\|G(x_k^\delta)-y^\delta\|^2]^\frac12\leq (c_G\E[\|e_k^\delta\|^2]^\frac12+1)\E[ \|K_{G}e_k^\delta\|^2]^\frac12+\|G(x^\dag)-y^\delta\|.
\end{equation*}
Further, if Assumption \ref{ass:sol}(iii) is fulfilled, then there holds
\begin{equation*}
\E[\|G(x_k^\delta)-y^\delta\|^2]^\frac12\leq (c_G\E[\|e_k^\delta\|^2]^\frac12+1)\E[ \|K_{G}e_k^\delta\|^2]^\frac12+C_{max}+\delta.
\end{equation*}
\end{lemma}
\begin{proof}
Following the technique used in the proof of Lemma \ref{lem:resG} and the triangle inequality, there holds
\begin{align*}
    \|G(x_k^\delta)-y^\delta\|
    \leq \|G(x_k^\delta)-G(x^\dag)\|+\|G(x^\dag)-y^\dag\|+\|y^\dag-y^\delta\|,
\end{align*}
where
\begin{align*}
    \|G(x_k^\delta)-G(x^\dag)\|
    \leq& \|\int_{0}^1 G'\big(x^\dag+t(x_k^\delta-x^\dag)\big)(x_k^\delta-x^\dag){\rm d} t\|
    \leq \int_{0}^1 \|R_{G,x^\dag+t(x_k^\delta-x^\dag)} K_{G} e_k^\delta\|{\rm d} t \\
    \leq& \int_{0}^1 (\|R_{G,x^\dag+t(x_k^\delta-x^\dag)}-I\|+1) \|K_{G} e_k^\delta\|{\rm d} t 
    \leq (c_G\|e_k^\delta\|+1) \|K_{G} e_k^\delta\|.
\end{align*}
Further, with Assumption \ref{ass:sol}(iii), we have
\begin{align*}
\|G(x_k^\delta)-y^\delta\|\leq(c_G\|e_k^\delta\|+1) \|K_{G} e_k^\delta\|+C_{max}+\delta.
\end{align*}
Finally, by taking full expectation, we obtain the desired assertion.
\end{proof}
Now, we shall give the proof of Lemma \ref{lem:bound-N}.}
\begin{proof}
First, we derive an estimate for $\E[\|N_{j,1}\|^2]^\frac12$. Under Assumption \ref{ass:sol}(v), using the definition of $N_{j,1}$ in \eqref{eqn:N} and the triangle inequality, we may bound $\E[\|N_{j,1}\|^2]^\frac12$ by
\begin{align*}
\E[\|N_{j,1}\|^2]^\frac12 & \leq \E[\|K_F e_j^\delta -K_{F,i_j}e_j^\delta\varphi_{i_j}\|^2]^\frac12 + \lambda_j^\delta\E[\| R^*(K_Ge_j^\delta -K_{G,i_j}e_j^\delta\varphi_{i_j})\|^2]^\frac12\\
& \leq \E[\|K_F e_j^\delta -K_{F,i_j}e_j^\delta\varphi_{i_j}\|^2]^\frac12 + c_R\lambda_j^\delta\E[\| K_Ge_j^\delta -K_{G,i_j}e_j^\delta\varphi_{i_j}\|^2]^\frac12.
\end{align*}
With the measurability of the data-driven SGD iterate error $e_j^\delta=x_j^\delta-x^\dag$ with respect to the filtration $\mathcal{F}_j$, it directly implies that $\E[K_{H,i_j}e_j^\delta\varphi_{i_j}|
\mathcal{F}_{j}]= K_H e_j^\delta$ for $H=F$ or $G$. Thus, by the bias-variance decomposition and the definitions of $K_H$ and $K_{H,i}$ in Assumption \ref{ass:sol}(iv), the conditional expectation $\E[\| K_He_j^\delta -K_{H,i_j}e_j^\delta\varphi_{i_j}\|^2|\mathcal{F}_j]$ can be bounded by
\begin{align*}
&\E[\| K_He_j^\delta -K_{H,i_j}e_j^\delta\varphi_{i_j}\|^2|\mathcal{F}_j]=\E[\|K_{H,i_j}e_j^\delta\varphi_{i_j}\|^2|\mathcal{F}_j] -\E[\| K_He_j^\delta\|^2|\mathcal{F}_j]\leq\frac1n \sum_{i=1}^n \|K_{H,i}e_j^\delta\varphi_{i}\|^2\\
=&\frac1n \sum_{i=1}^n (n\|K_{H,i}e_j^\delta\|^2)=\frac1n n^2\|K_{H}e_j^\delta\|^2=n\|K_{H}e_j^\delta\|^2.
\end{align*}
Together with Assumption \ref{ass:sol}(v), we derive the following estimate by taking full expectation,
\begin{align*}
\E[\|N_{j,1}\|^2]^\frac12  \leq& n^\frac12\E[\|K_{F}e_j^\delta\|^2]^\frac12+n^\frac12c_R\lambda_j^\delta\E[\|K_{G}e_j^\delta\|^2]^\frac12
=n^\frac12\E[\|K_{F}e_j^\delta\|^2]^\frac12+n^\frac12c_R\lambda_j^\delta\E[\|RK_{F}e_j^\delta\|^2]^\frac12\\
\leq& n^\frac12(1+c_R^2\lambda_j^\delta)\E[\|K_{F}e_j^\delta\|^2]^\frac12=n^\frac12(1+c_R^2\lambda_j^\delta)\E[\|B_{F}^\frac12 e_j^\delta\|^2]^\frac12.
\end{align*}
Similarly, using the telescopic expectation identity $\E_{\mathcal{F}_j}
[\E[v_{H,j,i_j}\varphi_{i_j}|\mathcal{F}_j]]=\E_{\mathcal{F}_{j}}[v_{H,j}]$ for $H=F$ or $G$, where $\E_{\mathcal{F}_j}$ denotes
taking expectation in $\mathcal{F}_j$, we obtain that
\begin{align*}
\E[\|\E[v_{H,j}]-v_{H,j,i_j}\varphi_{i_j}\|^2]^\frac12
\leq &  \E_{\mathcal{F}_j}[\E[\|v_{H,j,i_j}\varphi_{i_j}\|^2|\mathcal{F}_j]]^\frac12 = n^\frac12\E[\|v_{H,j}\|^2]^\frac12,
\end{align*}
and we may bound $\E[\|N_{j,2}\|^2]^\frac12$ by 
\begin{align*}
\E[\|N_{j,2}\|^2]^\frac12  \leq& \E[\|\E[v_{F,j}]-v_{F,j,i_j}\varphi_{i_j}\|^2]^\frac12+ \lambda_j^\delta \E[\|R^*(\E[v_{G,j}]-v_{G,k,i_j}\varphi_{i_j})\|^2]^\frac12\\
\leq& n^\frac12\E[\|v_{F,j}\|^2]^\frac12+c_R\lambda_j^\delta n^\frac12\E[\|v_{G,j}\|^2]^\frac12.
\end{align*}
Now, under Assumptions \ref{ass:sol}(i)(ii) and Assumption \ref{ass:stoch}, we estimate $\E[\|v_{F,j}\|^2]^\frac12$ and $\E[\|v_{G,j}\|^2]^\frac12$ one by one. For $v_{F,j}$ defined in \eqref{eqn:w_F}, by the triangle inequality, Assumptions \ref{ass:stoch} and Lemma \ref{lem:R}, there holds
\begin{align*}
\E[\|v_{F,j}\|^2]^\frac12\leq& \E[\|( R_{F,x_j^\delta}^{ *}-I )(F(x_j^\delta)-y^\delta)\|^2]^\frac12+ \E[\|\big(F(x_j^\delta)-F(x^\dag)-K_{F}(x_j^\delta-x^\dag)\big)\|^2]^\frac12+ \E[\|\xi\|^2]^\frac12\\
\leq& c_F\E[\|e_j^\delta\|^2]^\frac{\theta}{2} \E[\|F(x_j^\delta)-y^\delta\|^2]^\frac12+\frac{c_F}{1+\theta}\E[\|K_F e_j^\delta\|^2]^\frac12\E[\|e_j^\delta\|^2]^\frac{\theta}{2}+\delta.
\end{align*}
Further, by the triangle inequality and Lemma \ref{lem:linear}, there holds
\begin{align*}
\E[\|F(x_j^\delta)-y^\delta\|^2]^\frac12 \leq &\E[\|F(x_j^\delta)-F(x^\dag)\|^2]^\frac12+\E[\|\xi\|^2]^\frac12
\leq \frac{1}{1-\eta_F}\E[\|K_F e_j^\delta\|^2]^\frac12 +\delta,
\end{align*}
which implies that
\begin{align*}
\E[\|v_{F,j}\|^2]^\frac12
&\leq c_F\E[\|e_j^\delta\|^2]^\frac{\theta}{2} (\frac{1}{1-\eta_F}\E[\|K_F e_j^\delta\|^2]^\frac12 +\delta)+\frac{c_F}{1+\theta}\E[\|K_F e_j^\delta\|^2]^\frac12\E[\|e_j^\delta\|^2]^\frac{\theta}{2}+\delta\\
&\leq \frac{c_F(2+\theta-\eta_F)}{(1+\theta)(1-\eta_F)}\E[\|e_j^\delta\|^2]^\frac{\theta}{2}\E[\|K_F e_j^\delta\|^2]^\frac12 +(c_F\E[\|e_j^\delta\|^2]^\frac{\theta}{2}+1)\delta.
\end{align*}
Similarly, for $v_{G,j}$ defined in \eqref{eqn:w_G}, by Assumptions \ref{ass:sol}(iii) and \ref{ass:stoch}, and Lemma \ref{lem:R}, we obtian that
\begin{align*}
\E[\|v_{G,j}\|^2]^\frac12
\leq&\E[\|( R_{G,x_j^\delta}^{*}-I )(G(x_j^\delta)-y^\delta)\|^2]^\frac12+\E[\|G(x_j^\delta)-G(x^\dag)-K_{G}(x_j^\delta-x^\dag)\|^2]^\frac12 +C_{max}+\delta\\
\leq&c_G\E[\|e_j^\delta\|^2]^\frac{\theta}{2} \E[\|G(x_j^\delta)-y^\delta\|^2]^\frac12+\frac{c_G}{1+\theta}\E[\|K_G e_j^\delta\|^2]^\frac12\E[\|e_j^\delta\|^2]^\frac{\theta}{2}+C_{max}+\delta.
\end{align*}
{Further, by Lemma \ref{lem:resG_refined}, we have
\begin{align*}
\E[\|G(x_j^\delta)-y^\delta\|^2]^\frac12\leq &(c_G\E[\|e_j^\delta\|^2]^\frac12+1)\E[ \|K_{G}e_j^\delta\|^2]^\frac12+C_{max}+\delta
\end{align*}}
and thus
\begin{align*}
\E[\|v_{G,j}\|^2]^\frac12
\leq&{c_G(c_G\E[\|e_j^\delta\|^2]^\frac12+1+\tfrac{1}{1+\theta})}\E[\|e_j^\delta\|^2]^\frac{\theta}{2} \E[\|K_G e_j^\delta\|^2]^\frac12 +(c_G\E[\|e_j^\delta\|^2]^\frac{\theta}{2}+1)(C_{max}+\delta).
\end{align*}
Combining these two estimates gives the bound on $\E[\|N_{j,2}\|^2]^\frac12$ that
\begin{align*}
\E[\|N_{j,2}\|^2]^\frac12 
\leq& n^\frac12\E[\|v_{F,j}\|^2]^\frac12+c_R\lambda_j^\delta n^\frac12\E[\|v_{G,j}\|^2]^\frac12\\
\leq& n^\frac12 \frac{c_F(2+\theta-\eta_F)}{(1+\theta)(1-\eta_F)}\E[\|e_j^\delta\|^2]^\frac{\theta}{2}\E[\|K_F e_j^\delta\|^2]^\frac12+n^\frac12\big((c_F+c_R\lambda_j^\delta c_G)\E[\|e_j^\delta\|^2]^\frac{\theta}{2}+c_R\lambda_j^\delta+1\big)\delta\\ &+n^\frac12c_R\lambda_j^\delta {c_G(c_G\E[\|e_j^\delta\|^2]^\frac12+1+\tfrac{1}{1+\theta})}\E[\|e_j^\delta\|^2]^\frac{\theta}{2}\E[ \|K_{G}e_j^\delta\|^2]^\frac12+n^\frac12c_R\lambda_j^\delta (c_G\E[\|e_j^\delta\|^2]^\frac{\theta}{2}+1)C_{max}.
\end{align*}
Now, with Assumptions \ref{ass:sol}(iv)(v), we simplify this estimate as
\begin{align*}
\E[\|N_{j,2}\|^2]^\frac12 
\leq& n^\frac12 (\frac{c_F(2+\theta-\eta_F)}{(1+\theta)(1-\eta_F)}+c_R^2\lambda_j^\delta {c_G(c_G\E[\|e_j^\delta\|^2]^\frac12+1+\tfrac{1}{1+\theta})})\E[\|e_j^\delta\|^2]^\frac{\theta}{2}\E[ \|K_{F}e_j^\delta\|^2]^\frac12 \\
&+n^\frac12c_R\lambda_j^\delta (c_G\E[\|e_j^\delta\|^2]^\frac{\theta}{2}+1)C_{max}+n^\frac12\big((c_F+c_R\lambda_j^\delta c_G)\E[\|e_j^\delta\|^2]^\frac{\theta}{2}+c_R\lambda_j^\delta+1\big)\delta.
\end{align*}
The notation $B_F^\frac12=K_F$ completes the proof of the lemma.

\subsection{Proof of Lemma \ref{lem:constants}}\label{app:lem:constants}
By the definitions of $C_j, C_j^G, C_j^F, \tilde{C}_j, \tilde{C}_j^G$ and $\tilde{C}_j^F$ and the assumption $\lambda_j^\delta\leq\lambda_0^\delta\leq \min(c_R^{-2},c_R^{-1})$, for any $\theta\in(0,1]$ {and $\eta_F\in[0,1)$}, we derive the estimates
\begin{align*}
C_j=&\frac{3-\eta_F}{2(1-\eta_F)}c_F + {(c_G\E[\|e_j^\delta\|^2]^\frac12+\tfrac{3}{2})}c_Gc_R^2\lambda_j^\delta\leq \frac{2+\theta-\eta_F}{(1+\theta)(1-\eta_F)}c_F+ {(c_G\E[\|e_j^\delta\|^2]^\frac12+1+\frac{1}{1+\theta})}c_Gc_R^2\lambda_j^\delta=\tilde{C}_j\\
\leq& (\frac{1}{1+\theta}+\frac{1}{1-\eta_F})c_F+ {(c_G\E[\|e_j^\delta\|^2]^\frac12+1+\frac{1}{1+\theta})}c_G
\leq {(1+\frac{1}{1-\eta_F})c_F+ (2+c_G\E[\|e_j^\delta\|^2]^\frac12)c_G},\\
\max(&C^G_j,\tilde{C}^G_j)\leq\max(c_R(c_G\E[\|e_j^\delta\|^2]^\frac12+1),c_R(c_G \E[\|e_j^\delta\|^2]^\frac\theta2+1))\leq c_R (c_G\max(\E[\|e_j^\delta\|^2]^\frac12,1)+1)\\
&\qquad\quad\;\;\;\leq c_R (c_G(\E[\|e_j^\delta\|^2]^\frac12+1)+1),\\
\max(&C^F_j,\tilde{C}^F_j)\leq\max((c_{F}+c_G)\E[\|e_j^\delta\|^2]^\frac12+2, (c_{F}+c_G)\E[\|e_j^\delta\|^2]^\frac\theta2+2)\leq(c_{F}+c_G)\max(\E[\|e_j^\delta\|^2]^\frac12, 1)+2\\
&\qquad\quad\;\;\;\leq (c_{F}+c_G)(\E[\|e_j^\delta\|^2]^\frac12+1)+2.
\end{align*}
This completes the proof.
\end{proof}

\subsection{Estimates for Section \ref{ssec:rate}}
Now, we give a set of estimates employed in the analysis of convergence rate in Section \ref{ssec:rate}.
The next lemma gives a variant of a well known estimate on operator norms (see, e.g., \cite[Lemma 15]{LinRosasco:2017}).
\begin{lemma}\label{lem:estimate-B}
Under Assumptions \ref{ass:sol}(v) and \ref{ass:stepsize}, for any $j<k$ and $s\geq 0$, there holds
\begin{equation*}
  \phi_j^{s} = \|B_F^{s}\Pi_{j+1}^k(B)\|=\|B_F^{s}\prod_{i=j+1}^k\big(I-\eta_i(B_F+\lambda_i^\delta B_G)\big)\| 
\leq (se^{-1}(\sum_{i=j+1}^k\eta_i)^{-1})^s.
\end{equation*}
\end{lemma}
\begin{proof}
With the definitions $B_F=K_F^*K_F$ and $B_G=K_G^*K_G$, and the singular value decomposition of the operators $K_F$ and $K_G$ in Assumption \ref{ass:sol}(v), we have 
\begin{align*}
\phi_j^{s} = \|B_F^{s}\Pi_{j+1}^k(B)\|=\|B_F^{s}\prod_{i=j+1}^k\big(I-\eta_i(B_F+\lambda_i^\delta B_G)\big)\| =\sup_{t\geq1}\sigma_t^{2s}\prod_{i=j+1}^k\big(1-\eta_i(\sigma_t^2+\lambda_i^\delta \tilde\sigma_t^2)\big).
\end{align*}
Further, by using the fact $1-x\leq e^{-x}$ for any $x\in[0,1]$ and {Assumption \ref{ass:stepsize}(ii)} which implies $\eta_i(\sigma_t^2+\lambda_i^\delta \tilde\sigma_t^2)\in[0,1]$ for any $i, t\geq1$, we can derive that 
\begin{align*}
\phi_j^{s} =& \sup_{t\geq1}\sigma_t^{2s}\prod_{i=j+1}^k\big(1-\eta_i(\sigma_t^2+\lambda_i^\delta \tilde\sigma_t^2)\big)\leq \sup_{t\geq1}\sigma_t^{2s}\prod_{i=j+1}^k e^{-\eta_i(\sigma_t^2+\lambda_i^\delta \tilde\sigma_t^2)}=\sup_{t\geq1}\sigma_t^{2s}e^{-\sum_{i=j+1}^k(\eta_i(\sigma_t^2+\lambda_i^\delta \tilde\sigma_t^2))}\\
=&\sup_{t\geq1}\sigma_t^{2s}e^{-(\sum_{i=j+1}^k\eta_i)\sigma_t^2}e^{-(\sum_{i=j+1}^k\eta_i\lambda_i^\delta) \tilde\sigma_t^2}\leq \sup_{t\geq1}\sigma_t^{2s}e^{-(\sum_{i=j+1}^k\eta_i)\sigma_t^2}.
\end{align*}
For the function $g(x)=x^{s} e^{-ax}$, with some constant $a>0$,
the maximum is attained at $x=sa^{-1}$, with a maximum value $s^s(ea)^{-s}$.
Then setting $a=\sum_{i=j+1}^k\eta_i$ 
complete the proof of the lemma.
\end{proof}
Next we gather several useful estimates from \cite{JinZhouZou:2020} in Lemma \ref{lem:basicest}.
\begin{lemma}\label{lem:basicest}
If $\eta_j=\eta_0j^{-\alpha}$, $\alpha\in(0,1)$, $\beta\in[0,1]$ and $r\geq 0$, then there hold
\begin{align}
&\;\;\;\sum_{i=1}^k \eta_i  \geq (1-2^{\alpha-1})(1-\alpha)^{-1}\eta_0(k+1)^{1-\alpha},\label{eqn:bdd-phi}\\
&\;\;\;\sum_{j=1}^{k-1} \frac{\eta_j}{(\sum_{\ell=j+1}^k\eta_\ell)^{r}} j^{-\beta} \leq
\eta_0^{1-r} B(1-r,1-\alpha-\beta) k^{(1-r)(1-\alpha)-\beta}, \quad r\in[0,1),  \alpha+\beta<1,\label{eqn:bdd-phi-r}\\
&\;\;\;\sum_{j=1}^{[\frac k2]}\frac{\eta_j^2}{(\sum_{\ell=j+1}^k\eta_\ell)^r}j^{-\beta} \leq c_{\alpha,\beta,r} k^{-r(1-\alpha)+\max(0,1-2\alpha-\beta)},\label{eqn:bdd-phi-r-2}\\
&\sum_{j=[\frac k2]+1}^{k-1}\frac{\eta_j^2}{(\sum_{\ell=j+1}^k\eta_\ell)^r}j^{-\beta}\leq c'_{\alpha,\beta,r} k^{-((2-r)\alpha+\beta)+\max(0,1-r)},\label{eqn:bdd-phi-r-22}
\end{align}
where we slightly abuse $k^{-\max(0,0)}$ for $\ln k$,  $B(\cdot,\cdot)$ denotes the Beta function defined by
\begin{equation}\label{eqn:Beta}
B(a,b)=\int_0^1s^{a-1}(1-s)^{b-1}{\rm d} s \;\;\mbox{ for any }a,b>0,
\end{equation} 
and the constants $c_{\alpha,\beta,r}$ and $c'_{\alpha,\beta,r}$ are given by
\begin{align*}
c_{\alpha,\beta,r} =2^r\eta_0^{2-r}\left\{\begin{array}{ll}
1+(2\alpha+\beta-1)^{-1},   & 2\alpha+\beta>1,\\
2,                  & 2\alpha+\beta=1,\\
2^{2\alpha+\beta-1}(1-2\alpha-\beta)^{-1},     & 2\alpha+\beta<1,
\end{array}\right.\quad\mbox{and}\quad
   c'_{\alpha,\beta,r}  = 2^{2\alpha+\beta}\eta_0^{2-r}\left\{\begin{array}{ll}
     1+(r-1)^{-1}, & r>1,\\
     2, & r= 1,\\
   2^{r-1}(1-r)^{-1}, & r<1.
   \end{array}\right.
\end{align*}
\end{lemma}

The next result collects some lengthy estimates, following routine rather tedious computations, which are essential for the proof of Theorems \ref{thm:err-total-ex} and \ref{thm:err-total}. 
\begin{prop}\label{prop:est-ex}
Under the conditions in Theorems \ref{thm:err-total-ex} and \ref{thm:err-total}, especially the conditions $\|B_F\|\leq1$ and $ \eta_0\leq 1$,
the following estimates hold for any $\theta\in(0,\frac{1-\alpha}{\beta}-1)$ and $\epsilon\in(0,2\theta\beta)$, with $$\beta=\min(2\nu(1-\alpha),\alpha), \quad\gamma=\min((1+2\nu)(1-\alpha),1)\quad \mbox{and} \quad \zeta=1-\alpha-\frac{\gamma}{2}:$$
\begin{align}
\sum_{j=1}^k\eta_j\phi_j^\frac12 j^{-\frac{\gamma}{2}}&\leq 2^{\frac{\beta}{2}-1}\eta_0^{\frac12}(B(\frac12,\zeta)+2)(k+1)^{-\frac{\beta}{2}},\label{eqn:est1-ex}\\  \sum_{j=1}^k\eta_j^2(\phi_j^{\frac12})^2 j^{-\gamma}  &  \leq 2^{\beta-1}\eta_0((\alpha+\beta)^{-1}+4)(k+1)^{-\beta},\label{eqn:est2-ex}\\
\sum_{j=1}^k\eta_j\phi^1_j j^{-\frac{\gamma}{2}} 
&\leq  2^{\frac{\gamma}{2}-1} \eta_0^{\frac{\epsilon}{4(1-\alpha)}}(B(\tfrac{\epsilon}{4(1-\alpha)},\zeta)+2)(k+1)^{\frac{\epsilon}{4}-\frac{\gamma}{2}},\label{eqn:est3.1-ex}\\
\sum_{j=1}^k\eta_j\phi^1_j j^{-\frac{\gamma+\theta\beta}{2}}
&\leq 2^{\frac{\gamma}{2}-\frac12} \eta_0^{\tfrac{2\theta \beta-\epsilon}{4(1-\alpha)}}\big(B(\tfrac{2\theta \beta-\epsilon}{4(1-\alpha)},\zeta-\tfrac{\theta\beta}{2})+2\big)(k+1)^{-\frac{\epsilon}{4}-\frac{\gamma}{2}},\label{eqn:est3.2-ex}\\    
\sum_{j=1}^{k} \eta_j^2(\phi_j^1)^2j^{-\gamma}
&\leq 2^{\gamma+1}\eta_0^{1-\frac{\beta}{1-\alpha}}(
\alpha^{-1} +1)(k+1)^{-\gamma},\label{eqn:est4-ex}\\
\sum_{j=1}^{k}\eta_j \phi_j^\frac12 
&\leq 2^{-1}\eta_0^\frac12( B(\tfrac12,1-\alpha)+2) (k+1)^{\frac{1-\alpha}{2}},\label{eqn:est5-ex}\\
   \sum_{j=1}^{k}\eta_j\phi_j^1  
& \leq 2^{-1} \eta_0^{\frac{\epsilon}{4(1-\alpha)}}(B(\tfrac{\epsilon}{4(1-\alpha)},1-\alpha)+2)(k+1)^{\frac{\epsilon}{4}}.\label{eqn:est6-ex}
\end{align}
\end{prop}
\begin{proof}
A similar analysis can be found in \cite{JinZhouZou:2020}. We refined the analysis in order to derive the recursion of upper bounds on $a_{k+1}$ and $b_{k+1}$ in Theorems \ref{thm:err-total-ex} and \ref{thm:err-total}, where the estimates \eqref{eqn:est1-ex}, \eqref{eqn:est2-ex} and \eqref{eqn:est5-ex} are needed for $a_{k+1}$ while the others are for $b_{k+1}$. Now, we show the estimates one by one.
First, by Lemma \ref{lem:estimate-B}, \eqref{eqn:bdd-phi-r}, and the conditions $\|B_F\|\leq 1$ and $\eta_0\leq1$, we derive that
\begin{align*}
\sum_{j=1}^k\eta_j\phi_j^\frac12j^{-\frac{\gamma}{2}}
\leq& \sum_{j=1}^{k-1}\eta_j((2e)^{-\frac12}(\sum_{i=j+1}^k\eta_i)^{-\frac12})j^{-\frac{\gamma}{2}}+\eta_0\|B_F^\frac12\|k^{-\alpha-\frac{\gamma}{2}}
\leq (2e)^{-\frac12}\sum_{j=1}^{k-1}\frac{\eta_j}{(\sum_{i=j+1}^k\eta_i)^{\frac12}}j^{-\frac{\gamma}{2}}+\eta_0^{\frac12}k^{-\alpha-\frac{\gamma}{2}}\\
\leq& (2e)^{-\frac12}\eta_0^{\frac12} B(\frac12,1-\alpha-\frac{\gamma}{2}) k^{-\frac12(1-\alpha)+1-\alpha-\frac{\gamma}{2}}+\eta_0^{\frac12}k^{-\alpha-\frac{\gamma}{2}}.
\end{align*}
Using the relations $\beta=\gamma-(1-\alpha)$ and $\zeta=1-\alpha-\frac{\gamma}{2}$ follow directly from the definitions of $\beta$, $\gamma$ and $\zeta$, we further simplify the above bound by
\begin{align*}
\sum_{j=1}^k\eta_j\phi_j^\frac12j^{-\frac{\gamma}{2}}
\leq& (2e)^{-\frac12}\eta_0^{\frac12} B(\frac12,\zeta) k^{-\frac\beta2}+\eta_0^{\frac12}k^{-\frac{\beta}{2}}
\leq \eta_0^{\frac12}((2e)^{-\frac12} B(\frac12,\zeta)+1)k^{-\frac{\beta}{2}}
\leq 2^{-1}\eta_0^{\frac12}(B(\frac12,\zeta)+2)k^{-\frac{\beta}{2}}.
\end{align*}
Then the inequality $2k\geq k+1$ for $k\geq 1$ immediately implies the estimate \eqref{eqn:est1-ex}.
Similarly, it follows
from Lemma \ref{lem:estimate-B}, \eqref{eqn:bdd-phi-r-2} and \eqref{eqn:bdd-phi-r-22}, that
\begin{align*}
\sum_{j=1}^k\eta_j^2(\phi_j^{\frac12})^2 j^{-\gamma}  
\leq& (2e)^{-1}\sum_{j=1}^{k-1}\frac{\eta_j^2}{\sum_{i=j+1}^k\eta_i}j^{-\gamma}+\eta_0^2k^{-2\alpha-\gamma}\\
\leq& (2e)^{-1}(c_{\alpha,\gamma,1}\eta_0 k^{-(1-\alpha)+\max(0,1-2\alpha-\gamma)}+c'_{\alpha,\gamma,1}\eta_0 k^{-(\alpha+\gamma)+\max(0,0)})+\eta_0k^{-\gamma}.
\end{align*}
By the facts that $1-2\alpha-\gamma=-(\gamma-(1-\alpha)+\alpha)=-(\beta+\alpha)<0$ and $\alpha+\gamma=\alpha+\beta+(1-\alpha)=\beta+1$, there holds
\begin{align*}
\sum_{j=1}^k\eta_j^2(\phi_j^{\frac12})^2 j^{-\gamma} 
\leq& (2e)^{-1}(c_{\alpha,\gamma,1}\eta_0 k^{-(1-\alpha)}+c'_{\alpha,\gamma,1}\eta_0 k^{-(\beta+1)}\ln k)+\eta_0k^{-\gamma}\\
\leq& 2^\beta\eta_0\big((2e)^{-1}(c_{\alpha,\gamma,1} +c'_{\alpha,\gamma,1} k^{-1}\ln k)+1\big)(k+1)^{-\beta}.
\end{align*}
Then, by using the inequality, for any $s>0$ and $k\geq 1$,
\begin{equation}\label{eqn:bdd-log}
k^{-s}\ln k \leq (es)^{-1},
\end{equation}
and setting $s=1$, we derive that 
\begin{align*}
\sum_{j=1}^k\eta_j^2(\phi_j^{\frac12})^2 j^{-\gamma} 
\leq& 2^\beta\eta_0\big((2e)^{-1}(c_{\alpha,\gamma,1} +e^{-1}c'_{\alpha,\gamma,1})+1\big)(k+1)^{-\beta}
\end{align*}
Further, the definition of the constants $c_{\alpha,\gamma,1}$ and $c'_{\alpha,\gamma,1}$ in Lemma \ref{lem:basicest}, with the inequalities 
\begin{equation}\label{eqn:alpha-gamma}
1<2\alpha+\gamma\leq 2\alpha+(1+2\nu)(1-\alpha)=2-(1-2\nu)(1-\alpha)< 2,    
\end{equation} 
and the relation $2\alpha+\gamma-1=\alpha+\beta$,
implies the estimate \eqref{eqn:est2-ex}
\begin{align*}
\sum_{j=1}^k\eta_j^2(\phi_j^{\frac12})^2 j^{-\gamma} 
\leq& 2^\beta\eta_0\big((2e)^{-1}(2(1+(2\alpha+\gamma-1)^{-1}) +e^{-1}2^{1+2\alpha+\gamma})+1\big)(k+1)^{-\beta}\\
\leq& 2^{\beta-1}\eta_0(2e^{-1}+(2\alpha+\gamma-1)^{-1} +2^3e^{-2}+2)(k+1)^{-\beta}
\leq 2^{\beta-1}\eta_0((\alpha+\beta)^{-1}+4)(k+1)^{-\beta}.
\end{align*}
By noting the inequality that $\phi_j^1\leq\|B_F^{1-r}\|\phi_j^r\leq\phi_j^r$ for any $r\in[\frac12,1)$, with \eqref{eqn:bdd-phi-r} and the fact that, for any $\theta\in(0,\frac{1-\alpha}{\beta}-1)$ and $\epsilon\in(0,2\theta\beta)$, we derive that
\begin{align*}
\sum_{j=1}^k\eta_j\phi^1_j j^{-\frac{\gamma+\theta\beta}{2}} 
\leq& \sum_{j=1}^k\eta_j\phi^r_j j^{-\frac{\gamma+\theta\beta}{2}} 
\leq (\frac{r}{e})^r \eta_0^{1-r}B(1-r,1-\alpha-\tfrac{\gamma+\theta\beta}{2})k^{(1-r)(1-\alpha)-\frac{\gamma+\theta\beta}{2}}+\eta_0k^{-\alpha-\frac{\gamma+\theta\beta}{2}}\\
\leq& (\frac{r}{e})^r \eta_0^{1-r}B(1-r,\zeta-\tfrac{\theta\beta}{2})k^{(1-r)(1-\alpha)-\frac{\gamma+\theta\beta}{2}}+\eta_0k^{-\alpha-\frac{\gamma+\theta\beta}{2}}.
\end{align*}
Then setting $r=1-\tfrac{2\theta \beta-\epsilon}{4(1-\alpha)}
\in(\frac12,1)$, with the inequality $\eta_0\leq 1$ and the fact that the function $(\frac{r}{e})^r$ is decreasing in $r$ over the interval $[\frac12,1]$, gives
\begin{align*}
\sum_{j=1}^k\eta_j\phi^1_j j^{-\frac{\gamma+\theta\beta}{2}} 
\leq& 2^{-1} \eta_0^{\tfrac{2\theta \beta-\epsilon}{4(1-\alpha)}}B(\tfrac{2\theta \beta-\epsilon}{4(1-\alpha)},\zeta-\tfrac{\theta\beta}{2})k^{-\frac{\epsilon}{4}-\frac{\gamma}{2}}+\eta_0k^{-\frac{\epsilon}{4}-\frac{\gamma}{2}}
\leq 2^{-1} \eta_0^{\tfrac{2\theta \beta-\epsilon}{4(1-\alpha)}}\big(B(\tfrac{2\theta \beta-\epsilon}{4(1-\alpha)},\zeta-\tfrac{\theta\beta}{2})+2\big)k^{-\frac{\epsilon}{4}-\frac{\gamma}{2}}\\
\leq& 2^{\frac{\gamma}{2}+\frac{\epsilon}{4}-1} \eta_0^{\tfrac{2\theta \beta-\epsilon}{4(1-\alpha)}}\big(B(\tfrac{2\theta \beta-\epsilon}{4(1-\alpha)},\zeta-\tfrac{\theta\beta}{2})+2\big)(k+1)^{-\frac{\epsilon}{4}-\frac{\gamma}{2}}.
\end{align*}
Then, the fact that $\epsilon< 2\theta\beta<2(1-\alpha-\beta)<2$ yields the estimate \eqref{eqn:est3.2-ex}.
Similarly, for any $\epsilon\in(0,2\theta\beta)$, we have $1-\frac{\epsilon}{4(1-\alpha)}\in(\frac12,1)$ and
\begin{align*}
\sum_{j=1}^k\eta_j\phi^1_j j^{-\frac{\gamma}{2}} 
\leq& \sum_{j=1}^k\eta_j\phi^{1-\frac{\epsilon}{4(1-\alpha)}}_j j^{-\frac{\gamma}{2}} 
\leq (\frac{1-\frac{\epsilon}{4(1-\alpha)}}{e})^{1-\frac{\epsilon}{4(1-\alpha)}}\sum_{j=1}^{k-1}\frac{\eta_j}{(\sum_{i=j+1}^k\eta_i)^{1-\frac{\epsilon}{4(1-\alpha)}}}j^{-\frac{\gamma}{2}}+\eta_0k^{-\alpha-\frac{\gamma}{2}}\\
\leq& (2e)^{-\frac12}  \eta_0^{\frac{\epsilon}{4(1-\alpha)}}B(\tfrac{\epsilon}{4(1-\alpha)},1-\alpha-\tfrac{\gamma}{2})k^{\frac{\epsilon}{4}-\frac{\gamma}{2}}+\eta_0k^{-\alpha-\frac{\gamma}{2}}\\
\leq& 2^{-1} \eta_0^{\frac{\epsilon}{4(1-\alpha)}}(B(\tfrac{\epsilon}{4(1-\alpha)},\zeta)+2)k^{\frac{\epsilon}{4}-\frac{\gamma}{2}}
\leq 2^{\frac{\gamma}{2}-1} \eta_0^{\frac{\epsilon}{4(1-\alpha)}}(B(\tfrac{\epsilon}{4(1-\alpha)},\zeta)+2)(k+1)^{\frac{\epsilon}{4}-\frac{\gamma}{2}}.
\end{align*}
Now, we bound $\sum_{j=1}^{k} \eta_j^2(\phi_j^1)^2j^{-\gamma}$ by decomposing it into three parts
\begin{align*}
\sum_{j=1}^{k} \eta_j^2(\phi_j^1)^2j^{-\gamma}\leq &\sum_{j=1}^{[\frac k2]} \eta_j^2(\phi_j^r)^2j^{-\gamma} + \sum_{j=[\frac k2]+1}^{k-1}\eta_j^2(\phi_j^\frac12)^2j^{-\gamma}+\eta_0^2 k^{-2\alpha-\gamma}.
\end{align*}
Then by \eqref{eqn:bdd-phi-r-2}, \eqref{eqn:bdd-phi-r-22}, \eqref{eqn:alpha-gamma} and the equation $2\alpha+\gamma-1=\alpha+\beta$, we obtain that
\begin{align*}
\sum_{j=1}^{k} \eta_j^2(\phi_j^1)^2j^{-\gamma}\leq &(\frac{r}{e})^{2r} c_{\alpha,\gamma,2r}\eta_0 k^{-2r(1-\alpha)+\max(0,1-2\alpha-\gamma)} + (2e)^{-1} c'_{\alpha,\gamma,1}\eta_0 k^{-(\alpha+\gamma)+\max(0,0)}+\eta_0^2 k^{-\gamma}\\
\leq &\big((\frac{2r}{e})^{2r}\eta_0^{2-2r}
(1+(2\alpha+\gamma-1)^{-1}) k^{\gamma-2r(1-\alpha)} + (2e)^{-1} 2^{2\alpha+\gamma+1}\eta_0 k^{-\alpha}\ln k+\eta_0^2 \big)k^{-\gamma}\\
\leq &2^{\gamma}\big(\eta_0^{2-2r}
(1+(\alpha+\beta)^{-1}) k^{\gamma-2r(1-\alpha)} + 2\eta_0 k^{-\alpha}\ln k+\eta_0^2 \big)(k+1)^{-\gamma}
\end{align*}
By setting $r=\frac{\gamma}{2(1-\alpha)}=\frac{1}{2}+\frac{\beta}{2(1-\alpha)}\in(\frac12,1)$, together with the inequalities $\eta_0\leq 1$ and \eqref{eqn:bdd-log} with $s=\alpha$, the estimate \eqref{eqn:est4-ex} holds
\begin{align*}
\sum_{j=1}^{k} \eta_j^2(\phi_j^1)^2j^{-\gamma}
\leq &2^{\gamma}\big(\eta_0^{2-2r}
(1+(\alpha+\beta)^{-1}) + 2\eta_0(e\alpha)^{-1} +\eta_0^2 \big)(k+1)^{-\gamma}\\
\leq &2^{\gamma}\eta_0^{1-\frac{\beta}{1-\alpha}}(
(\alpha+\beta)^{-1} + \alpha^{-1} +2)(k+1)^{-\gamma}
\leq 2^{\gamma+1}\eta_0^{1-\frac{\beta}{1-\alpha}}(\alpha^{-1} +1)(k+1)^{-\gamma}.
\end{align*}
Finally, the following estimates \eqref{eqn:est5-ex} and \eqref{eqn:est6-ex}, for any $\epsilon\in(0,2\theta\beta)$, that
\begin{align*}
\sum_{j=1}^{k}\eta_j \phi_j^\frac12 
\leq& (2e)^{-\frac12}\sum_{j=1}^{k-1}\frac{\eta_j}{(\sum_{\ell=1}^k\eta_\ell)^\frac12} + \eta_0k^{-\alpha}
\leq (2^{-1}\eta_0^\frac12 B(\tfrac12,1-\alpha)+\eta_0) k^{\frac{1-\alpha}{2}}\\
\leq& 2^{-1}\eta_0^\frac12( B(\tfrac12,1-\alpha)+2) (k+1)^{\frac{1-\alpha}{2}},\\
\sum_{j=1}^{k}\eta_j\phi_j^1  
\leq & \sum_{j=1}^{k}\eta_j\phi^{1-\frac{\epsilon}{4(1-\alpha)}}_j  
\leq (\frac{1-\frac{\epsilon}{4(1-\alpha)}}{e})^{1-\frac{\epsilon}{4(1-\alpha)}}\sum_{j=1}^{k-1}\frac{\eta_j}{(\sum_{i=j+1}^k\eta_i)^{1-\frac{\epsilon}{4(1-\alpha)}}}+\eta_0k^{-\alpha}\\
\leq& (2e)^{-\frac12}  \eta_0^{\frac{\epsilon}{4(1-\alpha)}}B(\tfrac{\epsilon}{4(1-\alpha)},1-\alpha)k^{\frac{\epsilon}{4}}+\eta_0k^{-\alpha}
\leq 2^{-1} \eta_0^{\frac{\epsilon}{4(1-\alpha)}}(B(\tfrac{\epsilon}{4(1-\alpha)},1-\alpha)+2)(k+1)^{\frac{\epsilon}{4}},
\end{align*}
complete the proof.
\end{proof}

\bibliographystyle{abbrv}
\bibliography{sgd}
\end{document}